\documentclass[11pt,a4paper]{article}
\usepackage{amsmath,amssymb,amsthm,enumerate,graphicx,mathtools,titling}
\usepackage[margin=2.5cm]{geometry}

\usepackage[usenames, dvipsnames]{color}

\usepackage{xypic} 
\usepackage{mathrsfs} 
\usepackage{xspace}
\usepackage{authblk}

\usepackage[backref=page]{hyperref}

\newtheorem{thm}{Theorem}[section]

\newtheorem{prop}[thm]{Proposition}
\newtheorem{lem}[thm]{Lemma} 
\newtheorem{cor}[thm]{Corollary}

\theoremstyle{definition}
\newtheorem{defn}[thm]{Definition}
\newtheorem{quest}[thm]{Question}
\newtheorem{rmk}[thm]{Remark}

\newtheorem*{claim*}{Claim}

\newcommand{\mc}[1]{\mathcal{#1}}

\newcommand{\mf}[1]{\mathfrak{#1}}
\newcommand{\ms}[1]{\mathscr{#1}} 

\newcommand{\Zb}{\mathbb{Z}}
\newcommand{\Nb}{\mathbb{N}}

\DeclareMathOperator{\sph}{sph}

\newcommand{\Univ}{\mathbf{U}}

\newcommand{\Reg}{\widetilde{\mc{O}}_c}
\newcommand{\Oc}{\mc{O}_c}
\newcommand{\Der}{\mathrm{D}}

\newcommand{\Mon}{\mathrm{Mon}}
\newcommand{\NMon}{\mathrm{Mon}_{\mc{RE}}}

\newcommand{\Sclass}{\ms{S}}

\newcommand{\Bigsclass}{[\mathrm{Sim}]}

\newcommand{\Es}[1]{\ms{E}^{#1}}

\newcommand{\xip}[1]{\xi^{#1}}

\newcommand{\tdlc}{t.d.l.c.\@\xspace}
\newcommand{\tdlcsc}{t.d.l.c.s.c.\@\xspace}

\newcommand{\hji}{h.j.i.\@\xspace}
\newcommand{\defbold}{\textbf}

\newcommand{\rist}{\mathrm{rist}}

\newcommand{\inv}{^{-1}}
\newcommand{\triv}{\{1\}}

\newcommand{\Homeo}{\mathrm{Homeo}}

\newcommand{\QZ}{\mathrm{QZ}}

\newcommand{\CC}{\mathrm{C}}
\newcommand{\N}{\mathrm{N}}
\newcommand{\Z}{\mathrm{Z}}

\newcommand{\con}{\mathrm{con}}

\newcommand{\cgrp}[1]{\overline{\langle #1 \rangle}}

\newcommand{\grp}[1]{\langle #1 \rangle}
\newcommand{\ol}[1]{\overline{#1}}

\newcommand{\Sym}{\mathop{\rm Sym}\nolimits}
\newcommand{\Aut}{\mathop{\rm Aut}\nolimits}

\newcommand{\Rad}[1]{\mathop{\rm Rad}_{#1}\nolimits}
\newcommand{\Res}[1]{\mathop{\rm Res}_{#1}\nolimits}

\newcommand{\propS}{(\mathrm{S})}
\newcommand{\propP}[1]{(\mathrm{P}_{#1})}

\newcommand{\RadRE}{\Rad{\mc{RE}}}

\makeindex
\begin{document}

\title{A class of well-founded totally disconnected locally compact groups}

\preauthor{\large}
\DeclareRobustCommand{\authoring}{
\renewcommand{\thefootnote}{\arabic{footnote}}
\begin{center}Colin D. Reid\textsuperscript{1}\footnotetext[1]{Research supported by ARC grant FL170100032.}
\\ \bigskip
The University of Newcastle, School of Mathematical and Physical Sciences, Callaghan, NSW 2308, Australia.\\
\href{mailto:colin@reidit.net}{colin@reidit.net}
\end{center}
}
\author{\authoring}
\postauthor{\par}

\maketitle

\begin{abstract}
Motivated by the problem of finding a ``well-foundedness principle'' for totally disconnected, locally compact (\tdlc) groups, we introduce a class $\Es{\Sclass}$ of \tdlc groups, containing P. Wesolek's class $\Es{}$ of (regionally) elementary groups but also including many groups in the class $\Sclass$ of nondiscrete compactly generated topologically simple \tdlc groups.  The class $\Es{\Sclass}$ carries a well-behaved rank function and is closed under taking directed unions, open subgroups, closed normal subgroups, extensions and quotients.  The class $\Es{\Sclass}$ also includes other well-studied families of \tdlc groups that are not contained in $\Es{}$, including all locally linear \tdlc groups, all complete geometric Kac--Moody groups over finite fields, the Burger--Mozes groups $\Univ(F)$ where $F$ is primitive, and $2^{\aleph_0}$ more examples of groups in $\Sclass$ that arise as groups acting on trees with Tits' independence property $\propP{}$.  On the other hand, $\Es{\Sclass}$ excludes the Burger--Mozes groups $\Univ(F)$ where $F$ is nilpotent and does not act freely.  By contrast, a larger class $\Es{\Bigsclass}$ (with similar closure properties to $\Es{\Sclass}$) is closed under forming actions on trees with property $\propP{}$.
\end{abstract}

\tableofcontents


\section{Introduction}

\subsection{Background}

Arguments by induction on some invariant are ubiquitous in group theory.  For example, in finite group theory one can appeal to induction on the order; in the theory of soluble groups, induction on the derived length; for Lie groups, induction on the dimension; and so on.  Such a basis for induction is not known in the class of totally disconnected, locally compact (\tdlc) groups, and at first glance it seems unreasonable to hope for one, given that the class includes all groups with the discrete topology.

Recent results, however, have shown that if we can find a way to put questions about discrete groups to one side, then \tdlc groups do exhibit meaningful finiteness properties.  For technical reasons it is sometimes useful to make a topological countability assumption, e.g. to consider only \tdlc second-countable (\tdlcsc) groups; this is not such a consequential restriction in practice.  The general picture that has emerged is as follows:
\begin{enumerate}[(1)]
\item Any \tdlc group $G$ is the directed union of the set $\Oc(G)$ of compactly generated open subgroups of $G$.  So we can obtain insights about the general case by combining two ingredients: an understanding of the structure of compactly generated \tdlc groups, and an understanding of which properties pass from compactly generated open subgroups (or shared properties of all sufficiently large compactly generated open subgroups) to the ambient group.  This is the ``regional approach'' to \tdlc groups, where ``regional'' refers to those properties shared by all sufficiently large compactly generated open groups.  (The analogous adjective in abstract group theory is ``local'', e.g. ``locally finite'', but in the context of topological groups, we prefer to restrict the use of the word ``local'' to the topological sense, that is, pertaining to neighbourhoods of the identity.)
\item Once we specialize to compactly generated \tdlc groups $G$, there is significant extra structure.  The first key result is Abels' discovery (\cite[Beispiel~2.7]{Abels}) that $G$ acts properly (in particular, with compact kernel) and vertex-transitively on a connected locally finite graph $\Gamma$; the quasi-isometry class of the graph is then uniquely determined.  Thus most of the methods of geometric group theory apply; see \cite{CorDLH}.  Moreover, unlike the case of a finitely generated group acting on its Cayley graph, there is additional structure coming from the action of $G_v$ on the neighbours of $v$, where $v$ is a vertex of $\Gamma$.  This was used in \cite{RW-EC} to show that $G$ admits a finite normal series in which each factor is compact, discrete, or a chief factor of $G$.  Here it is also natural to focus on the second-countable case, since every compactly generated \tdlc group is compact-by-(second-countable).
\item Given the last point, one might hope for a decomposition theory of compactly generated \tdlcsc groups similar to finite groups or Lie groups (at least if one is prepared to ignore the compact and discrete factors), where the group decomposes into factors of known types.  Additional difficulties emerge, however, because in general, the chief factors are not themselves compactly generated.  So one is led to study the compactly generated open subgroups of chief factors; these compactly generated groups themselves have chief factors (as long as they are sufficiently far from compact and discrete groups); and so on.  For an effort to take a systematic approach to the decomposition theory of chief factors, see \cite{RW-Polish} and \cite{RW-DenseLC}.  It is not clear if this process terminates in general, but if it does not, that implies the existence of a certain kind of descending chain $(G_i)_{i \in \Nb}$ of closed subgroups, where each term is open or normal in the previous term.
\item If the process of alternating between taking chief factors and passing to compactly generated open subgroups does terminate in some \tdlc group $H$, then $H$ is compactly generated and we are in one of two cases.  The first case is that $H$ has no nondiscrete chief factors, which means it has a finite normal series with only compact and discrete factors; we can regard this as analogous to the soluble case of Lie theory.  The second is that $H$ is nondiscrete but has no proper nontrivial closed normal subgroups, so $H$ belongs to the class $\Sclass$ of nondiscrete compactly generated, topologically simple \tdlc groups.  There is little prospect of classification of groups in $\Sclass$ at present, but some general properties are known (see \cite{CRW-Part2}).  The most successful approach to date in describing general properties of $H \in \Sclass$ has been via the local structure of $H$, and indeed there are nontrivial restrictions on the local isomorphism type of groups in $\Sclass$.  This suggests that one might be able to continue the analysis by understanding properties of the proper open subgroups of $H$.\end{enumerate}

The last two points suggest a route to proving general structural results about \tdlc groups, so long as we have sufficiently powerful well-foundedness properties, that is, properties that forbid certain kinds of infinite descending chains of subgroups.  Any infinite profinite group has an infinite descending chain of open normal subgroups, and there is little hope of controlling the subgroup structure of discrete groups, so we have to be careful in specifying which descending chains of subgroups to consider.

In \cite{WesEle}, P. Wesolek introduced a large class of \tdlc groups that have such a well-foundedness property.  The class of \defbold{elementary \tdlcsc groups} $\Es{}_{\aleph_0}$ is the smallest class of \tdlcsc groups that contains the profinite and discrete \tdlcsc groups and is closed under \tdlcsc group extensions and increasing unions of open subgroups.  One can easily extend to the \defbold{regionally elementary groups} $\Es{}$, which has the same definition, but within the class of first-countable (equivalently, metrizable) \tdlc groups instead of second-countable \tdlc groups.  The classes $\Es{}_{\aleph_0}$ and $\Es{}$ are also closed under closed subgroups and quotients, and admit an ordinal-valued rank function $\xi$ that can be characterized as follows: if $G = \triv$ then, $\xi(G) = 1$; if $G$ is nontrivial and compactly generated, then $\xi(G) = \xi(\Res{}(G)) +1$; and otherwise
\[
\xi(G) = \sup \{ \xi(\Res{}(O)) \mid O \in \Oc(G) \} + 1,
\]
where $\Res{}(O)$ is the intersection of all open normal subgroups of $O$.  The class $\Es{}$ can equivalently be characterized among first-countable \tdlc groups as follows: $G \in \Es{}$ if and only if there does not exist an infinite descending chain $(G_i)_{i \in \Nb}$ of closed compactly generated subgroups such that $G_{i+1} \le \Res{}(G_i)$ for all $i$.  The class also contains many \tdlcsc groups of interest; for example, every locally soluble \tdlcsc group is elementary (\cite[Theorem~8.1]{WesEle}) and at present, all known examples of amenable \tdlcsc groups are elementary.  More broadly, the elementary decomposition rank gives a way of formalizing that a group is ``far from compact or discrete'', for example in the context of chief factors, where some complications can only occur with elementary groups of small rank.

The class $\Es{}$ and its decomposition rank are thus very useful for the general theory, but they cannot account for all \tdlcsc groups.  In particular, the class $\Sclass$ is disjoint from $\Es{}$, and hence any \tdlcsc group having a group in $\Sclass$ as a quotient of a subgroup is not elementary.  At the time of writing, it is an important open question in the structure theory of \tdlcsc groups whether the class $\Sclass$ is the only obstacle to a \tdlcsc group belonging to $\Es{}$.  Specifically, let us say a class $\mc{C}$ of \tdlc groups satisfies the \defbold{elementary dichotomy} if for all $G \in \mc{C} \smallsetminus \Es{}$, there is some compactly generated closed subgroup $H$ of $G$ and a closed normal subgroup $K$ of $H$, such that $H/K \in \Sclass$.

\begin{quest}[{\cite[19.70]{Kour}}]\label{que:dichotomy}
Does the class $\mc{C}$ of all \tdlcsc groups satisfy the elementary dichotomy?
\end{quest}

\subsection{A larger class of well-founded groups}

Given the observations above, the goal of the present paper is to introduce a class $\Es{\Sclass}$ with good closure properties, properly containing $\Es{}$ and containing some groups in $\Sclass$, for which a well-foundedness principle still applies, and which satisfies the elementary dichotomy.  Given the overall approach of recent work, it is important that the class be closed under the following (at least within the class of \tdlcsc groups): directed unions, open subgroups, closed normal subgroups, extensions and quotients.  So, for example, we cannot simply force $\Sclass \subseteq \Es{\Sclass}$, since this would not account for the complexity of open subgroups of groups in $\Sclass$.  It should be stressed though that the class $\Es{\Sclass}$ given below is not intended as the definitive notion of well-foundedness for \tdlc groups; as will see later, $\Es{\Sclass}$ does not include all \tdlcsc groups.  It should instead be taken as a detailed illustration of an approach to the theory of \tdlc groups, that can either be generalized to encompass more kinds of ``well-founded'' \tdlc groups, or else contribute to an understanding of the limitations of such an approach.

\begin{defn}
Given a compactly generated \tdlcsc group $H$, we recall (\cite[Proposition~3.19]{RW-DenseLC}) that among closed normal subgroups $N$ of $H$ such that $H/N$ is elementary with $\xi(H/N)<\omega$, there is a smallest one, denoted $N = \Res{\omega}(H)$.

Let $G$ be a first-countable \tdlc group.  We define a partial order $\ll$ on the set $\Oc(G)$ of compactly generated open subgroups of $G$ as follows: say $H \ll K$ if $K$ is noncompact, $H \le K$, and the following two conditions are met:
\begin{enumerate}[(a)]
\item $\ol{H\Res{\omega}(K)}/\Res{\omega}(K)$ is compact;
\item For all closed normal subgroups $N$ of $\Res{\omega}(K)$ such that $\Res{\omega}(K)/N \in \Sclass$ and $|K:\N_K(N)|<\infty$, then $\Res{\omega}(H) \le N$.
\end{enumerate}
We say $G$ is \defbold{$\Sclass$-well-founded}, and write $G \in \Es{\Sclass}$, if the poset $(\Oc(G),\ll)$ is well-founded, that is, every nonempty open subset has at least one minimal element.  We then define $\xip{\Sclass}(G)$ to be the well-founded rank of $(\Oc(G),\ll)$ (see Section~\ref{sec:poset}), and given an ordinal $\alpha$, write $\Es{\Sclass}(\alpha)$ for the class of $\Sclass$-well-founded groups $G$ with $\xip{\Sclass}(G) \le \alpha$.
\end{defn}

For many results, we can replace the class $\Sclass$ with a larger class of topologically simple groups, such as $\Bigsclass \smallsetminus \Es{}$ where $\Bigsclass$ consists of all topologically simple \tdlcsc groups (not necessarily compactly generated), and analogously define the class $\Es{\Bigsclass}$ of $\Bigsclass$-well-founded groups.  The latter class is strictly larger than $\Es{\Sclass}$, but it should be noted that groups in $\Bigsclass \smallsetminus (\Sclass \cup \Es{})$ are not ``irreducible'' from the regional perspective, since they are directed unions of proper open subgroups.

The rank is not stable under passage to arbitrary closed subgroups, but it is stable with respect to a large class of subgroups introduced in \cite{ReidDistal}.

\begin{defn}
We say $H$ is a \defbold{RIO subgroup} of $G$ if $H$ is closed and every $K \in \Oc(H)$ is an intersection of open subgroups of $G$.
\end{defn}

Not all closed subgroups are RIO.  However, being a RIO subgroup is a transitive property (that is, all RIO subgroups of a RIO subgroup are RIO) and both open subgroups and closed normal subgroups are RIO, so it is a sufficiently rich class of subgroups for developing the decomposition theory of \tdlc groups.

Given an ordinal $\alpha$, write $\alpha^+$ for the least successor ordinal $\beta$ such that $\beta \ge \alpha$.  Given a set $S$ of ordinals, define $\sup^+(S) := (\sup(S))^+$.  Given a successor ordinal $\alpha$, we write $\alpha-1$ for the immediate predecessor of $\alpha$.

We now give a number of closure properties of the class $\Es{\Sclass}$, including rank inequalities.

\begin{thm}\label{intro:closure}
Let $G$ be a first-countable \tdlc group.
\begin{enumerate}[(i)]
\item(Proposition~\ref{prop:elementary}) If $G$ is regionally elementary, then $G \in \Es{\Sclass}$ and
\[
\xip{\Sclass}(G) \le \xi(G) \le \omega.\xip{\Sclass}(G) + 1.
\]
\item(Lemma~\ref{lem:rank_inequalities}(i)) Suppose that $G \in \Es{\Sclass}$ and $H$ is a RIO subgroup of $G$.  Then $H \in \Es{\Sclass}$ and $\xip{\Sclass}(H) \le \xip{\Sclass}(G)$.
\item(Lemma~\ref{lem:rank_inequalities}(iii)) Suppose that $\mc{D}$ is a family of RIO subgroups of $G$, directed upwards by inclusion, such that $D \in \Es{\Sclass}$ for all $D \in \mc{D}$ and $\bigcup_{D \in \mc{D}}D$ is dense in $G$.  Then $G \in \Es{\Sclass}$ and $\xip{\Sclass}(G) = {\sup_{D \in \mc{D}}}^+\xip{\Sclass}(D)$.
\item(See Theorem~\ref{thm:extensions}) Let $N$ be a closed normal subgroup of $G$.  Then $G \in \Es{\Sclass}$ if and only if $N,G/N \in \Es{\Sclass}$.  If $G \in \Es{\Sclass}$, then
\[
\max\{\xip{\Sclass}(N),\xip{\Sclass}(G/N)\} \le \xip{\Sclass}(G) \le (\xip{\Sclass}(N)-1) + \xip{\Sclass}(G/N).
\]
\item(See Proposition~\ref{prop:residual}) Let $\mc{N}$ be a family of closed normal subgroups of $G$.  Suppose that $\bigcap_{N \in \mc{N}}N = \triv$ and that $G/N \in \Es{\Sclass}$ for all $N \in \mc{N}$.  Let $\alpha = \sup\{\xip{\Sclass}(G/N) \mid N \in \mc{N}\}$.  Then $G \in \Es{\Sclass}$ and
\[
\alpha \le \xip{\Sclass}(G) \le 1+ \alpha+1.
\]
\end{enumerate}
\end{thm}

We also show that $\Es{\Sclass}$ is closed under a certain construction that can be used to produce groups in $\Es{\Sclass}$ of rank at least $\omega^2+1$: see Section~\ref{sec:ended_tree}.

For the notion of an $\Sclass$-well-founded group to be useful, there should be a systematic way, given a compactly generated $G \in \Es{\Sclass}$, to obtain a closed normal subgroup $N$ of $G$ of smaller rank, such that $G/N$ is of a special form.  For the class $\Es{}$, the natural statement of this form is the following: if $G \in \Es{}$ is compactly generated, then $\xi(\Res{}(G)) < \xi(G)$, while $G/\Res{}(G)$ is a \defbold{SIN group}, that is, it has a base of neighbourhoods of the identity consisting of compact open normal subgroups.

Here is the analogous decomposition coming from the partial order $\ll$.

\begin{thm}[See Theorem~\ref{thm:rank_reduction}]\label{intro:rank_reduction}
Let $G$ be a noncompact compactly generated \tdlcsc group such that $G \in \Es{\Sclass}$.  Then $G$ admits a finite series
\[
G_0 \le R_n \le \dots R_1 \le R_0 = G
\]
of closed characteristic subgroups with the following properties:
\begin{enumerate}[(i)]
\item $\xip{\Sclass}(G_0) < \xip{\Sclass}(G)$;
\item $R_n = \Res{\omega}(G)$ and $G_0$ is expressible as the intersection of $R_n$ with a finite number (possibly zero) of closed normal subgroups $N$ of $R_n$ such that $R_n/N \in \Sclass$ and $|G:\N_G(N)|<\infty$;
\item For $1 \le i \le n$, and given $H \in \Oc(R_{i-1}/R_i)$, then $H$ is a SIN group.
\end{enumerate}
\end{thm}

In particular, one can characterize the class $\Es{\Sclass}$ as follows, given the compactly generated nondiscrete topologically simple groups in the class.

\begin{thm}[See Theorem~\ref{thm:Esp_decomp}]\label{intro:Esp_decomp}
Write $\Sclass^* = \Es{\Sclass} \cap \Sclass$.  Then the class $\Es{\Sclass}$ is the smallest class of \tdlc groups such that
\begin{enumerate}[(i)]
\item $\Es{\Sclass}$ contains $\Sclass^*$, the discrete groups and the first-countable profinite groups; and 
\item $\Es{\Sclass}$ is closed under extensions that result in a \tdlc group and under directed unions of open subgroups.
\end{enumerate}
\end{thm}

\subsection{Relationship to known sources of groups in $\Sclass$}\label{sec:intro:special}

For all of the closure properties given so far, there are analogous closure properties already known for the class $\Es{}$.  In light of the description of the class $\Es{\Sclass}$ in Theorem~\ref{intro:Esp_decomp}, the difference between $\Es{}$ and $\Es{\Sclass}$ comes down to the class $\Sclass^*$ of nondiscrete, compactly generated, topologically simple $\Sclass$-well-founded groups, and so it is natural to turn to the literature on known sources of groups in $\Sclass$ to see which provide more examples of groups in $\Sclass^*$, and which give examples of \tdlcsc groups outside of $\Es{\Sclass}$.  We obtain results in some cases and briefly discuss what might be interesting classes of examples for future work (see Section~\ref{sec:other} for the latter).

In the literature there are several results restricting the open subgroups of $G$, under hypotheses that do not imply that $G$ is (regionally) elementary; we can use these results to deduce that $G \in \Es{\Sclass}$, with a bound on the rank.  Here are three families where such results are known, all of which include examples of groups in $\Sclass$.

\begin{thm}[See Section~\ref{sec:trees}]\label{intro:tree}
Let $T$ be a countable tree and let $G$ be a locally compact closed subgroup of $\Aut(T)$ with property $\propP{}$.  Given $v \in VT$, let $G_{v,1}$ be the stabilizer of the $1$-ball around $v$.  Suppose that for all $v \in VT$, the local action $G(v) := G_v/G_{v,1}$ is an $\Sclass$-well-founded group with primitive but not regular action on the neighbours of $v$.  Let $x$ and $y$ be adjacent vertices of $T$.  Then $G \in \Es{\Sclass}$ and
\[
\xip{\Sclass}(G) = \max\{\xip{\Sclass}(G(x)),\xip{\Sclass}(G(y)),2\},
\] 
except in the following case: If $\xip{\Sclass}(G(x)) = \xip{\Sclass}(G(y)) = 2$ and at least one of $G(x)$ and $G(y)$ is an infinite discrete group, then $\xip{\Sclass}(G)=3$.
\end{thm}

\begin{thm}[See Section~\ref{sec:locally_linear}]\label{intro:linear}
Let $G$ be a first-countable \tdlc group with a linear open subgroup.  Then $G \in \Es{\Sclass}(4)$.
\end{thm}

\begin{thm}[See Section~\ref{sec:KM}]\label{intro:KM}
Let $G$ be a complete geometric Kac--Moody group over a finite field.  Then $G \in \Es{\Sclass}(\omega)$.
\end{thm}

For Theorem~\ref{intro:tree} in particular, primitive actions are critical to the proof that the group is $\Sclass$-well-founded.  To contrast with this, we consider the case of a group $G$ acting on a tree with property $\propP{}$, such that the local actions are nilpotent but with some nontrivial point stabilizer; this is effectively a strong negation of primitive local actions, since the only primitive nilpotent permutation groups are those acting regularly.  In this case, we end up with a very different structure on $(\Oc(G),\ll)$ (see Theorem~\ref{thm:nilpotent_tree}).  In particular, we have the following family of examples of compactly generated \tdlcsc groups that are \emph{not} in $\Es{\Sclass}$.

\begin{thm}[See Corollary~\ref{cor:nilpotent_tree}]\label{intro:notEsp}
Let $d$ be a natural number, $d \ge 3$, and let $F$ be a nilpotent subgroup of $\Sym(d)$, such that some point stabilizer is nontrivial.  Then the Burger--Mozes group $\Univ(F)$ is not $\Sclass$-well-founded, and no quotient of a RIO subgroup of $\Univ(F)$ belongs to $\Sclass$.
\end{thm}

The groups $\Univ(F)$ in Theorem~\ref{intro:notEsp} do not belong to $\Sclass$, however they can be used to obtain non-$\Sclass$-well-founded groups in $\Sclass$.  Specifically, we can appeal to a construction due to W. Lederle (\cite{Lederle}): for each $F \le \Sym(d)$, there is an associated \defbold{coloured Neretin group} $\mc{N}_F = \mathrm{F}(\Univ(F))$, consisting of all homeomorphisms of the boundary of the tree that can be represented as piecewise combinations of elements of $\Univ(F)$. Then $\mc{N}_F$ carries a group topology in which $\Univ(F)$ is embedded as an open subgroup; moreover, by \cite[Theorem~1.2]{Lederle}, $\mc{N}_F$ is compactly generated and simple-by-(finite abelian).

\begin{cor}\label{cor:notEsp}
Let $d$ be a natural number, $d \ge 3$, and let $F$ be a nilpotent subgroup of $\Sym(d)$, such that some point stabilizer is nontrivial.  Let $\mc{N}_F$ be the associated coloured Neretin group.  Then the derived group of $\mc{N}_F$ is in $\Sclass$ but not $\Es{\Sclass}$.
\end{cor}

\subsection{Open questions}

Theorem~\ref{intro:notEsp} does not rule out the possibility that some other structural properties are sufficient to ensure $\Sclass$-well-foundedness.

Many of the $\Sclass$-well-founded groups described by Theorems~\ref{intro:tree} and \ref{intro:linear}, and all of those described by Theorem~\ref{intro:KM}, are also \defbold{Noetherian}, meaning they are \tdlc groups $G$ such that every open subgroup of $G$ is compactly generated.  By contrast, none of the groups $\Univ(F)$ described in Theorem~\ref{intro:notEsp} are Noetherian.  In fact, any Noetherian group $G$ acting on a locally finite tree with Tits' independence property belongs to $\Es{\Sclass}(3)$ (see Proposition~\ref{prop:tree_Noetherian}).   Certainly there are compactly generated $\Sclass$-well-founded groups that are not Noetherian, for instance the examples constructed in Section~\ref{sec:ended_tree} are not Noetherian.  However, the converse remains open.

\begin{quest}\label{que:Noetherian}
Let $G$ be a Noetherian \tdlcsc group.  Must $G$ be $\Sclass$-well-founded?
\end{quest}

Another interesting special case to consider is the class $\ms{P}$ of permutation groups that are nondiscrete and \defbold{primitive closed subdegree-finite}, or in other words, \tdlc groups $G$ possessing a maximal subgroup $U$, such that $U$ is compact open and the intersection of conjugates of $U$ is trivial.  Given recent work of S. Smith (\cite{SmithPrimitive}), it seems the base case for understanding the class $\ms{P}$ is to understand those $G \in \ms{P}$ where $G$ is one-ended and belongs to $\Sclass$.  Theorems~\ref{intro:linear} and \ref{intro:KM} account for many of the known examples of groups of the latter type, while the groups $G$ appearing in Theorem~\ref{intro:tree} are a major source of groups in $\ms{P} \cap \Sclass$ with infinitely many ends.

\begin{quest}
Let $G$ be a primitive closed subdegree-finite group.  Must $G$ be $\Sclass$-well-founded?  What if we also assume $G$ is in $\Sclass$ and/or one-ended?
\end{quest}

As mentioned above, the class $\Es{\Bigsclass}$ (where in the partial order $\ll$, one replaces $\Sclass$ with the class $\Bigsclass$ of all topologically simple \tdlcsc groups) is strictly larger than $\Es{\Sclass}$.  Specifically, in contrast to Theorem~\ref{intro:notEsp}, the class $\Es{\Bigsclass}$ contains all groups $G$ acting on countable trees with property $\propP{}$, such that arc stabilizers are compact and the local actions belong to $\Es{\Bigsclass}$ (see Proposition~\ref{prop:tree_Noetherian}); in particular, $\Es{\Bigsclass}$ contains all Burger--Mozes groups $\Univ(F)$ on locally finite trees.  The class $\Es{\Bigsclass}$ also contains all \tdlc groups that are locally isomorphic to a just infinite profinite group (Corollary~\ref{cor:ld_wellfounded}).  The following question remains open.

\begin{quest}\label{que:BigS}
Is every \tdlcsc group $\Bigsclass$-well-founded?
\end{quest}

Note that if the answer to Question~\ref{que:BigS} is ``yes'', then the answer to Question~\ref{que:Noetherian} is also ``yes'': given a Noetherian group $G$, there is no distinction between the $\Bigsclass$-well-founded rank and the $\Sclass$-well-founded rank, because the simple factors that are relevant to the partial order associated to $\Bigsclass$ on $\Oc(G)$ all belong to $\Sclass$.

\subsection{Structure of the article}

Section~\ref{sec:preliminaries} is a preliminary section in which we recall relevant parts of the known structure theory of \tdlc groups, with a few additional lemmas adapted for the purposes of the present article.  In Section~\ref{sec:ll} we introduce the partial order $\ll$ on a family of subgroups of a first-countable \tdlc group, with primary focus on the compactly generated open subgroups (up to finite index), and show how $\ll$ is preserved under important classes of homomorphisms of first-countable \tdlc groups.  In Section~\ref{sec:well-founded}, we establish general properties of the class $\Es{\Sclass}$ of first-countable \tdlc groups $G$ such that $(\Oc(G),\ll)$ is well-founded; in particular, the properties of $\ll$ established in Section~\ref{sec:ll} are used to show that $\Es{\Sclass}$ has desirable closure properties and to establish various inequalities in terms of the rank function associated to $(\Oc(G),\ll)$.  Finally, in Section~\ref{sec:special}, we consider the families of groups recalled in Section~\ref{sec:intro:special}.  The subsections of Section~\ref{sec:special} depend on the previous sections but can be read independently of one another.

\section{Preliminaries}\label{sec:preliminaries}

\subsection{Definitions and notation}

Let $G$ be a topological group.

Given subgroups $A$ and $B$ of $G$, define $[A,B] = \grp{aba\inv b\inv \mid a \in A, b \in B}$.

The \defbold{discrete residual} $\Res{}(G)$ is the intersection of all open normal subgroups of $G$.  Given a group $H$ acting on $G$, then $\Res{G}(H)$ is the intersection of all open $H$-invariant subgroups of $G$.  Unless otherwise specified, when $G$ and $H$ are subgroups of some other group $L$ such that $H \le \N_L(G)$, we let $H$ act on $G$ by conjugation.

The \defbold{quasi-centre} $\QZ(G)$ of $G$ is the set of elements $g \in G$ such that $\CC_G(g)$ is open.  We say $G$ is \defbold{quasi-discrete} if $\QZ(G)$ is dense in $G$.

Given subgroups $H$ and $K$ of $G$, we say $H$ \defbold{virtually contains} $K$ if $H$ contains a finite index subgroup of $K$.  We say $H$ and $K$ are \defbold{commensurate}, and write $H \sim_f K$, if each of $H$ and $K$ virtually contains the other.  A \defbold{commensurated} subgroup of $G$ is one that is commensurate with all of its $G$-conjugates.

Let $\mc{Q}$ be a set of closed normal subgroups of $G$.  Given a subset $I$ of $\mc{Q}$, write $G_I := \cgrp{N \in i}$.  Then $G$ is a \defbold{quasi-product} of $\mc{Q}$, or $\mc{Q}$ is a \defbold{quasi-direct factorization} of $G$, if 
\[
G = G_{\mc{Q}} \text{ and } \bigcap_{N \in \mc{Q}}G_{\mc{Q} \smallsetminus \{N\}} = \triv.
\]

Let $(G_i)_{i \in I}$ be a sequence of topological groups and for each $i \in I$, let $U_i$ be an open subgroup of $G_i$.  The \defbold{local direct product} (or \defbold{restricted product}) of $(G_i,U_i)_{i \in I}$, denoted $\bigoplus_{i \in I}(G_i,U_i)$, consists of all sequences $(g_i)_{i \in I}$ such that $g_i \in G_i$ for all $i \in I$ and $g_i \in U_i$ for all but finitely many $i \in I$.  The group operations on $\bigoplus_{i \in I}(G_i,U_i)$ are defined pointwise, and the topology is the unique group topology such that the natural embedding of $\prod_{i \in I}U_i$ is continuous and open.

Suppose now that $G$ is a \tdlc group.  A property of $G$ holds \defbold{locally} if the property is satisfied by any sufficiently small open subgroup, and \defbold{regionally} if the property holds for any sufficiently large compactly generated open subgroup: that is, $G$ is \defbold{regionally-$P$} if there is a compactly generated open subgroup $U$ of $G$, such that every compactly generated open subgroup $O$ of $G$ containing $U$ has $P$.

Write $\Oc(G)$ for the class of open subgroups of $G$ that are compactly generated.  The \defbold{regional poset} of $G$ is $\Reg(G) = \Oc(G)/\sim_f$, carrying the partial order $\le$ induced by inclusion.  (We will later equip $\Reg(G)$ with an additional partial order.)  The $\sim_f$-class of $H \in \Oc(G)$ is denoted $[H]_f$, or just $[H]$ if it is clear from the context that we are talking about an element of $\Reg(G)$.

\subsection{Regionally elementary groups}

The class $\Es{}_{\aleph_0}$ of elementary \tdlcsc groups is a class introduced by Wesolek in \cite{WesEle}.

\begin{defn}
The class $\Es{}_{\aleph_0}$ is the smallest class of \tdlcsc groups with the following properties:
\begin{enumerate}[(i)]
\item $\Es{}_{\aleph_0}$ contains all countable discrete groups and second-countable profinite groups;
\item Given a \tdlcsc group $G$ such that $G = \bigcup_{i \in I}O_i$, where $(O_i)_{i \in I}$ is a net of open subgroups directed under inclusion that belong to $\Es{}_{\aleph_0}$, then $G \in \Es{}_{\aleph_0}$;
\item Given a short exact sequence
\[
\triv \rightarrow N \rightarrow G \rightarrow Q \rightarrow 1
\]
of \tdlcsc groups such that $N,Q \in \Es{}_{\aleph_0}$, then $G \in \Es{}_{\aleph_0}$.
\end{enumerate}
Given a \tdlcsc group $G$, we say $G$ is \defbold{elementary} if $G \in \Es{}_{\aleph_0}$.
\end{defn}

For our purposes it is more natural to work in the class of first-countable \tdlc groups (equivalently, metrizable \tdlc groups) rather than second-countable \tdlc groups, since second-countability is neither a local nor a regional property, whereas first-countability is both.  One then considers the larger class $\Es{}$ of \defbold{regionally elementary} groups, consisting of those \tdlc groups $G$ such that every compactly generated open subgroup of $G$ is elementary.  The class $\Es{}$ inherits from $\Es{}_{\aleph_0}$ the properties of being closed under directed unions and extensions that result in a first-countable \tdlc group.  Note that a regionally elementary group is always first-countable, but need not be $\sigma$-compact: for example, a discrete group of any cardinality is regionally elementary. 

Another perspective is that the regionally elementary groups are exactly the first-countable \tdlc groups that are ``well-founded'' with respect to taking discrete residuals of compactly generated subgroups.  The class of regionally elementary groups admits a well-behaved ordinal-valued rank function $\xi$, called the \defbold{decomposition rank}, which can be characterized as follows:
\begin{enumerate}
\item $\xi(\triv) = 1$.
\item If $G$ is compactly generated, then $\xi(G) = \xi(\Res{}(G))+1$.
\item If $G$ is not compactly generated, then
\[
\xi(G) = \sup \{\xi(\Res{}(O)) \mid O \in \Oc(G)\} + 1.
\]
\end{enumerate}
Henceforth we will also refer to the decomposition rank as the \defbold{elementary rank}, to emphasize the connection with the class of (regionally) elementary groups and to distinguish it from the $\Sclass$-well-founded rank that we will define later.  (The other rank function defined in \cite{WesEle}, the construction rank, is not relevant for our purposes.)

Write $\Es{}(\alpha)$ for the class of regionally elementary groups $G$ with $\xi(G) \le \alpha$.  Note that the elementary rank can only take successor ordinals as values, so for example $\Es{}(\omega)$ is the class of regionally elementary groups of finite elementary rank.

Given a first-countable \tdlc group, we have $G \in \Es{}(2)$ if and only if $\Res{}(H) = \triv$ for all $H \in \Oc(G)$.  The following characterization of $\Es{}(2)$ follows from \cite[Corollary~4.1]{CM}.

\begin{lem}\label{lem:elementary2}
Let $G$ be a first-countable \tdlc group.  Then $G \in \Es{}(2)$ if and only if for every $H \in \Oc(G)$, there is a base of neighbourhoods of the identity consisting of compact open normal subgroups of $H$.
\end{lem}

The regionally elementary groups (of bounded rank) have some important closure properties that are not needed to characterize them.  We will use the following without further comment.

\begin{lem}\label{lem:elementary_closure}
Let $G \in \Es{}(\alpha)$.  Then $\Es{}(\alpha)$ contains all closed subgroups of $G$ and all quotients of $G$ by closed normal subgroups.
\end{lem}

\begin{proof}
See \cite[Corollary~4.10 and Theorem~4.19]{WesEle}; the arguments easily generalize from second-countable to first-countable \tdlc groups, recalling that a first-countable \tdlc group is a directed union of second-countable open subgroups.
\end{proof}

We also note the following quantitative version of the fact that $\Es{}$ is closed under extensions.  (Again the reference is for second-countable groups, but the proof easily generalizes to the first-countable case.)

\begin{lem}[{\cite[Lemma~3.8]{RW-DenseLC}}]\label{lem:elementary_extension}
Let
\[
\triv \rightarrow N \rightarrow G \rightarrow Q \rightarrow 1
\]
be a short exact sequence of first-countable \tdlc groups such that $N,Q \in \Es{}$.  Then $G \in \Es{}$ and
\[
\xi(G) \le \xi(N-1) + \xi(Q).
\]
In particular, for all ordinals $\alpha$, the class $\Es{}(\omega^\alpha)$ is closed under extensions within the class of first-countable \tdlc groups.
\end{lem}

The regionally elementary groups are exactly the first-countable \tdlc groups for which the recursion implied by the characterization of $\xi$ terminates.  More precisely, one has the following dichotomy:

\begin{prop}[{\cite[Proposition~6.2.2]{CRWes}}]\label{prop:elementary_wellfounded}
Let $G$ be a first-countable \tdlc group.  Then exactly one of the following holds:
\begin{enumerate}[(i)]
\item $G \in \Es{}$;
\item There is an infinite descending sequence $(K_i)_{i \in \Nb}$ of nontrivial compactly generated closed subgroups of $G$, such that $K_{i+1} \le \Res{}(K_i)$ for all $i$.
\end{enumerate}
\end{prop}

Note that in Proposition~\ref{prop:elementary_wellfounded}(ii), we do not necessarily have $K_{i+1} < K_i$.  For example, given $G \in \Sclass$, we have $\Res{}(G)=G$ in this case, so the constant sequence $K_i = G$ for all $i \in \Nb$ satisfies (ii).  This will be a key distinction from the larger class of $\Sclass$-well-founded groups that we will define later.

If $G$ is a \tdlcsc group that is compactly generated but not compact, and $G/\Res{}(G)$ is compact, then $G$ cannot be elementary.  There are two ways to see this.  First, using the decomposition rank, one has the following.

\begin{lem}[{\cite[Lemma~3.10]{RW-DenseLC}}]\label{lem:elementary:cocompact}
Let $G$ be a \tdlcsc group and let $N$ be a closed normal subgroup such that $N$ is elementary and $G/N$ is compact.  Then $G$ is elementary, and either $N$ is trivial or $\xi(N) = \xi(G)$.  In particular, if $N = \Res{}(G)$, then $N = \triv$, so $G$ is compact.
\end{lem}

Second, the following theorem of P.-E. Caprace and N. Monod shows that in fact $G$ must involve groups in $\Sclass$.

\begin{thm}[{See \cite[Theorem~A]{CM}}]\label{thm:CM}
Let $G$ be a compactly generated \tdlc group such that $G$ is not compact and $G/\Res{}(G)$ is compact.  Then $\Res{}(\Res{}(G)) = \Res{}(G)$, and the set $\mc{N}$ of closed normal subgroups $N$ of $\Res{}(G)$ such that $\Res{}(G)/N \in \Sclass$ is finite and nonempty.  Moreover, every proper closed normal subgroup of $\Res{}(G)$ is contained in some $N \in \mc{N}$.
\end{thm}

Given a \tdlcsc group $G$ and an ordinal $\alpha$, the \defbold{rank $\alpha$ residual} $\Res{\alpha}(G)$ is the intersection of all closed normal subgroups $N$ such that $G/N \in \Es{}(\alpha)$.  In particular, note that if $G$ is compactly generated, then $\Res{2}(G) = \Res{}(G)$: compactly generated \tdlc groups of elementary rank $2$ are residually discrete, while nontrivial discrete groups have elementary rank $2$.  More interesting is the \defbold{finite elementary rank residual} $\Res{\omega}(G)$.  By \cite[Theorem~3.25]{RW-DenseLC}, $\Res{\omega}(G)$ is the intersection of the \defbold{elementary rank-$2$ series} $(G_i)$ of $G$, which is formed by setting $G_0 = G$ and thereafter $G_{i+1} = \Res{2}(G_i)$.

\begin{lem}\label{lem:finite_rank_residual}
Let $G$ be a compactly generated \tdlcsc group.
\begin{enumerate}[(i)]
\item For every countable ordinal $\alpha$, then $G/\Res{\alpha}(G) \in \Es{}(\alpha)$.  In particular, $G/\Res{\omega}(G)$ has finite elementary rank.
\item Let $H$ be a closed subgroup of $G$ such that $\Res{\omega}(G) \le H$.  Then $\Res{\omega}(G) = \Res{\omega}(H)$.  In particular, $\Res{\omega}(\Res{\omega}(G)) = \Res{\omega}(G)$, that is, $\Res{\omega}(G)$ has no nontrivial quotient with finite elementary rank.
\end{enumerate}
\end{lem}

\begin{proof}
Part (i) is given by \cite[Proposition~3.19]{RW-DenseLC}.

Given a \tdlcsc group $A$, write $(A_i)$ for the elementary rank-$2$ series of $A$.  Given a closed subgroup $B$ of $A$, it is straightforward to see (for example using \cite[Corollary~3.21]{RW-DenseLC}) that $B_i \le A_i$ for all $i$; in particular, $\Res{\omega}(B) \le \Res{\omega}(A)$.

Let $K = \Res{\omega}(G)$ and form the elementary rank-$2$ series $(G_i)$ and $(K_i)$ of $G$ and $K$ respectively.  Then by (i), we have $K = G_{j_0}$ for some $j_0 \ge 0$.  We then see that in fact $K_i = G_{j_0+i} = K$ for all $i \ge 0$, so $\Res{\omega}(K) = K$.  Now given $H \le G$ such that $K \le H$, then on the one hand $\Res{\omega}(K) \le \Res{\omega}(H)$, so $K \le \Res{\omega}(H)$, but on the other hand $\Res{\omega}(H) \le \Res{\omega}(G)$, so $\Res{\omega}(H) \le K$.  Thus $\Res{\omega}(H) = K$, proving (ii).
\end{proof}

Note that all first-countable abelian \tdlc groups belong to $\Es{}(2)$.  In fact, the following holds:

\begin{lem}[{See \cite[Lemma~3.16]{RW-DenseLC}}]\label{lem:quasi-discrete}
Let $G$ be a first-countable \tdlc group that is quasi-discrete.  Then $G \in \Es{}(2)$.
\end{lem}

Thus in any compactly generated \tdlcsc group $G$, then $\Res{\omega}(G)$ has no nontrivial quasi-discrete Hausdorff quotient, and in particular $\Res{\omega}(G)$ is topologically perfect.

Given a nondiscrete elementary \tdlcsc group $G$, then all possible ranks are witnessed by $H \in \Oc(G)$, subject to the basic restrictions that $\xi(H) = \alpha+2$ for some $\alpha$ and $1 < \xi(H) \le \xi(G)$.

\begin{lem}[{See \cite[Theorem~4.19]{ReidDistal}}]\label{lem:ele_rank_witnesses}
Let $G$ be a nondiscrete elementary \tdlcsc group and let $\alpha$ be an ordinal such that $\alpha+2 \le \xi(G)$.  Then there is $H \in \Oc(G)$ such that $\xi(H) = \alpha+2$.
\end{lem}

We note that in the case of characteristically simple groups, only certain ranks can occur.

\begin{lem}[{\cite[Theorem~3.17]{RW-DenseLC}}]\label{lem:chief_rank}
Let $G$ be an elementary \tdlcsc group that is nontrivial and topologically characteristically simple.  Then $\xi(G)$ is either $2$ or the successor of an infinite limit ordinal.
\end{lem}

As a special case of the regionally elementary groups, we have the first-countable regionally elliptic \tdlc groups, where a locally compact group is \defbold{regionally elliptic} if every compact subset is contained in a compact subgroup.  Write $\Es{}(\mathrm{r})$ for the first-countable regionally elliptic \tdlc groups; note that $\Es{}(\mathrm{r}) \subseteq \Es{}(2)$.

By a theorem of V. Platonov (\cite{Platonov}), in any locally compact group $G$ there is a unique largest regionally elliptic closed normal subgroup $\RadRE(G)$, the \defbold{regionally elliptic radical} of $G$, and 
\[
\RadRE(G/\RadRE(G)) = \triv.
\]

\subsection{Some classes of topologically simple groups}\label{sec:simple_def}

Recall that $\Sclass$ is the class of \tdlc groups that are nondiscrete, compactly generated and topologically simple.  Note that all groups in $\Sclass$ are also second-countable: indeed, every compactly generated \tdlc group is $\sigma$-compact, and it is straightforward to verify, given Van Dantzig's theorem, that every $\sigma$-compact \tdlc group has a second-countable quotient with compact kernel.

For many results we can generalize from $\Sclass$ to a possibly larger class $\mc{S}$ of topologically simple groups.  To make sure $\mc{S}$ has the right properties, let us say that a class $\mc{S}$ of \tdlcsc groups has the property $\propS$ if the following holds:
\begin{enumerate}[(a)]
\item $\mc{S} \smallsetminus \Es{}$ consists of topologically simple groups;
\item We have $\Sclass \subseteq \mc{S}$;
\item Given a continuous homomorphism $\psi: G \rightarrow H$ of \tdlcsc groups with dense normal image:
\begin{enumerate}[(i)]
\item If $G$ is a topologically perfect central extension of a group in $\mc{S} \smallsetminus \Es{}$, then the same is true of $H$.
\item If $H \in \mc{S} \smallsetminus \Es{}$ and $\psi$ is injective, then $\ol{[G,G]} \in \mc{S} \smallsetminus \Es{}$.
\end{enumerate}
\end{enumerate}

For the application we will exclude elementary groups, so in the context where we use property $\propS$, the distinction between $\mc{S}$ and $\mc{S} \smallsetminus \Es{}$ will not matter.

In particular, let $\Bigsclass$ be the class of all topologically simple \tdlcsc groups.  An interesting class between $\Sclass$ and $\Bigsclass$ comes from the results of \cite{CRWes}.  To define this class we need to recall a few definitions.

\begin{defn}
A \tdlc group $G$ is \defbold{expansive} if there is a neighbourhood $U$ of the identity in $G$ such that $\bigcap_{g \in G}gUg\inv = \triv$.  We say $G$ is \defbold{regionally expansive} if there exists $O \in \Oc(G)$ such that $O$ is expansive.  A \tdlc group $G$ is \defbold{robustly monolithic} if the intersection $\Mon(G)$ of all nontrivial closed normal subgroups of $G$ is nondiscrete, regionally expansive and topologically simple.  Write $\ms{R}$ for the class of robustly monolithic groups.
\end{defn}

The class $\ms{R}$ was introduced in \cite{CRWes}, where it was introduced as a class containing $\Sclass$ and with the following additional property: given a continuous injective dense homomorphism $H \rightarrow G$ of \tdlc groups such that $G \in \ms{R}$ and $H$ is nondiscrete, then $H \in \ms{R}$ (\cite[Theorem~5.4.1]{CRWes}).  

We will see in the next subsection (Lemma~\ref{lem:semisimple:basic}) that each of the classes $\Sclass$, $\ms{R} \cap \Bigsclass$ and $\Bigsclass$ has $\propS$.

\subsection{Normal compressions}

A \defbold{normal compression} of topological groups is a continuous injective group homomorphism $\phi: G \rightarrow H$ with dense normal image.  In this subsection we recall some results on normal compressions obtained by the author and Wesolek in \cite{RW-Polish} and \cite{RW-DenseLC}.

The following was stated for normal compressions in \cite[Proposition~3.5]{RW-Polish}, however the hypothesis that $\psi(G)$ is dense in $H$ was not used in the proof.

\begin{prop}\label{prop:compression_continuous}
Let $G$ and $H$ be Polish groups and let $\psi: G \rightarrow H$ be a continuous injective homomorphism such that $\psi(G)$ is normal in $H$.  Then there is a unique jointly continuous action of $H$ on $G$ such that $\psi$ is equivariant with respect to this action and the conjugation action of $H$ on itself.
\end{prop}

Under the hypotheses of Proposition~\ref{prop:compression_continuous}, we write $G \rtimes_\psi H$ for the semidirect product associated to the given action of $H$ on $G$; note that $G \rtimes_\psi H$ is then also a Polish group.

\begin{thm}[{\cite[Theorem~3.6]{RW-Polish}}]\label{thm:psi-compression_factor_rel}
Let $G$ and $H$ be Polish groups, let $\psi: G \rightarrow H$ be a normal compression, and let $O\leq H$ be an open subgroup. Then the following hold:
\begin{enumerate}[(i)]
\item The map $\pi:G\rtimes_{\psi}O\rightarrow H$ via $(g,o)\mapsto \psi(g)o$ is a continuous surjective homomorphism with $\ker(\pi)=\{(g\inv,\psi(g))\mid g\in \psi^{-1}(O)\}$, and if $O = H$, then $\ker(\pi)\cong G$ as topological groups.
\item Writing $\iota: G \rightarrow G\rtimes_{\psi}O$ for the usual inclusion map, then $\psi = \pi \iota$.
\item We have $G\rtimes_{\psi}O =\ol{\iota(G)\ker (\pi)}$, and the subgroups $\iota(G)$ and $\ker(\pi)$ are closed normal subgroups of $G\rtimes_{\psi}O$ with trivial intersection.
\end{enumerate}
\end{thm}

\begin{lem}\label{lem:compression_normal_subgroup}
Let $G$ and $H$ be Polish groups and let $\psi: G \rightarrow H$ be a continuous homomorphism with dense normal image.  Let $K$ be a normal subgroup of $G$.  Then $\ol{\psi(K)}$ is normal in $H$; if $\ker(\psi) K$ is closed in $G$, then $\psi(K)$ is normal in $H$.
\end{lem}

\begin{proof}
Let $L = \ol{\ker(\psi) K}$.  Then $L/\ker(\psi)$ is a closed normal subgroup of $G/\ker(\psi)$, and $\psi$ induces a normal compression from $G/\ker(\psi)$ to $H$.  By \cite[Corollary~3.7]{RW-Polish}, we have $\psi(L) \unlhd H$.  We see by continuity that $\ol{\psi(K)} = \ol{\psi(L)}$, so $\ol{\psi(K)}$ is normal in $H$.  If $\ker(\psi) K$ is closed, we have $\psi(L) = \psi(\ker(\psi) K) = \psi(K)$, so $\psi(K)$ is normal in $H$. 
\end{proof}

The next lemma is a variant of \cite[Proposition~3.8]{RW-Polish}.

\begin{lem}\label{lem:normal_compression:derived}
Let $G$ and $H$ be Polish groups, let $\psi: G \rightarrow H$ be a normal compression, let $K$ be a closed normal subgroup of $G$ and let $L = \ol{\psi(K)}$.  Then
\[
\ol{[G,\psi\inv(L)]} \le K \text{ and } [H,\psi(G) \cap L] \le \psi(K).
\]
\end{lem}

\begin{proof}
Form the semidirect product $G\rtimes_{\psi} H$, let $\iota:G\rightarrow G\rtimes_{\psi}H$ be the usual inclusion, and let $\pi:G\rtimes_{\psi}H \rightarrow H$ be the map $(g,h) \mapsto \psi(g)h$.  Note that by Theorem~\ref{thm:psi-compression_factor_rel}, $\pi$ is a quotient map and $\psi = \pi \iota$.  In particular, $M = \psi\inv(L)$ can be recovered from the semidirect product as follows:
\[
M = \iota\inv(\ol{\iota(K)\ker\pi}).
\]
We first claim that $[G,M] \le K$.  Take $g \in G$ and $m \in M$; it suffices to show $[g,m] \in K$.  We see that $\iota(m) \in \ol{\iota(K)\ker\pi}$, so there are nets $(k_i)_{i \in I}$ and $(r_i)_{i \in I}$, with $k_i \in \iota(K)$ and $r_i \in \ker\pi$, such that $k_ir_i \rightarrow \iota(m)$.  By Theorem~\ref{thm:psi-compression_factor_rel}, $\iota(G)$ commutes with $\ker\pi$, so for each $i \in I$ we can write $[\iota(g),k_ir_i] = [\iota(g),k_i]$.  In particular, we see that $[\iota(g),k_i]$ converges to $[\iota(g),\iota(m)]$ as $i \rightarrow \infty$.  Since $\iota(K)$ is a closed normal subgroup of $G \rtimes_{\psi} H$, it follows that $[\iota(g),\iota(m)] \in \iota(K)$, and hence $[g,m] \in K$.

Now we claim that
\[
[H,\psi(G) \cap L] \le \psi(K).
\]
Consider again the semidirect product $G\rtimes_{\psi} H$.  By Lemma~\ref{lem:compression_normal_subgroup} we see that $\iota(M)$ and $\iota(K)$ are both closed normal subgroups of the semidirect product, so there is a quotient $G/K \rtimes H$, which has a closed normal subgroup $M/K \rtimes H$.  The action of $H$ on $M/K$ is again continuous, but there is a dense subgroup $\psi(G)$ of $H$ that acts trivially on $M/K$.  It follows that $H$ acts trivially on $M/K$.  Converting back to subgroups of $H$, it follows that $[H,\psi(M)] \le \psi(K)$.  The conclusion follows by noting that $\psi(M) = \psi(G) \cap L$.
\end{proof}

Normal compressions of \tdlcsc groups preserve elementary rank.

\begin{lem}[{\cite[Proposition~5.4]{RW-DenseLC}}]\label{lem:compression_rank}
Let $G$ and $H$ be \tdlcsc groups, let $\psi: G \rightarrow H$ be a normal compression and let $\alpha$ be an ordinal.  Then $G \in \Es{}(\alpha)$ if and only if $H \in \Es{}(\alpha)$.
\end{lem}

We recall the property $\propS$ defined in Section~\ref{sec:simple_def}, and confirm it indeed applies to $\Sclass$, $\Bigsclass$ and $\ms{R} \cap \Bigsclass$.

\begin{lem}\label{lem:semisimple:basic}
Let $\mc{S}$ be either the class $\Sclass$ of nondiscrete compactly generated topologically simple \tdlc groups, or the class $\Bigsclass$ of nondiscrete topologically simple \tdlcsc groups, or the class $\ms{R} \cap \Bigsclass$.  Then $\mc{S}$ has property $\propS$.
\end{lem}

\begin{proof}
The only nontrivial part of the definition of $\propS$ to check is condition (c), namely, if we are given a continuous homomorphism $\psi: G \rightarrow H$ of \tdlcsc groups with dense normal image, then the following holds:
\begin{enumerate}[(i)]
\item If $G$ is a topologically perfect central extension of a group in $\mc{S} \smallsetminus \Es{}$, then the same is true of $H$.
\item If $H \in \mc{S} \smallsetminus \Es{}$ and $\psi$ is injective, then $\ol{[G,G]} \in \mc{S} \smallsetminus \Es{}$.
\end{enumerate}
We see that in (i) we have $G \not\in \Es{}$ and in (ii) we have $H \not\in \Es{}$.  By Lemma~\ref{lem:compression_rank} it follows that $G,H \not\in \Es{}$ in both cases.  Consequently, we see that given a quotient $K$ of $G$ or $H$ with kernel in $\Es{}$, or a closed normal subgroup $K$ of $G$ or $H$ such that the quotient is in $\Es{}$, then $K \not\in \Es{}$.  In particular, by Lemma~\ref{lem:quasi-discrete} we have $\ol{\QZ(K)} < K$ in this case.

For (i), suppose $G$ is topologically perfect and $G/\Z(G) \in \mc{S}$.  Then given a closed normal subgroup $K$ of $G$ such that $K \nleq \Z(G)$, we have $G = \ol{K\Z(G)}$, so $G/K$ is abelian and hence $K = G$.  In other words, in the terminology of \cite{RW-Polish}, $G$ is of semisimple type with only one component.  By \cite[Proposition~5.15(1)]{RW-Polish}, $H$ is likewise of semisimple type with one component, so $H$ is topologically perfect and $H/\Z(H)$ is topologically simple.  In particular, $H/\Z(H) \in \Bigsclass$.  It is then clear that we have a normal compression from $G/\Z(G)$ to $H/\Z(H)$.  If $G/\Z(G) \in \ms{R}$, then $H/\Z(H) \in \ms{R}$ by \cite[Proposition~5.3.1]{CRWes}.  If $G/\Z(G)$ is also compactly generated, we see that the same is true of $H/\Z(H)$; hence $H/\Z(H) \in \Sclass$ in this case.

For (ii), suppose $H \in \mc{S}$ and $\psi$ is injective.  We see that every nontrivial closed normal subgroup of $G$ has dense image in $H$; it is then clear that $G$ has no nontrivial abelian normal subgroup.  By Lemma~\ref{lem:normal_compression:derived}, we deduce that $K:=\ol{[G,G]}$ is the smallest nontrivial closed normal subgroup of $G$; in particular, $K$ is topologically perfect.  Applying the same argument again to the normal compression from $K$ to $H$, we see that $K$ is topologically simple, so $K \in \Bigsclass$.  If $H \in \ms{R}$, then $K \in \ms{R}$ by \cite[Theorem~5.4.1]{CRWes}.  Now suppose that $H \in \Sclass$, that is, $H$ is compactly generated.  Then by \cite[Corollary~5.3]{RW-DenseLC}, $K$ is compactly generated, so $K \in \Sclass$.
\end{proof}

In this article we will be considering first-countable \tdlc groups; such topological groups are only regionally Polish in general.  However, with an additional assumption about the homomorphisms under consideration, we can still apply results about normal compressions of Polish groups to homomorphisms between first-countable \tdlc groups.

\begin{defn}
Let $G$ and $H$ be \tdlc groups.  A \defbold{regionally normal map} from $G$ to $H$ is a continuous homomorphism $\phi: G \rightarrow H$ such that for every $O \in \Oc(G)$, then $\N_H(\phi(O))$ is open in $H$.  We say $\phi$ is \defbold{dense} if it has dense image; a \defbold{regionally normal compression} is an injective dense regionally normal map.
\end{defn}

Notice that continuous open homomorphisms (for example, quotient homomorphisms) are regionally normal, and that the restriction of a regionally normal map to an open subgroup is still regionally normal.  As the next lemma shows, every continuous homomorphism between \tdlcsc groups with dense normal image is regionally normal; it is also enough to check that every \emph{sufficiently large} $O \in \Oc(G)$ has locally normal image in $H$.

\begin{lem}\label{lem:sigma-normal}
Let $G$ and $H$ be first-countable \tdlc groups and let $\phi: G \rightarrow H$ be a continuous homomorphism.  Suppose that for all $A \in \Oc(G)$, there is a second-countable open subgroup $B$ of $G$ containing $A$ such that $\N_H(\phi(B))$ is open in $H$.  Then $\phi$ is regionally normal.
\end{lem}

\begin{proof}
We can decompose $\phi$ as $\phi = \psi \pi$, where $\pi: G \rightarrow G/N$ is a quotient map and $\psi: G/N \rightarrow H$ is injective.  Given $O \in \Oc(G)$, then $ON/N \in \Oc(G/N)$, and every element of $\Oc(G/N)$ occurs in this way; we also have $\phi(O) = \psi(ON/N)$.  So we may assume without loss of generality that $\phi$ is injective.

Let $O \in \Oc(G)$; we aim to show that $\N_H(\phi(O))$ is open in $H$.  By hypothesis there is a second-countable open subgroup $B$ of $G$ containing $O$ such that $\N_H(\phi(B))$ is open in $H$.  Without loss of generality we may assume $B = G$; we can then also replace $H$ with the open subgroup $\grp{\phi(G),V}$ of $H$, where $V$ is some compact open subgroup of $\N_H(\phi(B))$.  Under these assumptions, $G$ and $H$ are Polish and $\phi(G)$ is normal in $H$.

By Proposition~\ref{prop:compression_continuous}, we can make $G$ and $H$ closed subgroups of the semidirect product $G \rtimes_\phi H$; in particular, it follows that $\N_H(\phi(O))$ is closed and the set $\{\phi(A) \mid A \in \Oc(G)\}$ is invariant under conjugation in $H$.  Since $G$ is second-countable, $\Oc(G)$ is countable, and hence $\N_H(\phi(O))$ has countable index.  By the Baire Category Theorem, it follows that $\N_H(\phi(O))$ is open in $H$.
\end{proof}

\subsection{Chief blocks}

We recall some definitions and properties of chief blocks from \cite{RW-Polish} and \cite{RW-DenseLC}.

\begin{defn}
A \defbold{(closed) normal factor} of a topological group $G$ is a pair of closed normal subgroups $(K,L)$ of $G$ such that $L \le K$.  We will generally write the pair as $K/L$ to emphasize the group structure of $K$ modulo $L$.  We say $K/L$ is a \defbold{chief factor} of $G$ if there are no closed normal subgroups of $G$ strictly between $L$ and $K$.

The \defbold{centralizer} $\CC_G(K/L)$ of a normal factor is the subgroup $\{g \in G \mid \forall k \in K: [g,k] \in L\}$.  Two nonabelian chief factors are \defbold{associated} if they have the same centralizer.  A \defbold{chief block} of $G$ is an association class of nonabelian chief factors.  Given a chief block $\mf{a}$, we write $\CC_G(\mf{a})$ for the centralizer of any representative of $\mf{a}$.  The quotient $G/\CC_G(\mf{a})$ then has a nontrivial closed normal subgroup $G^{\mf{a}}/\CC_G(\mf{a})$ representing $\mf{a}$, called the \defbold{uppermost representative} of $\mf{a}$.

Given a chief block $\mf{a}$ and a closed normal factor $K/L$, we say $K$ \defbold{covers} $\mf{a}$ if $L \le \CC_G(\mf{a})$ and $K \not\le \CC_G(\mf{a})$; otherwise, $K$ \defbold{covers} $\mf{a}$.  We say $K$ covers or avoids $\mf{a}$ according to whether $K/\triv$ covers or avoids $\mf{a}$.  A chief block $\mf{a}$ is \defbold{minimally covered} if, on taking the intersection $K$ of all closed normal subgroups of $G$ that cover $\mf{a}$, then $K$ covers $\mf{a}$.  In this case, we write $K = G_{\mf{a}}$.

Given a closed subgroup $H$ of $G$ and a chief block $\mf{a}$ of $H$, we say $\mf{a}$ \defbold{extends to $G$} if there is a chief block $\mf{b}$ such that for every closed normal subgroup $K$ of $G$, then $K$ covers $\mf{b}$ if and only if $K \cap H$ covers $\mf{a}$.  We then write $\mf{b} = \mf{a}^G$ and call $\mf{b}$ the \defbold{extension of $\mf{a}$}.
\end{defn}

Note that if $K/L$ is a representative of some chief block $\mf{a}$, then there is a normal compression from $K/L$ to the uppermost representative.  In particular, any property of topological groups that is invariant under normal compressions is an association invariant of chief blocks; in this case, we will say $\mf{a}$ has a property of topological groups, to mean its representatives all have that property.  Thus by Lemma~\ref{lem:compression_rank}, it makes sense to say a chief block $\mf{a}$ of a \tdlcsc group is elementary, and if so to define $\xi(\mf{a})$.  Similarly, by \cite[Theorem~5.8]{RW-DenseLC}, being quasi-discrete is an association invariant of chief blocks.  In particular, note that quasi-discrete chief blocks $\mf{a}$ are elementary with $\xi(\mf{a})=2$, by Lemma~\ref{lem:quasi-discrete}.  In a compactly generated \tdlcsc group, we refer to all chief blocks that are not quasi-discrete as \defbold{robust}.  In a general \tdlcsc group, a chief block is \defbold{regionally robust} if it is the extension of a robust chief block of a compactly generated open subgroup.

\begin{thm}[{See \cite[\S8]{RW-DenseLC}}]\label{thm:robust_block}
Let $G$ be a \tdlcsc group and let $H \in \Oc(G)$.
\begin{enumerate}[(i)]
\item Robust chief blocks of $H$ do not have any compact or discrete representatives.
\item There are only finitely many robust chief blocks of $H$, and hence only countably many regionally robust chief blocks of $G$.
\item Let $\mf{a}$ be a robust chief block of $H$.  Then $\mf{a}$ extends to $G$, and the extension $\mf{a}^G$ is minimally covered.
\item Every nonelementary chief block of $G$ is regionally robust.
\end{enumerate}
\end{thm}

\begin{proof}
Part (i) follows from the fact that all topologically characteristically simple profinite or discrete groups are quasi-discrete, together with Lemma~\ref{lem:quasi-discrete}.  For part (ii), the fact that $H$ has only finitely many robust chief blocks follows from \cite[Corollary~4.19]{RW-EC}, and then $G$ has only countably many regionally robust chief blocks because $\Oc(G)$ is a countable set. For part (iii): $\mf{a}$ extends to $G$ by \cite[Lemma~8.13]{RW-DenseLC}; $\mf{a}$ is minimally covered by \cite[Proposition~4.10]{RW-EC}; and then $\mf{a}^G$ is minimally covered by \cite[Lemma~8.6]{RW-DenseLC}.  Part (iv) follows from \cite[Corollary~8.18]{RW-DenseLC}.  
\end{proof}

In particular, we have the following.

\begin{cor}\label{cor:robust_block:countable}
Let $G$ be a \tdlcsc group.  Then $G$ has at most countably many nonelementary chief blocks, all of which are minimally covered.
\end{cor}

\subsection{Reduced envelopes and RIO subgroups}

To finish the preliminaries, we recall some results on reduced envelopes and RIO subgroups from \cite{ReidDistal}.

\begin{defn}
Given a \tdlc group $G$ and a subgroup $H$ of $G$, a \defbold{reduced envelope} for $H$ in $G$ is an open subgroup $E$ of $G$ such that $H \le E$, and whenever $|H:O \cap H|$ is finite, then $|E:O \cap E|$ is finite.  Write $E_G(H)$ for the set of open subgroups of $G$ commensurate with a reduced envelope of $H$.
\end{defn}

Reduced envelopes exist for every compactly generated subgroup.

\begin{thm}[{See \cite[Theorem~B]{ReidDistal}}]\label{thm:reduced_envelope}
Let $G$ be a \tdlc group and let $H$ be a compactly generated subgroup of $G$.  Then there is an open subgroup $E$ of $G$ with the following properties:
\begin{enumerate}[(i)]
\item\label{reduced_envelope:1} $E$ is a reduced envelope for $H$ in $G$;
\item\label{reduced_envelope:2} $E = H\Res{G}(H)U$, where $U$ is a compact open subgroup of $G$;
\item\label{reduced_envelope:3} $E$ is compactly generated;
\item\label{reduced_envelope:4} $\Res{G}(H) = \Res{}(E)$, so in particular, $\Res{G}(H)$ is normal in $E$.
\end{enumerate}
\end{thm}

In particular, given a \tdlc group $G$ and a compactly generated subgroup $H$ of $G$, then the reduced envelope of $H$ in $G$ is clearly unique up to finite index, so $E_G(H)$ is a well-defined element of $\Reg(G)$.  For any continuous homomorphism between \tdlc groups, there is then an induced map on the regional posets.  We will refer to the map $\theta$ defined by the next theorem as the \defbold{envelope map} of $\phi$.  In the case that $G$ is a closed subgroup of $H$ and $\phi$ is just inclusion, we refer to $\theta$ as the envelope map of $G$ in $H$.

\begin{thm}[{See \cite[Theorem~D]{ReidDistal}}]\label{thm:embedding_envelopes}
Let $G$ and $H$ be \tdlc groups and let $\phi: G \rightarrow H$ be a continuous homomorphism.  Let $K_1$ and $K_2$ be compactly generated subgroups of $G$ such that $E_G(K_1) = E_G(K_2)$.  Then $E_H(\phi(K_1)) = E_H(\phi(K_2))$. 

In particular, there is a well-defined map
\[
\theta: \Reg(G) \rightarrow \Reg(H); \; [O] \mapsto E_H(\phi(O)).
\]
If $H = \phi(G)X$ for some compact subset $X$ of $H$, then $\theta$ is surjective.
\end{thm}

Even when $\phi$ is an inclusion map, its envelope map can be far from injective, because a closed subgroup of a \tdlc group $H$ can have a complicated regional structure that is not witnessed by the open subgroups of $H$.  However, there is a class of closed subgroups that are specified by open subgroups in a certain sense, and for these closed subgroups, the envelope map will turn out to carry much more information.

\begin{defn}\label{def:RIO}
Let $G$ be a topological group and let $\mc{C}(G)$ be a class of closed subgroups of $G$.  Define $\mc{IC}(G)$ (the class of \defbold{intersection-$\mc{C}$} subgroups) to be the class of subgroups that are intersections of groups in $\mc{C}(G)$.  Define $\mc{RC}(G)$ (the class of \defbold{regionally-$\mc{C}$} subgroups) as follows: we have $H \in \mc{RC}(G)$ if there is a compactly generated open subgroup $U$ of $H$ such that, for every compactly generated group $K$ such that $U \le K \le H$, then $K \in \mc{C}(G)$.

Let $\mc{O}(G)$ denote the class of open subgroups of $G$.  If $H \in \mc{IO}(G)$, say $H$ is an \defbold{IO subgroup} of $G$, and if $H \in \mc{RIO}(G)$, say $H$ is a \defbold{RIO subgroup}.

Write $\mc{IO}_c(G)$ for the class of IO subgroups of $G$ that are compactly generated.
\end{defn}

Not all closed subgroups of a \tdlc group are RIO subgroups, but the RIO subgroups nevertheless form a large and robust class.

\begin{thm}[{\cite[Theorem~4.11]{ReidDistal}}]\label{thm:rio_closure_properties}
Let $G$ be a \tdlc group.  All of the following are RIO subgroups of $G$:
\begin{enumerate}[(i)]
\item\label{intro:rio_closure_properties:1} any closed subgroup $H$ that acts distally on $G/H$ by translation (for example, any closed subnormal subgroup of $G$);
\item\label{intro:rio_closure_properties:2} any RIO subgroup of a RIO subgroup of $G$;
\item\label{intro:rio_closure_properties:3} any intersection of RIO subgroups of $G$;
\item\label{intro:rio_closure_properties:4} the closure of any pointwise limit inferior of RIO subgroups of $G$;
\item\label{intro:rio_closure_properties:5} any closed subgroup $H$ of $G$, such that there is $K \le H \le G$ with $K \in \mc{RIO}(G)$ and $H/K$ compact.
\end{enumerate}
\end{thm}

Among compactly generated closed subgroups, those that are intersections of open subgroups can characterized in several ways.

\begin{lem}[{See \cite[Lemma~4.7]{ReidDistal}}]\label{lem:IO_equi_gen}
Let $G$ be a \tdlc group and let $H$ be a compactly generated closed subgroup of $G$.  Then the following are equivalent:
\begin{enumerate}[(i)]
\item $\Res{G}(H) = \Res{}(H)$;
\item $\Res{G}(H) \le H$;
\item $H \in \mc{RIO}(G)$;
\item $H \in \mc{IO}(G)$;
\item There is an open subgroup $E$ of $G$ such that $\Res{}(E) \le H \le E$.
\end{enumerate}
Moreover, if (i)--(v) hold, then in fact there is $E \in \Oc(G)$ such that $\Res{}(E) = \Res{}(H)$, $E = HU$ for a compact open subgroup $U$ of $E$, and the finite index open subgroups of $E$ that contain $H$ form a base of neighbourhoods of the trivial coset in $E/H$.\end{lem}

We recall here a standard fact about compactly generated \tdlc groups that will be useful later.

\begin{lem}\label{lem:cocompact_gen}
Let $G$ be a compactly generated \tdlc group, let $U$ be a compact open subgroup of $G$ and let $H$ be a symmetric subset, not necessarily closed, such that $G = HU$.  Then $G = \grp{S}U$ for a finite subset $S$ of $H$.
\end{lem}

\begin{proof}
Let $A$ be a compact symmetric generating set for $G$.  Then the product $UAU$ is compact, so there is a finite subset $S_1$ of $G$ such that $UAU = \bigcup_{s \in S_1}sU$.  Since $G = HU$ we are free to take $S_1 \subseteq H$.  Similarly, there is a finite subset $S_2$ of $G$ such that $UAU = \bigcup_{s \in S_2}Us$; by symmetry, we have $G = UH$, so we can take $S_2 \subseteq H$.  Now set $S$ to be the following finite subset of $H$:
\[
S = S_1 \cup S\inv_1 \cup S_2 \cup S\inv_2.
\]
Observe now that the following equations hold:
\[
SU = US = USU = UAU.
\]
From here on, the proof that $G = \grp{S}U$ is given by \cite[Proposition~4.1(ii)]{CRW-Part2}.
\end{proof}

Here are some more closure properties of the class of RIO subgroups that we will use later.

\begin{lem}\label{lem:RIO_quotient}
Let $G$ be a \tdlc group, let $N$ be a closed normal subgroup of $G$ and let $H \in \mc{RIO}(G)$.  Then $\ol{HN} \in \mc{RIO}(G)$ and $\ol{HN/N} \in \mc{RIO}(G/N)$.
\end{lem}

\begin{proof}
Let $K \in \Oc(\ol{HN})$.  By Lemma~\ref{lem:cocompact_gen}, we see that $K \le UK'N$, where $U$ is a compact open subgroup of $K$ and $K' \in \Oc(H)$.  Let $\mc{W}'$ be the set of open subgroups of $G$ that contain $K'$ as a cocompact subgroup.  Since $H \in \mc{RIO}(G)$, we see that $K' \in \mc{IO}_c(G)$, so by Lemma~\ref{lem:IO_equi_gen}, the set $\{W/K' \mid W \in \mc{W}'\}$ forms a base of neighbourhoods of the identity in $G/K'$.  Consequently, the set $\{WN/\ol{K'N} \mid W \in \mc{W}'\}$ is a base of neighbourhoods of the identity in $G/\ol{K'N}$, so $\ol{K'N} \in \mc{IO}(G)$; by Theorem~\ref{thm:rio_closure_properties}(v) it follows that $UK'N \in \mc{RIO}(G)$, and then since $K$ is a compactly generated open subgroup of $UK'N$, we have $K \in \mc{IO}(G)$.  Given the freedom of choice of $K$, we conclude that $\ol{HN} \in \mc{RIO}(G)$.

Now consider $L/N \in \Oc(\ol{HN}/N)$.  Then $L = KN$ for some $K \in \Oc(\ol{HN})$.  Similar to the previous paragraph, letting $\mc{W}$ be the set of open subgroups of $G$ that contain $K$ as a cocompact subgroup, then $\{W/K \mid W \in \mc{W}\}$ is a base of neighbourhoods of the identity in $G/K$ and hence $\{WL/L \mid W \in \mc{W}\}$ is a base of neighbourhoods of the identity in $G/L$, so $L \in \mc{IO}(G)$ and hence $L/N \in \mc{IO}(G/N)$.  We conclude that $\ol{HN/N} \in \mc{RIO}(G/N)$.
\end{proof}

\begin{lem}\label{lem:RIO_sigma-normal}
Let $G$ and $H$ be first-countable \tdlc groups, let $\phi: G \rightarrow H$ be a regionally normal map such that $\ol{\phi(G)} \in \mc{RIO}(H)$, and let $K$ be a RIO subgroup of $G$.
\begin{enumerate}[(i)]
\item We have $\ol{\phi(K)} \in \mc{RIO}(H)$.
\item Suppose that $K$ is compactly generated.  Then for $2 \le \alpha \le \omega$, we have $\ol{\phi(\Res{\alpha}(K))} = \Res{\alpha}(\ol{\phi(K)})$.  In particular, $\Res{}(\ol{\phi(K)}) = \ol{\phi(\Res{}(K))}$.
\end{enumerate}
\end{lem}

\begin{proof}
Since the RIO subgroup relation is transitive, without loss of generality we may replace $H$ with $\ol{\phi(G)}$, and so assume that $\phi$ has dense image.

We can write $K = \bigcup_{i \in I}K_i$, where $(K_i)_{i \in I}$ is an ascending net of compactly generated open subgroups of $K$.  It follows that 
\[
\ol{\phi(K)} = \ol{\bigcup_{i \in I}\phi(K_i)} = \ol{\liminf_{i \in I}(\ol{\phi(K_i)})}.
\]
Hence by Theorem~\ref{thm:rio_closure_properties}(iv), to show $\ol{\phi(K)} \in \mc{RIO}(H)$, it suffices to show that $\ol{\phi(K_i)} \in \mc{RIO}(H)$ for all $i \in I$.  Thus for the rest of the proof, we may assume that $K$ is compactly generated.  Fix a reduced envelope $E$ for $K$ in $G$ and note that $K$ is cocompact in $E$.

We have $E \in \Oc(G)$, so by the regionally normal property, $\N_H(\phi(E))$ is open in $H$.  In particular, $\ol{\phi(E)}$ is a RIO subgroup of $H$, so to show $\ol{\phi(K)} \in \mc{RIO}(H)$, it is enough to show $\ol{\phi(K)} \in \mc{RIO}(\ol{\phi(E)})$.  Thus without loss of generality, we may replace $G$ with $E$ and $H$ with $\ol{\phi(E)}$.  Under these hypotheses, $G$ and $H$ are both second-countable and hence Polish.  By Lemma~\ref{lem:RIO_quotient}, $\ol{K\ker\phi}/\ker\phi$ is a RIO subgroup of $G/\ker\phi$.  Thus by replacing $G$ with $G/\ker\phi$, we may assume that $\phi$ is injective, and hence a normal compression of \tdlc Polish groups.

By Theorem~\ref{thm:psi-compression_factor_rel}, we can now factorize $\phi$ as $\phi = \pi \theta$, where $\theta: G \rightarrow G \rtimes H$ is a closed embedding and $\pi: G \rtimes H \rightarrow H$ is a quotient map.  Given $K \in \mc{RIO}(G)$, we then have $\theta(K) \in \mc{RIO}(G \rtimes H)$ by Theorem~\ref{thm:rio_closure_properties}(ii), and then $\ol{\phi(K)} \in \mc{RIO}(H)$ by Lemma~\ref{lem:RIO_quotient}.  This completes the proof of (i).
 
For (ii), we assume that $K$ is compactly generated.  Recall that in the context of compactly generated \tdlc groups, the discrete residual is the same as the elementary rank-$2$ residual.  By Lemma~\ref{lem:IO_equi_gen} we have $\Res{}(E) = \Res{}(K)$.  Using the elementary rank-$2$ series we then see that $\Res{\alpha}(E) = \Res{\alpha}(K)$ for $2 \le \alpha \le \omega$.  Similarly, considering $\ol{\phi(K)}$ as a subgroup of $\ol{\phi(E)}$, we see from part (i) that $\ol{\phi(K)}$ is an intersection of open subgroups of $H$ and hence also of $\ol{\phi(E)}$; since $\ol{\phi(K)}$ is clearly cocompact in $\ol{\phi(E)}$, applying Lemma~\ref{lem:IO_equi_gen} shows that $\Res{}(\ol{\phi(E)}) = \Res{}(\ol{\phi(K)})$, and hence also $\Res{\alpha}(\ol{\phi(E)}) = \Res{\alpha}(\ol{\phi(K)})$ for $2 \le \alpha \le \omega$.  Thus without loss of generality we may assume $K=E$, that is, $K$ is open in $G$.  In that case, $\N_H(\phi(K))$ is open in $H$.  In particular, $\phi(K)$ is normal in its closure.

We consider elementary quotients of $K$ and of $\ol{\phi(K)}$ of rank at most some ordinal $\alpha$ where $2 \le \alpha \le \omega$; let $R_{\alpha} = \Res{\alpha}(\ol{\phi(K)})$.  Note that $\ol{\phi(\Res{\alpha}(K))}$ is normal in $\ol{\phi(K)}$ by Lemma~\ref{lem:compression_normal_subgroup}.  We then have a dense regionally normal map from $K/\Res{\alpha}(K)$ to $\ol{\phi(K)}/\ol{\phi(\Res{\alpha}(K))}$; by Lemmas~\ref{lem:finite_rank_residual}(i) and~\ref{lem:compression_rank}, it follows that $\xi(\ol{\phi(K)}/\ol{\phi(\Res{\alpha}(K))}) \le \alpha$, and hence $R_{\alpha} \le \ol{\phi(\Res{\alpha}(K))}$.  On the other hand, we have a normal compression from $K/\phi\inv(R_{\alpha})$ to $\ol{\phi(K)}/R_{\alpha}$, and hence $\xi(K/\phi\inv(R_{\alpha})) \le \alpha$ by Lemma~\ref{lem:finite_rank_residual}(i) and Lemma~\ref{lem:compression_rank}, showing that $\phi\inv(R_{\alpha}) \ge \Res{\alpha}(K)$, or in other words $\phi(\Res{\alpha}(K)) \le R_{\alpha}$.  Hence $R_{\alpha} = \ol{\phi(\Res{\alpha}(K))}$.  In particular, in the case $\alpha=2$ we have $\Res{}(\ol{\phi(K)}) = \ol{\phi(\Res{}(K))}$.
\end{proof}

\section{A partial order on the regional poset}\label{sec:ll}

\subsection{Definition of the strong ordering}

In this section we set up the partial order on $\Reg(G)$ that will be used to define the class of $\mc{S}$-well-founded groups, where $\mc{S}$ is a class of \tdlcsc groups with property $\propS$.

\begin{defn}\label{def:ll}
Let $\mc{S}$ be a class of \tdlcsc groups with property $\propS$.  We define the \defbold{strong ordering} on compactly generated \tdlcsc groups (with respect to $\mc{S}$) as follows: $H \ll_{\mc{S}} K$ if $H$ is a RIO subgroup of $K$, $K$ is noncompact, and the following two conditions are met:
\begin{enumerate}[(a)]
\item $\ol{H\Res{\omega}(K)}/\Res{\omega}(K)$ is compact;
\item For all closed normal subgroups $N$ of $\Res{\omega}(K)$ such that $\Res{\omega}(K)/N \in \mc{S} \smallsetminus \Es{}$ and $|K:\N_K(N)|<\infty$, then $\Res{\omega}(H) \le N$.
\end{enumerate}
\end{defn}

Condition (a) is based on the characterization of regionally elementary groups given by Proposition~\ref{prop:elementary_wellfounded}.  Note that condition (b) is vacuously satisfied when $K$ is elementary, since in that case there cannot be any quotient of a subgroup of $K$ that belongs to $\mc{S} \smallsetminus \Es{}$.  Condition (b) is motivated by Theorem~\ref{thm:CM}.  See Remark~\ref{rmk:omega} below for why the finite elementary rank residual is used here instead of the discrete residual.

For the rest of this section, we fix the class $\mc{S}$ (the reader is safe to assume $\mc{S}$ is any of $\Sclass$, $\ms{R} \cap \Bigsclass$, or $\Bigsclass$) and write $\ll$ for $\ll_{\mc{S}}$.

The next lemma places some useful restrictions on the quotients under consideration in condition (b).

\begin{lem}\label{lem:S_quotient_count}
Let $G$ be a compactly generated \tdlcsc group.  Let $\mc{L}$ be the set of closed normal subgroups $M$ of $\Res{\omega}(G)$ such that $\Res{\omega}(G)/M \in \mc{S} \smallsetminus \Es{}$, and let $\mc{M}$ be the set of $M \in \mc{L}$ such that $|G:\N_G(M)|<\infty$.
\begin{enumerate}[(i)]
\item The set $\mc{L}$ is countable.  Moreover, for each $M \in \mc{L}$ the normalizer of $M$ in $G$ is open.
\item There are only finitely many orbits on $\mc{L}$ under conjugation; in particular, $\mc{M}$ is finite.
\item If $G/\Res{}(G)$ is compact but $G$ is not compact, then $\Res{}(G) = \Res{\omega}(G)$ and $\mc{L} = \mc{M} \neq \emptyset$.
\end{enumerate}
\end{lem}

\begin{proof}
Let $N = \bigcap_{M' \in \mc{L}}M'$.  Given $M \in \mc{L}$, we see that $\Res{\omega}(G)/M$ is a nonelementary chief factor of $\Res{\omega}(G)/N$, with centralizer $M/N$, representing the chief block $\mf{a}_M$ of $\Res{\omega}(G)/N$.

For (i), we see by Corollary~\ref{cor:robust_block:countable} that $\mc{L}$ is countable.  In particular, the normalizer of $M$ has countable index in $G$; since normalizers of closed subgroups are closed, it follows by the Baire Category Theorem that $\N_G(M)$ is open in $G$.

For (ii), consider again $M \in \mc{L}$.  By Corollary~\ref{cor:robust_block:countable} the chief block $\mf{a}_M$ is minimally covered, that is, there is a unique smallest closed normal subgroup $K_M/N$ of $\Res{\omega}(G)/N$ such that $K_M \nleq M$.  Given $M' \in \mc{L} \smallsetminus \{M\}$, we see that $M' \nleq M$, so $K_M \leq M'$; in particular, it follows that $K_M \cap M = N$, so in fact $K_M$ is minimal among nontrivial closed normal subgroups of $G/N$ and $K_M/N$ is the lowermost representative of $\mf{a}_M$.  We also see that $K_M \cap K_{M'} = N$, and hence also $[K_M,K_{M'}]\le N$, for all $M' \in \mc{L} \smallsetminus \{M\}$.  Since $K_M/N$ is a representative of $\mf{a}_M$ we have a normal compression from $K_M/N$ to $\Res{\omega}(G)/M$.  Clearly $K_M/N$ is nonabelian, so by minimality it must be topologically perfect; it follows by property $\propS$ that in fact $K_M/N \in \mc{S} \smallsetminus \ms{E}{}$.  We now see that $\cgrp{gK_Mg\inv/N \mid g \in G}$ is a quasi-product of copies of a group in $\mc{S}$, giving rise to a chief block $\mf{b}_M$ of $G$; this chief block is robust since $K_M/N \not\in \ms{E}{}$ .  Since $K_M/N$ and $K_{M'}/N$ commute for any distinct $M,M' \in \mc{M}$, it is easy to see that $\mf{b}_M = \mf{b}_{M'}$ if and only if $K_M$ is conjugate in $G$ to $K_{M'}$; since $M$ is the centralizer in $\Res{\omega}(G)$ of $K_M/N$, this can happen only if $M$ is conjugate in $G$ to $M'$.  By Theorem~\ref{thm:robust_block}, there are only finitely many robust chief blocks of $G$, so there are only finitely many $G$-conjugacy classes in $\mc{L}$.  Now $\mc{M}$ is the union of finitely many finite orbits of $G$, so it is finite.

For (iii), suppose that $K = \Res{}(G)$ is cocompact in $G$ and $G$ is not compact.  Then $K$ has the same discrete residual acting on $G$ as $G$ itself has, that is, $\Res{G}(K) = \Res{}(G)$; since $K$ is normal, we have $\Res{G}(K) = \Res{}(K)$.  Thus the elementary rank-$2$ series of $G$ terminates at $\Res{}(G)$, that is, $\Res{}(G) = \Res{\omega}(G)$.

Under the given hypotheses, we see by Theorem~\ref{thm:CM} that $\Res{}(G)$ has $n$ quotients
\[
\Res{}(G)/M_1,\dots,\Res{}(G)/M_n
\]
in $\Sclass$, where $0 < n < \infty$.  Moreover, in this situation we see that $\Res{}(G)$ is compactly generated, so all of its quotients are compactly generated, and hence every quotient in $\mc{S}$ is in fact in $\Sclass$.  We then have $|\mc{L}|=n$, and it is clear that $\mc{L} = \mc{M}$.
\end{proof}

\subsection{Comparing RIO subgroups}

We will mostly be using $\ll$ to compare elements of $\Oc(G)$ or $\Reg(G)$ for a given first-countable \tdlc group $G$.  We first establish that $\ll$ defines a strict partial order in both cases.

\begin{lem}\label{lem:strong_order:basic}
Let $G$ be a first-countable \tdlc group and let $H,K \in \Oc(G)$ such that $H \ll K$.  Then $H$ has infinite index in $K$.  Moreover, we have $H' \ll K'$ whenever $K' \in [K]_f$ and $H' \in \Oc(G)$ is a subgroup of $K'$ such that $[H']_f \le [H]_f$.
\end{lem}

\begin{proof}Suppose for a contradiction that $H \ll K$ and $H$ has finite index in $K$.  In this case, $\Res{}(H) = \Res{}(K)$ and hence also $\Res{\omega}(H) = \Res{\omega}(K)$.  Since $H/\Res{\omega}(K)$ is compact, we see that $K/\Res{}(K)$ is compact, while $K$ is not compact.  By Lemma~\ref{lem:S_quotient_count}(iii), there is a quotient map $\pi: \Res{\omega}(K) \rightarrow S$ to some $S \in \mc{S} \smallsetminus \Es{}$, such that $|K:\N_K(\ker\pi)|<\infty$; we then have $\Res{\omega}(H) \le \ker\pi < \Res{\omega}(K)$, a contradiction.

The last conclusion follows by observing how the definition of $\ll$ is expressed in terms of finite elementary rank residuals, and noting that if $L,L' \in \Oc(G)$ are such that $L$ virtually contains $L'$, then $\Res{}(L) \ge \Res{}(L')$ and hence $\Res{\omega}(L) \ge \Res{\omega}(L')$.\end{proof}

For $[H],[K] \in \Reg(G)$, we write $[H] \ll [K]$ if $H \ll K$ for some representatives $H$ and $K$.

From the perspective of regional theory, the natural extension of $\ll$ to $\mc{RIO}(G)$ is as follows: given $H,K \in \mc{RIO}(G)$, we say $H \ll K$ if 
\begin{equation}\label{eq:RIO_ll}
\forall H' \in \Oc(H) \; \exists K' \in \Oc(K) \; \forall K'' \in \Oc(K): \; K' \le K'' \Rightarrow H' \ll K''.
\end{equation}
Lemma~\ref{lem:strong_order:basic} ensures that this definition agrees with Definition~\ref{def:ll} in the case that both $H$ and $K$ are compactly generated.  However, in practice we will mostly be able to focus on the poset $(\Reg(G),\ll)$, because the latter poset effectively captures the partial order on compactly generated IO subgroups.  This will be shown in the next two results.

\begin{prop}\label{prop:ll_RIO}
Let $G$ be a first-countable \tdlc group and let $H_1$ and $H_2$ be compactly generated IO subgroups of $G$ such that $H_1 \le H_2$.  Let $H'_1$ be a reduced envelope of $H_1$ in $H_2$; let $H''_2$ be a reduced envelope of $H_2$ in $G$; and let $H''_1$ be a reduced envelope of $H_1$ in $H''_2$.  Then the following are equivalent:
\begin{enumerate}[(i)]
\item $H_1 \ll H_2$;
\item $H'_1 \ll H_2$;
\item $H_1 \ll H''_2$;
\item $H'_1 \ll H''_2$;
\item $H''_1 \ll H''_2$.
\end{enumerate}
\end{prop}

\begin{proof}
Let $H'_2 = H_2$ and let $R_i = \Res{\omega}(H_i)$.  We note that $\Res{}(H_i) = \Res{}(H'_i) = \Res{}(H''_i)$, so $R_i = \Res{\omega}(H'_i) = \Res{\omega}(H''_i)$ for $i=1,2$.  Without loss of generality we can choose $H''_1$ to contain $H'_1$.  Then $H_i$ is cocompact in $H'_i$, while $H'_i$ is cocompact in $H''_i$, for $i=1,2$.  It follows that all five comparisons (i)--(v) are equivalent as far as condition (a) is concerned.  For condition (b), it suffices to show the following: given $N \unlhd R_2$ such that $R_2/N \in \mc{S} \smallsetminus \Es{}$, then $|H''_2:\N_{H''_2}(N)| < \infty$ if and only if $|H_2:\N_{H_2}(N)| <\infty$.  If $|H''_2:\N_{H''_2}(N)| < \infty$ then clearly $|H_2:\N_{H_2}(N)| <\infty$.  Conversely, suppose that $|H_2:\N_{H_2}(N)| <\infty$.  Then by Lemma~\ref{lem:S_quotient_count}, $\N_{H''_2}(N)$ is open in $H''_2$.  Now $H_2$, and hence its finite index subgroup $\N_{H_2}(N)$, is cocompact in $H''_2$; thus $\N_{H''_2}(N)$ is both open and cocompact in $H''_2$, so it is of finite index.  We now see that condition (b) is equivalent for the comparisons (i)--(v), and hence that (i)--(v) are equivalent.
\end{proof}

\begin{cor}\label{cor:res_RIO}
Let $G$ be a first-countable \tdlc group and let $H$ be a RIO subgroup of $G$.  Then the envelope map of $H$ in $G$ is a $(\le,\ll)$-order embedding.  If $G = HX$ for some compact $X \subseteq G$, then the envelope map is a $(\le,\ll)$-order isomorphism.
\end{cor}

\begin{proof}
Let $\theta: \Reg(H) \rightarrow \Reg(G)$ be the envelope map.  Let $K_1, K_2 \in \Oc(H)$ and let $O_i$ be a reduced envelope for $K_i$ in $G$.  By \cite[Theorem~F]{ReidDistal} we can choose $O_i$ so that $K_i = H \cap O_i$.  If $O_2$ virtually contains $O_1$ then it follows immediately that $K_2$ virtually contains $K_1$.  Conversely if $K_2$ virtually contains $K_1$, then $O_2$ virtually contains $K_1$; since $O_1$ is a reduced envelope of $K_1$, it follows that $O_2$ virtually contains $O_1$.  So $\theta$ is an order embedding for $\le$.  The fact that $\theta$ is a $\ll$-order embedding is clear from Lemma~\ref{lem:strong_order:basic} and Proposition~\ref{prop:ll_RIO}.

If $G = HX$ for some compact subset $X$ of $G$, then $\theta$ is surjective by Theorem~\ref{thm:embedding_envelopes}, and hence $\theta$ is a $(\le,\ll)$-order isomorphism in this case.
\end{proof}

For Corollary~\ref{cor:res_RIO}, it is critical that $H$ be a RIO subgroup of $G$.  We will see later (Remark~\ref{rmk:bad_rank}) that the following can happen: $G$ is a compactly generated \tdlcsc group, $H$ is a closed cocompact subgroup of $G$, $|\Reg(G)|=2$, but $(\Reg(H),\ll)$ is not even well-founded.

\subsection{Quotient maps}

A map $\phi: (A,\le) \rightarrow (B,\le)$ of partially ordered sets is a \defbold{quotient map} if $\phi$ is surjective, weakly order-preserving, and whenever $b_1 \le b_2$ in $B$, there are $a_i \in A$ such that $a_1 \le a_2$ and $\phi(a_i) = b_i$.  In this subsection we obtain a sufficient condition for the envelope map of a homomorphism of \tdlc groups to be a $(\le,\ll)$-quotient map.

To keep track of condition (b) in the definition of $\ll$, we will need to understand how normal factors in $\mc{S}$ behave under dense normal maps of compactly generated \tdlcsc groups.

\begin{lem}\label{lem:S_compression}
Let $G$ and $H$ be compactly generated \tdlcsc groups and let $\psi: G \rightarrow H$ be a continuous homomorphism with dense normal image.
\begin{enumerate}[(i)]
\item Let $M$ be a closed $G$-invariant subgroup of $\Res{\omega}(G)$ such that $\Res{\omega}(G)/M \in \mc{S} \smallsetminus \Es{}$.  Then $\ol{\psi(M)}$ is an $H$-invariant subgroup of $\Res{\omega}(H)$.  If $\psi(M)$ is not dense in $\Res{\omega}(H)$, then $\Res{\omega}(H)/\ol{\psi(M)}$ is a topologically perfect central extension of a group in $\mc{S} \smallsetminus \Es{}$, and writing $Z/\ol{\psi(M)} = \Z(\Res{\omega}(H)/\ol{\psi(M)})$, then the $\psi$-preimage of $Z$ in $\Res{\omega}(G)$ is exactly $M$.
\item Let $N$ be a closed $H$-invariant subgroup of $\Res{\omega}(H)$ such that $\Res{\omega}(H)/N \in \mc{S} \smallsetminus \Es{}$.  Then $M := \psi\inv(N) \cap \Res{\omega}(G)$ is a $G$-invariant closed normal subgroup of $\Res{\omega}(G)$ such that $\Res{\omega}(G)/M \in \mc{S} \smallsetminus \Es{}$.
\end{enumerate}
\end{lem}

\begin{proof}
For the proof it is convenient to factorize $\psi$ as $\psi = \phi\pi$, where $\pi: G \rightarrow G/\ker\psi$ is a quotient map and $\phi: G/\ker\psi \rightarrow H$ is a normal compression.  Note that by Lemma~\ref{lem:RIO_sigma-normal}(ii) we have
\[
\Res{\omega}(\pi(G)) = \ol{\pi(\Res{\omega}(G))} \text{ and } \Res{\omega}(H) = \ol{\psi(\Res{\omega}(G))} =  \ol{\phi(\Res{\omega}(\pi(G)))};
\]
by Lemma~\ref{lem:finite_rank_residual}, all the groups in the line above are topologically perfect.

(i)
The fact that $\ol{\psi(M)}$ is normal in $H$ follows from Lemma~\ref{lem:compression_normal_subgroup}.  Since $\Res{\omega}(H) = \ol{\psi(\Res{\omega}(G))}$, we have $\ol{\psi(M)} \le \Res{\omega}(H)$.

For the rest of the proof of this part, we may suppose that $\ol{\psi(M)} < \Res{\omega}(H)$.  Since $\Res{\omega}(G)$ is topologically perfect, it follows that $\ol{\psi(M)}$ does not contain the derived group of $\psi(\Res{\omega}(G))$.  In particular, $\Res{\omega}(H)/\ol{\psi(M)}$ is not abelian, so $Z < \Res{\omega}(H)$, where $Z/\ol{\psi(M)} = \Z(\Res{\omega}(H)/\ol{\psi(M)})$.  We see that
\[
M = \psi\inv(\ol{\psi(M)}) \cap \Res{\omega}(G) \le \psi\inv(Z) \cap \Res{\omega}(G) < \Res{\omega}(G);
\]
since all the terms in the above inequalities are closed normal subgroups of $\Res{\omega}(G)$ and $\Res{\omega}(G)/M$ is topologically simple, in fact
\[
M = \psi\inv(\ol{\psi(M)}) \cap \Res{\omega}(G) = \psi\inv(Z) \cap \Res{\omega}(G).
\]
We have a normal compression from $\Res{\omega}(G)/M$ to $\Res{\omega}(H)/\ol{\psi(M)}$; the latter group is therefore a topologically perfect central extension of a group in $\mc{S} \smallsetminus \Es{}$, by property $\propS$.

(ii)
Let $M' = \pi\inv(\Res{\omega}(\pi(G)) \cap \phi\inv(N))$.  Then $\phi$ restricts to a continuous injective map $\phi'$ from $\Res{\omega}(\pi(G))/\pi(M')$ to $\Res{\omega}(H)/N$ with dense image; in fact $\phi'$ is a normal compression by Lemma~\ref{lem:compression_normal_subgroup}.  By property $\propS$ we deduce that $\Res{\omega}(\pi(G))/\pi(M') \in \mc{S} \smallsetminus \Es{}$.  We then have a natural isomorphism of topological groups between $\Res{\omega}(\pi(G))/\pi(M')$ and $\Res{\omega}(G/M')$, so $\Res{\omega}(G/M') \in \mc{S} \smallsetminus \Es{}$.

Now let $M = \psi\inv(N) \cap \Res{\omega}(G)$; notice that $M = M' \cap \Res{\omega}(G)$.  We then have a normal compression from $\Res{\omega}(G)/M$ to $\ol{\Res{\omega}(G)M'}/M'$; applying Lemma~\ref{lem:RIO_sigma-normal}(ii), in fact $\ol{\Res{\omega}(G)M'}/M' = \Res{\omega}(G/M')$.  Since $\Res{\omega}(G)$ is topologically perfect, we conclude that $\Res{\omega}(G)/M \in \mc{S} \smallsetminus \Es{}$ via property $\propS$.
\end{proof}

We come now to the main theorem of this subsection, ensuring that the partial orders $\le$ and $\ll$ are well-behaved under the envelope map of a dense regionally normal map.

\begin{thm}\label{thm:reg_quotient}
Let $G$ and $H$ be first-countable \tdlc groups, let $\phi: G \rightarrow H$ be a dense regionally normal map, and let $\theta$ be the envelope map of $\phi$.  Then $\theta$ is a $(\le,\ll)$-quotient map from $\Reg(G)$ to $\Reg(H)$.  Moreover, the following holds:
\begin{enumerate}[(i)]
\item Given $O_1,O_2 \in \Oc(G)$ such that $O_1 \ll O_2$, then either $\ol{\phi(O_2)}$ is compact (implying that $\theta[O_1] = \theta[O_2]$), or $\ol{\phi(O_1)} \ll \ol{\phi(O_2)}$.
\item Given $E_1,E_2 \in \Oc(H)$ such that $E_1 \ll E_2$, and $O_2$ in $\Oc(G)$ such that $E_2$ is a reduced envelope of $\phi(O_2)$, there exists $O_1 \in \Oc(G)$ such that $O_1 \ll O_2$ and $E_1$ is a reduced envelope of $\phi(O_1)$.
\end{enumerate}
\end{thm}

\begin{proof}
Note first that $\phi$ has dense image in $H$, so $H = \phi(G)U$ for any compact open subgroup $V$ of $H$.  Consequently $\theta$ is surjective by Theorem~\ref{thm:embedding_envelopes}.  It is also clear that $\theta$ is weakly $\le$-order preserving.  By Lemma~\ref{lem:RIO_sigma-normal}(i), given $K \in \mc{RIO}(G)$, then $\ol{\phi(K)} \in \mc{RIO}(H)$; in particular, given $O \in \Oc(G)$, then $\ol{\phi(O)} \in \mc{IO}_c(H)$.

Consider now $E_1,E_2 \in \Oc(H)$ such that $E_1 \le E_2$; we aim to show that the $\le$-ordering of $[E_1]$ and $[E_2]$ is induced by a $\le$-ordered pair of preimages in $\Reg(G)$.  Let $O_i \in \Oc(G)$ be such that $\theta[O_i] = [E_i]$; after passing to subgroups of finite index we may ensure $\phi(O_i) \le E_i$.  We see that $E_2 = \grp{E_1,E_2}$ is also a reduced envelope for $\phi(O'_2)$, where $O'_2 = \grp{O_1,O_2}$.  Let $L_2 = \ol{\phi(O'_2)}$, and let $L_1$ be a reduced envelope of $\phi(O_1)$ as a subgroup of $L_2$.  We now have a dense homomorphism from $O'_2$ to $L_2$, and $L_1 \in \Oc(L_2)$.  Applying Theorem~\ref{thm:embedding_envelopes}, there is $O'_1 \in \Oc(O'_2)$ such that $L_1$ is commensurate with a reduced envelope of $\phi(O'_1)$ in $L_2$.  We now have $O'_1,O'_2 \in \Oc(G)$ such that $O'_1 \le O'_2$, and we see that $\theta[O'_2] = [E_2]$ and
\[
\theta[O'_1] = E_H(\phi(O'_1)) = E_H(L_1) = E_H(\phi(O_1)) = [E_1].
\]
This completes the proof that $\theta$ is a $\le$-quotient map.

For the strong ordering, suppose that we have $O_1, O_2 \in \Oc(G)$ such that $O_1 \ll O_2$.  If $\ol{\phi(O_2)}$ is compact, then $\theta[O_2]$ is the smallest element $[W]$ of $\Reg(H)$, where $W$ is any compact open subgroup of $H$; also $\theta[O_1] = [W]$, so $\theta[O_1] = \theta[O_2]$.  From now on we may assume that $\ol{\phi(O_2)}$ is not compact.

Since $O_1 \ll O_2$ we see that $O_1 \subseteq Y\Res{\omega}(O_2)$ for some compact open subset $Y$ of $O_1$.  Let $R_i = \Res{\omega}(\ol{\phi(O_i)})$.  By Lemma~\ref{lem:RIO_sigma-normal}(ii), we have $R_i = \ol{\phi(\Res{\omega}(O_i))}$; in particular, $\ol{\phi(O_1)} \subseteq \ol{\phi(Y)}R_2$, and $\ol{\phi(Y)}$ is a compact subset of $\ol{\phi(O_2)}$.  Hence condition (a) is satisfied by the pair $(\ol{\phi(O_1)}, \ol{\phi(O_2)})$.  In addition, we note that $O_2$ and $\ol{\phi(O_2)}$ are second-countable and hence Polish, and $\phi$ restricts to a dense regionally normal map from $O_2$ to $\ol{\phi(O_2)}$.

We now consider condition (b) of the definition of the $\ll$-ordering as it relates to $\ol{\phi(O_1)}$ and $\ol{\phi(O_2)}$.  Let $N$ be a closed normal subgroup of $R_2$ such that $R_2/N \in \mc{S} \smallsetminus \Es{}$ and such that $N$ has finitely many conjugates in $\ol{\phi(O_2)}$.  By passing to the normalizer of $N$ in $\ol{\phi(O_2)}$, and replacing $O_1$ and $O_2$ with their intersections with the preimage of the normalizer, we may assume that $N$ is normal in $\ol{\phi(O_2)}$; note that this is a change of finite index, so it will not affect $\ll$, nor will it change $R_i$ or $\Res{\omega}(O_i)$.  Applying Lemma~\ref{lem:S_compression}(ii), we then see that $M = \phi\inv(N) \cap \Res{\omega}(O_2)$ is an $O_2$-invariant normal subgroup of $\Res{\omega}(O_2)$ such that $\Res{\omega}(O_2)/M \in \mc{S} \smallsetminus \Es{}$.  In particular, since $O_1 \ll O_2$, we have $\Res{\omega}(O_1) \le M$, so $\phi(\Res{\omega}(O_1)) \le N$, and hence $R_1 \le N$.  Given the freedom of choice of $N$, we conclude that $\ol{\phi(O_1)} \ll \ol{\phi(O_2)}$.  This completes the proof of the claim (i).  We also deduce via Corollary~\ref{cor:res_RIO} that $\theta[O_1] \ll \theta[O_2]$.

\

We now move onto the assertion (ii), which will finish the proof of the theorem.

Let $E_1,E_2 \in \Oc(H)$ such that $E_1 \ll E_2$ and let $R_i = \Res{\omega}(E_i)$.  Choose some $O_2 \in \Oc(G)$ such that $E_2$ is a reduced envelope of $\phi(O_2)$; note that $E_2$ cannot be compact, and hence $O_2$ cannot be compact.  The map $\phi$ restricts to a dense regionally normal map from $O_2$ to $\ol{\phi(O_2)}$, and the latter group is locally normal (in particular, RIO) in $E_2$; since $E_2$ is also a reduced envelope, in fact $\ol{\phi(O_2)}$ is cocompact in $E_2$ and $R_2 = \Res{\omega}(\ol{\phi(O_2)})$.  In turn, $\Res{\omega}(\ol{\phi(O_2)}) = \ol{\phi(\Res{\omega}(O_2))}$ by Lemma~\ref{lem:RIO_sigma-normal}(ii).

Let $\mc{M}$ be the set of closed normal subgroups $M$ of $\Res{\omega}(O_2)$ such that $\Res{\omega}(O_2)/M \in \mc{S} \smallsetminus \Es{}$ and $\N_{O_2}(M)$ has finite index in $O_2$.  By Lemma~\ref{lem:S_quotient_count}, $\mc{M}$ is finite, so by replacing $O_2$ with a finite index open subgroup, we may assume that each $M \in \mc{M}$ is normal in $O_2$.  We partition $\mc{M}$ into two sets $\mc{M}_1$ and $\mc{M}_2$, where $M \in \mc{M}_1$ if $\phi(M)$ is dense in $R_2$ and $M \in \mc{M}_2$ otherwise.

Let $M \in \mc{M}_2$.  Then Lemma~\ref{lem:S_compression}(i) ensures that $R_2/\ol{\phi(M)}$ has a quotient in $\mc{S} \smallsetminus \Es{}$, with central kernel $C_M/\ol{\phi(M)}$ say, such that $\phi\inv(C_M) \cap \Res{\omega}(O_2) = M$.  Since $E_1 \ll E_2$ we then have $R_1 \le C_M$ and hence $\phi\inv(R_1) \cap \Res{\omega}(O_2) \le M$.  We have now shown that
\begin{equation}\label{eq:res_component}
\forall M \in \mc{M}_2: \; \phi\inv(R_1) \cap \Res{\omega}(O_2) \le M.
\end{equation}

\emph{Claim: If $\phi$ is injective, then $\mc{M} = \mc{M}_2$.}

Suppose $\phi$ is injective and let $M \in \mc{M}$.  Since we have a normal compression from $O_2$ to $\ol{\phi(O_2)}$, it follows by Lemma~\ref{lem:compression_normal_subgroup} that $\phi(M)$ is normal in $\ol{\phi(O_2)}$, and in particular $\phi(M)$ is normal in $R_2$.  Since $M$ does not contain the derived group of $\Res{\omega}(O_2)$, we see by Lemma~\ref{lem:normal_compression:derived} that $\phi(M)$ is not dense in $R_2$.  This proves the claim.

\

To finish the proof of (ii), we will prove it in some special cases that will lead to a general solution.  In particular, we are free to factorize $\phi$ as some product $\phi = \phi_n\phi_{n-1}\dots \phi_1$ and then prove (ii) for each map $\phi_i$ in turn.  Let $K$ be the kernel of $\phi$.

\emph{Case 1: $\mc{M} = \mc{M}_2$.}

Let $U$ be a compact open subgroup of $O_2$ and let $A = \phi\inv(E_1) \cap (\Res{\omega}(O_2)U)$.  We see that $\phi(A) = E_1 \cap \phi(\Res{\omega}(O_2)U)$; taking the closure on either side, and using the fact that $E_1$ is clopen, it follows that $\ol{\phi(A)} = E_1 \cap R_2\ol{\phi(U)}$.  In particular, $\ol{\phi(A)}$ contains $E_1 \cap R_2$; since $E_1 \ll E_2$ it follows that $\ol{\phi(A)}$ is cocompact in $E_1$.  By Theorem~\ref{thm:embedding_envelopes}, there is therefore some $O_1 \in \Oc(A)$ such that $E_1$ is a reduced envelope of $\phi(O_1)$.

It remains to check that $O_1 \ll O_2$.  We see that $O_1\Res{\omega}(O_2)/\Res{\omega}(O_2)$ is compact, since $O_1 \le \Res{\omega}(O_2)U$.  By Lemma~\ref{lem:RIO_sigma-normal}(ii) we have $\Res{\omega}(\ol{\phi(O_1)}) = \ol{\phi(\Res{\omega}(O_1))}$, so in particular, $\phi(\Res{\omega}(O_1)) \le R_1$ and hence $\Res{\omega}(O_1) \le \phi\inv(R_1)$.  Since $O_1 \le O_2$ we also have $\Res{\omega}(O_1) \le \Res{\omega}(O_2)$; hence by \eqref{eq:res_component}, $\Res{\omega}(O_1) \le M$ for all $M \in \mc{M}$.  Thus $O_1 \ll O_2$.  This completes the proof of (ii) in Case 1.

By the Claim, Case 1 applies whenever $\phi$ is injective. 

\emph{Case 2: $\phi$ is a quotient map; $G = O_2$; and $|\mc{M}_1| = 1$.}

We can identify $H$ with the quotient $G/K$ of $G$.  Write $\mc{M}_1 = \{M^*\}$.  Then $M^*K/K$ is dense in $R_2$, so if we write $R_i = \widetilde{R}_i/K$, then $M^*K$ is dense in $\widetilde{R}_2$.  In particular, $\widetilde{R}_2 \le UM^*K$ for every open subgroup $U$ of $O_2$.  Write $E_1 = \widetilde{E}_1/K$.  The fact that $E_1R_2/R_2$ is compact means that $U\widetilde{R}_2$ contains a finite index subgroup of $\widetilde{E}_1$ for all compact open subgroups $U$ of $\widetilde{E}_1$; without loss of generality, we can replace $\widetilde{E}_1$ with its intersection with $U\widetilde{R}_2$, and so $\widetilde{E}_1 \le U\widetilde{R}_2$.  Since $\widetilde{R}_2 \le UM^*K$ and $U,K \le \widetilde{E}_1$, we see that
\[
\widetilde{E}_1 = UM^*K \cap \widetilde{E}_1 = U(M^*K \cap \widetilde{E}_1) = U(M^* \cap \widetilde{E}_1)K.
\]
Applying Lemma~\ref{lem:cocompact_gen} to the compactly generated quotient $E_1$ of $\widetilde{E}_1$, we see that there is a finite subset $S$ of $M^* \cap \widetilde{E}_1$ such that
\[
\widetilde{E}_1 = U\grp{S}K.
\]
We now set $O_1 = \grp{U,S}$.  Clearly $O_1 \in \Oc(G)$ and $O_1K/K = E_1$; it remains to check that $O_1 \ll O_2$.  Let $T = O_1 \cap M^*$.  Then we see that $T$ is normal in $O_1$ and $S \subseteq T$, so $O_1 = UT$ and hence $O_1/T$ is residually discrete; we see via Lemma~\ref{lem:finite_rank_residual}(ii) that $\Res{\omega}(O_1) = \Res{\omega}(T)$.  Clearly $T \le \Res{\omega}(O_2)$, so $O_1 \le U\Res{\omega}(O_2)$, and hence condition (a) of the partial order $\ll$ is satisfied.  For condition (b), consider $M \in \mc{M}$.  If $M = M^*$ then certainly $\Res{\omega}(T) \le M$, so $\Res{\omega}(O_1) \le M$.  If $M \neq M^*$ then $M \in \mc{M}_2$; then \eqref{eq:res_component} shows that $\widetilde{R_1} \cap \Res{\omega}(O_2) \le M$, and it is clear that $\Res{\omega}(O_1)$ is contained in both $\widetilde{R_1}$ and $\Res{\omega}(O_2)$, so again $\Res{\omega}(O_1) \le M$.  Thus condition (b) of the partial order is satisfied, so $O_1 \ll O_2$.  This completes the proof of (ii) in Case 2.

\emph{Case 3: $\phi$ is a quotient map.}

Again we identify $H$ with the quotient $G/K$ of $G$.  In this context we see that $O_2K/K$ is a finite index open subgroup of $E_2$; without loss of generality, we can replace $E_2$ with $O_2K/K$ and so assume that $O_2K/K = E_2$.  In turn, $O_2K/K$ is naturally isomorphic to $O_2/(O_2 \cap K)$, so $E_2$ can also be regarded as a quotient of $O_2$; since all the groups under consideration will be subgroups of $O_2$ or $E_2$, there is also no loss of generality to assume that $G = O_2$.

Consider now how $\mc{M}_1$ depends on $K$.  For each $M \in \mc{M}$, the factor $\Res{\omega}(G)/M$ is a chief factor of $G$, representing a chief block $\mf{a}_M$ of $G$; we see that $M \in \mc{M}_1$ if and only if $K \cap \Res{\omega}(G) \nleq M$, that is, if and only if $K$ covers $\mf{a}_M$.  Moreover, since $\Res{\omega}(G)/M$ is not elementary, we see that $\mf{a}_M$ is a robust chief block of $G$, and hence minimally covered by Corollary~\ref{cor:robust_block:countable}.  Writing $L_M := G_{\mf{a}_M}$, we then have $L_M \le K$ if and only if $M \in \mc{M}_1$.  Finally, given a quotient $G'$ of $G$, if we have a closed $G'$-invariant subgroup $N$ of $\Res{\omega}(G')$ such that $\Res{\omega}(G')/N \in \mc{S} \smallsetminus \Es{}$, then Lemma~\ref{lem:S_compression}(ii) shows that $\Res{\omega}(G')/N$ arises as the normal compression of $\Res{\omega}(G)/M$ induced by the quotient map for some $M \in \mc{M}$.  So we can effectively factorize $\phi$ so as to cover the elements of $\mc{M}_1$ one at a time: say $\mc{M}_1 = \{M_1,\dots,M_n\}$, then we take $K_i = \cgrp{L_{M_1},\dots,L_{M_i}}$.  We then have a sequence
\[
G \rightarrow G/K_1 \rightarrow G/K_2 \rightarrow \dots \rightarrow G/K_n \rightarrow G/K
\]
of quotient maps.  Part (ii) of the theorem is now proved for each quotient map in this sequence by Case 1 or Case 2; we deduce that (ii) also holds in Case 3.

In the general case, $\phi$ can be factorized as the composition of a quotient map and a regionally normal compression.  Thus (ii) follows by combining Cases 1 and 3.\end{proof}

One consequence of Theorem~\ref{thm:reg_quotient}(i) is that in a compactly generated \tdlcsc group $G$, open subgroups strongly below $G$ (in the sense of the order defined by \eqref{eq:RIO_ll}) are ``robustly of infinite index'' in a way that passes to quotients.

\begin{cor}\label{cor:robust_index}
Let $G$ be a compactly generated \tdlcsc group.  Let $\mc{H}$ be a set of compactly generated open subgroups such that $H_1,H_2 \in \mc{H} \Rightarrow \grp{H_1,H_2} \in \mc{H}$ and such that $H \ll G$ for all $H \in \mc{H}$.  Let $K = \grp{\mc{H}}$ and let $N$ be a closed normal subgroup of $G$ such that $G/N$ is not compact.  Then $KN$ has infinite index in $G$.
\end{cor}

\begin{proof}
Suppose that $KN$ has finite index in $G$.  By replacing $G$ with $KN$, we may assume $G = KN$.  Since $G$ is compactly generated and $\mc{H}$ consists of open subgroups, we see that in fact $G = HN$ for some $H \in \mc{H}$.  Since $H \ll G$, by Theorem~\ref{thm:reg_quotient}(i) the quotient $G/N$ must be compact, a contradiction.
\end{proof}

\subsection{Robustness of the partial order}

To summarize this section, here are the most important properties we have shown for the partial order $\ll$ on the regional poset of a first-countable \tdlc group $G$:
\begin{enumerate}[(A)]
\item Given $[A],[B] \in \Reg(G)$ such that $[A] \ll [B]$, then every open normal subgroup of $B$ contains a finite index subgroup of $A$.  (Clear from the definition)
\item Given $[A],[B],[C] \in \Reg(G)$ such that $[A] \le [B] \ll [C]$, then $[A] \ll [C]$.  (Lemma~\ref{lem:strong_order:basic})
\item Given $H \in \mc{RIO}(G)$, then the envelope map from $\Reg(H)$ to $\Reg(G)$ is a $\ll$-order embedding.  (Corollary~\ref{cor:res_RIO})
\item Given $\phi: G \rightarrow H$ be a dense regionally normal map to a first-countable \tdlc group $H$, and letting $\theta$ be the envelope map of $\phi$, then the following holds:
\begin{enumerate}[(i)]
\item Given $[O_1],[O_2] \in \Oc(G)$ such that $[O_1] \ll [O_2]$, then either $\theta([O_1]) = \theta([O_2])$ is the smallest element of $\Reg(H)$, or $\theta([O_1]) \ll \theta([O_2])$;
\item Given $[E_1],[E_2] \in \Oc(H)$ such that $[E_1] \ll [E_2]$, and $[O_2] \in \theta\inv([E_2])$, there exists $[O_1] \in \theta\inv([E_1])$ such that $[O_1] \ll [O_2]$.
\end{enumerate}
(Theorem~\ref{thm:reg_quotient})
\end{enumerate}
Say that a strict partial order $<^*$ on $\Reg(G)$ is \defbold{robust} if it has the properties (A)--(D).  Most of the results we will prove in later sections are also valid for any robust partial order on $\Reg(G)$, and for the rank function associated to it.  The main exceptions are as follows.  For the associated rank function $\xi^*$, the function $f$ such that $f(\xi^*(G)) \ge \xi(G)$ depends on the specifics of $<^*$ (Proposition~\ref{prop:elementary} below; the inequality $\xi^*(G) \le \xi(G)$ can be deduced from (A)).  The construction, given a well-founded compactly generated \tdlcsc group $G$, of a closed characteristic subgroup of $G$ of strictly smaller rank also depends on the specifics of $<^*$.

\begin{rmk}\label{rmk:omega}
Given a first-countable \tdlc group $G$, one could define a partial order $\ll'$ on $\Oc(G)$ as follows: we have $H \ll' K$ if $H \le K$, $K$ is noncompact, and the following two conditions are satisfied:
\begin{enumerate}[(a)']
\item $H\Res{}(K)/\Res{}(K)$ is compact;
\item For all closed normal subgroups $N$ of $\Res{}(K)$ such that $\Res{}(K)/N \in \Sclass$ and $|K:\N_K(N)|<\infty$, then $\Res{}(H) \le N$.
\end{enumerate}
Condition (a)' can be equivalently phrased as follows: every open normal subgroup of $K$ virtually contains $H$.  In other words, it is the weakest condition that implies property (A) above.  On the class $\Es{}$, this condition gives rise to a rank function that differs from $\xi$ by at most $1$.

The difficulty is in the application of the argument of Lemma~\ref{lem:semisimple:basic}, since in general, unlike the finite elementary rank residual, the discrete residual of a compactly generated \tdlc group need not be topologically perfect.  This difficulty could be avoided if we could be sure that every group in $\Sclass$ were abstractly simple, since this would mean that a normal compression with range in $\Sclass$ is necessarily an isomorphism, from which one would deduce that its domain also belongs to $\Sclass$.  However, it is still an open question whether all groups in $\Sclass$ are abstractly simple.  It seems likely that $\ll'$ is also robust, but at present we do not have the tools to prove it.

One possible way to avoid explicitly invoking (a class of) topologically simple factor groups would be to take the conclusion of Corollary~\ref{cor:robust_index} and incorporate it into the definition: instead of (b) or (b)', one could impose the condition, given $H,K \in \Oc(G)$, that for every closed subgroup $N$ of $K$ such that $|K:\N_K(N)|<\infty$ and $\N_K(N)/N$ is noncompact, then $\N_H(N)N$ has infinite index in $K$.  However, it is not clear if this would result in a robust partial order.\end{rmk}

\section{The class of $\mc{S}$-well-founded groups}\label{sec:well-founded}

\subsection{Well-founded posets}\label{sec:poset}

A poset $P$ is \defbold{well-founded} if every nonempty subset has a minimal element; equivalently, given Zorn's lemma, $P$ is \defbold{Artinian}, meaning it has no infinite descending chain.  Given a well-founded poset $P$, we can define an ordinal-valued rank function $\rho_P$ on $P$ by recursion.

\begin{lem}\label{lem:poset_rank}
Let $P$ be a well-founded poset, with strict partial order $<$.  Then there is a unique ordinal-valued function $\rho_P$ with the following properties:
\begin{enumerate}[(a)]
\item If $p \in P$ is minimal, then $\rho_P(p)=0$;
\item If $p \in P$ is not minimal, then
\[
\rho_P(p) = \sup \{\rho_P(q) \mid q < p\} + 1.
\]
\end{enumerate}
\end{lem}

\begin{proof}
We construct a sequence of subsets $P_{\alpha}$ of $P$ as follows: 
\begin{enumerate}[(i)]
\item $P_0$ is the set of minimal elements;
\item $P_{\alpha+1}$ is the set $\{p \in P \mid q<p \Rightarrow q \in P_{\alpha}\}$;
\item If $\lambda$ is an infinite limit ordinal, $P_{\lambda} = \bigcup_{\alpha < \lambda}P_{\alpha}$.
\end{enumerate}
It is easy to see that $P_{\alpha}$ is downward-closed for all $\alpha$, and that $P_{\alpha} \subseteq P_{\alpha'}$ whenever $\alpha \le \alpha'$.

The sequence eventually repeats, that is, we reach an ordinal $\beta$ such that $P_{\beta} = P_{\beta+\gamma}$ for some $\gamma > 0$; it then follows that $P_{\beta} = P_{\beta+1}$.  We claim that in this case $P_{\beta} = P$; suppose for a contradiction that $P \smallsetminus P_{\beta}$ is nonempty.  Then since $P$ is well-founded, the set $P \smallsetminus P_{\beta}$ has a minimal element $p$.  But then $q < p \Rightarrow q \in P_{\beta}$, so $p \in P_{\beta+1}$, a contradiction.

Now set
\[
\rho_P(p) = \min \{\alpha: p \in P_{\alpha}\}.
\]
Since the sets $P_{\alpha}$ exhaust $P$, this is a well-defined ordinal-valued function; moreover, from the construction of the sets $P_{\alpha}$, we see that $\rho_P$ cannot take any infinite limit ordinal as a value.  Clearly $\rho_P(p)=0$ if and only if $p$ is minimal; suppose that $p$ is not minimal. Then $\rho_P(p) = \alpha+1$ for some ordinal $\alpha$.  In particular, we have $q<p \Rightarrow q \in P_{\alpha}$, so $\sup \{\rho_P(q) \mid q < p\} \le \alpha$.  At the same time, given $\beta < \alpha$, then $p \not\in P_{\beta+1}$, so there exists $q_{\beta} < p$ such that $q_{\beta} \not\in P_{\beta}$, and hence $\sup \{\rho_P(q) \mid q < p\} > \beta$.  Thus $\sup \{\rho_P(q) \mid q < p\} = \alpha$, and hence
\[
\rho_P(p) = \alpha+1 = \sup \{\rho_P(q) \mid q < p\} + 1.
\]

The uniqueness of $\rho_P$ follows by reversing the construction: if $\rho'_P(p)$ is any ordinal-valued function satisfying the conditions (a) and (b), then the sequence $(P'_{\alpha})$ defined by
\[
P'_{\alpha} = \{p \in P \mid \rho'_P(p) \le \alpha\}
\]
satisfies the conditions (i)--(iii) at the start of the proof, from which we deduce that $P'_{\alpha} = P_{\alpha}$ for all $\alpha$ and hence $\rho'_P = \rho_P$.
\end{proof}

We then set
\[
\rho(P) := \sup\{ \rho_P(p) \mid p \in P\} + 1.
\]
The fact that $P$ is well-founded ensures that $\rho(P)$ is well-defined.  A consequence of the definition is that the rank of a well-founded poset is always a successor ordinal.

More generally, if the poset $P$ is not well-founded, we can still define $\rho_P(p)$, as follows.  Define the \defbold{well-founded part} $Q$ of $P$ to consist of all elements $p \in P$ that do not belong to any infinite descending chain.  Then $Q$ is well-founded, and for $p \in Q$ we write $\rho_P(p) = \rho_Q(p)$.

\begin{rmk}A common convention in the literature is to define $\rho_P(p) =  \sup \{\rho_P(q)+1 \mid q < p\}$, so that for instance if $P$ is a von Neumann ordinal, one would have $\rho_P(\alpha) = \alpha$ for all $\alpha \in P$.  We depart from this convention in order to make some induction arguments work more smoothly.\end{rmk}

\subsection{Defining the class of $\mc{S}$-well-founded groups}

\begin{defn}
Let $\mc{S}$ be a class of \tdlcsc groups with property $\propS$.  A \tdlc group $G$ is \defbold{$\mc{S}$-well-founded} if $G$ is first-countable and $(\Reg(G),\ll_\mc{S})$ is well-founded.

Given a first-countable \tdlc group $G$, we define $\xip{\mc{S}}(G) = \rho(\Reg(G),\ll_\mc{S})$ and write $\rho_G$ for the rank function $\rho_{(\Reg(G),\ll_\mc{S})}$ (where these are defined).

Write $\Es{\mc{S}}$ for the class of $\mc{S}$-well-founded \tdlc groups, and given an ordinal $\alpha$, write $\Es{\mc{S}}(\alpha)$ for the class of $\mc{S}$-well-founded groups $G$ with $\xip{\mc{S}}(G) \le \alpha$.
\end{defn}

Let us note how the class $\Es{\mc{S}}$ and the associated rank changes if we replace $\mc{S}$ with a smaller or larger class of groups.

\begin{lem}\label{lem:rank_S}
Let $\mc{S}$ and $\mc{S}'$ be classes of \tdlcsc groups with property $\propS$, such that $\mc{S} \subseteq \mc{S}'$, and let $\alpha$ be an ordinal.  Then $\Es{\mc{S}}(\alpha) \subseteq \Es{\mc{S}'}(\alpha)$.
\end{lem}

\begin{proof}
Let $G$ be a first-countable \tdlc group.  Given $H,K \in \Oc(G)$, we see from the definitions that
\[
H \ll_{\mc{S}'} K \Rightarrow H \ll_{\mc{S}} K,
\]
in other words, $\ll_{\mc{S}}$ is a finer partial order than $\ll_{\mc{S}'}$.  The analogous statement then follows for $\Reg(G)$.  Writing $\rho$ and $\rho'$ for the (partial) rank functions on $\Reg(G)$ induced by $\mc{S}$ and $\mc{S}'$ respectively, we deduce (by induction on $\alpha$) that
\[
\forall [H] \in \Reg(G)\; \forall \alpha: \; \rho[H] \le \alpha \Rightarrow \rho'[H] \le \alpha.
\]
The conclusions are now clear. 
\end{proof}

For the rest of this section, we fix once more a class $\mc{S}$ of \tdlcsc groups with property $\propS$, and write $\ll$ for the associated partial order $\ll_{\mc{S}}$.

Here is a translation of what the definition of $\mc{S}$-well-founded groups means in terms of compactly generated open subgroups (without taking the equivalence relation of commensurability).

\begin{lem}\label{lem:pseudo_chain}
Let $G$ be a first-countable \tdlc group.  Then exactly one of the following holds:
\begin{enumerate}[(i)]
\item $G$ is $\mc{S}$-well-founded;
\item There exists an infinite sequence $(G_i)_{i \in \Nb}$, $G_i \in \Oc(G)$, such that $G_{i+1} \ll G_i$ for all $i \in \Nb$.
\end{enumerate}
Moreover, the equation 
\[
\rho_{(\Oc(G),\ll)}(H) = \rho_{(\Reg{}(G),\ll)}[H]
\]
holds for all $H \in \Oc(G)$ such that either side of the equation is defined.
\end{lem}

\begin{proof}
If an infinite sequence exists as in (ii), then $([G_i])_{i \in \Nb}$ is an infinite descending chain in $(\Reg(G),\ll)$, so $G$ is not $\mc{S}$-well-founded.

Conversely, suppose that $G$ is not $\mc{S}$-well-founded, that is, by Zorn's lemma we have an infinite descending chain $([G_i])_{i \in \Nb}$ in $(\Reg(G),\ll)$, where $G_i \in \Oc(G)$.  Then $G_i$ virtually contains $G_{i+1}$ for all $i$; so by setting $H_i = \bigcap_{j \le i}G_j$, we obtain a sequence $(H_i)_{i \in \Nb}$ in $\Oc(G)$ such that $[H_i] = [G_i]$.  By Lemma~\ref{lem:strong_order:basic}, we have $H_{i+1} \ll H_i$ for all $i \in \Nb$, so (ii) holds.

Let $\pi$ be the map from $\Oc(G)$ to $\Reg(G)$ given by $\pi(H) = [H]$.  For the equation
\[
\rho_{(\Oc(G),\ll)}(H) = \rho_{(\Reg{}(G),\ll)}[H],
\]
we construct the rank sets 
\[
P_{\alpha} := \{H \in \Oc(G) \mid \rho_{(\Oc(G),\ll)}(H) \le \alpha\} \text{ and } \widetilde{P}_{\alpha} := \{[H] \in \Reg(G) \mid \rho_{(\Reg(G),\ll)}([H]) \le \alpha\},
\]
showing by induction on $\alpha$ that $\pi\inv(\widetilde{P}_{\alpha}) = P_{\alpha}$ for each ordinal $\alpha$.  We note that, given $H \in \Oc(G)$, then every $[K] \in \Reg(G)$ such that $[K] \ll [H]$ is witnessed by some $K \in \Oc(G)$ such that $K \ll H$, and conversely if $K \ll H$ in $(\Oc(G),\ll)$, then $[K] \ll [H]$ in $(\Reg(G),\ll)$.  It is then immediately clear that $P_{0} = \pi\inv(\widetilde{P}_{0})$.  Similarly, we have $[H] \in \widetilde{P}_{\alpha+1}$ if and only if, for all $K \ll H$, we have $[K] \in \widetilde{P}_{\alpha}$, or equivalently (by the inductive hypothesis) $K \in P_{\alpha}$; thus $[H] \in \widetilde{P}_{\alpha+1}$ if and only if $H \in P_{\alpha+1}$.  The fact that $P_{\alpha} = \pi\inv(\widetilde{P}_{\alpha})$ when $\alpha$ is a limit ordinal follows from the inductive hypothesis and the observation that the preimage of a union is the union of the preimages.
\end{proof}

\begin{cor}
A first-countable \tdlc group $G$ is $\mc{S}$-well-founded if and only if $(\Oc(G),\ll)$ is well-founded.
\end{cor}

Given Lemma~\ref{lem:pseudo_chain}, we can also write $\rho_G$ for the rank function $\rho_{(\Oc(G),\ll)}$ without any danger of confusion (once the class $\mc{S}$ is specified).  Some caution is required, however, in the relationship between $\rho_G(H)$ and $\xip{\mc{S}}(H)$ for $H \in \Oc(G)$.  For instance, even if $G$ is compactly generated, we are not claiming that $\rho_G(G)$ is the largest value of $\rho_G(H)$ for $H$ in $\Oc(G)$.

\subsection{Rank inequalities for RIO subgroups}\label{sec:rank_inequalities}

Our convention, given a first-countable \tdlc group $G$, is that any ordinal upper bound on $\xip{\mc{S}}(G)$ is also an assertion that $G \in \Es{\mc{S}}$, but specific upper bounds will also be useful, for instance in arguments by induction.  Here are some general rank inequalities that we can prove with respect to $G$ and its RIO subgroups.

\begin{lem}\label{lem:rank_inequalities}
Let $G$ be a first-countable \tdlc group.
\begin{enumerate}[(i)]
\item Suppose that $G \in \Es{\mc{S}}$ and $H \in \mc{RIO}(G)$.  Then $H \in \Es{\mc{S}}$ and $\xip{\mc{S}}(H) \le \xip{\mc{S}}(G)$.
\item Suppose that $H$ is a closed cocompact RIO subgroup of $G$ such that $H \in \Es{\mc{S}}$.  Then $G \in \Es{\mc{S}}$ and $\xip{\mc{S}}(G) = \xip{\mc{S}}(H)$.
\item Suppose that $\mc{D}$ is a family of RIO subgroups of $G$, directed upwards by inclusion, such that $D \in \Es{\mc{S}}$ for all $D \in \mc{D}$ and $\bigcup_{D \in \mc{D}}D$ is dense in $G$.  Then $G \in \Es{\mc{S}}$ and $\xip{\mc{S}}(G) = {\sup_{D \in \mc{D}}}^+\xip{\mc{S}}(D)$.
\item Suppose $G \in \Es{\mc{S}}$, and let $H$ and $K$ RIO subgroups of $G$ such that $H \ll K$, where $K$ is compactly generated and $H$ is compactly generated or open in $K$.  Then $\xip{\mc{S}}(H) < \xip{\mc{S}}(K)$.
\end{enumerate}
\end{lem}

\begin{proof}
Before we start, note that if $H \in \Oc(G)$ then the functions $\rho_G$ and $\rho_H$ agree, that is, $\rho_G(K) = \rho_H(K)$ for all $K \in \Oc(H)$.

Parts (i) and (ii) are clear from Corollary~\ref{cor:res_RIO}.

For (iii), we first consider $H \in \Oc(G)$.  Let $U$ be a compact open subgroup of $G$ and let $X = \bigcup_{D \in \mc{D}}D$.  Since $X$ is dense in $G$, it has dense intersection with $H$, and hence by Lemma~\ref{lem:cocompact_gen}, we have $H = \grp{S}U$ for a finite subset $S$ of $X$.  Since $\mc{D}$ is directed, in fact $S \subseteq D \in \mc{D}$, and then there is $K \in \Oc(D \cap H)$ such that $S \subseteq K$, so that $H = KU$.  Now $D \cap H$ is a RIO subgroup of $G$, so $K$ is a compactly generated IO subgroup of $G$ and hence of $D$ and $H$.  Parts (i) and (ii) now imply that $H \in \Es{\mc{S}}$ with $\xip{\mc{S}}(H) = \xip{\mc{S}}(K) \le \xip{\mc{S}}(D)$.  In particular, it is now clear that $(\Reg(G),\ll)$ has no infinite descending chain, so $G \in \Es{\mc{S}}$.  

Let $\alpha+1 =  {\sup_{D \in \mc{D}}}^+\xip{\mc{S}}(D)$.  Since for all $H \in \Oc(G)$ there exists $D \in \mc{D}$ such that $\xip{\mc{S}}(H) \le \xip{\mc{S}}(D)$, we see that
\[
\forall H \in \Oc(G): \sup\{\rho_G(K) \mid K \in \Oc(H)\} \le \alpha,
\]
and hence
\[
\xip{\mc{S}}(G) = \sup\{\rho_G(K) \mid K \in \Oc(G)\} + 1 \le \alpha+1.
\]
On the other hand, by (i) we have $\xip{\mc{S}}(D) \le \xip{\mc{S}}(G)$ for all $D \in \mc{D}$.  From the way the rank is calculated, we know that $\xip{\mc{S}}$ only takes successor ordinal values, so $\alpha+1 \le \xip{\mc{S}}(G)$.

For (iv), in the case that $H$ is compactly generated, first take a reduced envelope $H'$ of $H$ in $K$.  Then $\xip{\mc{S}}(H) \le \xip{\mc{S}}(H')$ by (i), and by Proposition~\ref{prop:ll_RIO} we have $H' \ll K$.  Thus without loss of generality, $H$ is open in $K$.  We then compare the elements of $\Oc(H)$ with $K$.  Given $O \in \Oc(H)$, then from the fact that $H \ll K$, there is $O' \in \Oc(H)$ such that $O \le O' \ll K$; hence by Lemma~\ref{lem:strong_order:basic}, we have $O \ll K$.  Consequently,
\[
\rho_K(K) \ge \sup\{\rho_K(O) \mid O \in \Oc(H)\} + 1 = \sup\{\rho_H(O) \mid O \in \Oc(H)\} + 1 = \xip{\mc{S}}(H),
\]
and hence $\xip{\mc{S}}(K) \ge \xip{\mc{S}}(H)+1 > \xip{\mc{S}}(H)$.
\end{proof}

Given a compactly generated $G \in \Es{\mc{S}}$, we can now obtain a characteristic closed normal subgroup of $G$ of smaller rank.

\begin{thm}\label{thm:rank_reduction}
Let $G$ be a noncompact compactly generated \tdlcsc group such that $G \in \Es{\mc{S}}$.  Then $G$ admits a finite series
\[
G_0 \le R_n \le \dots R_1 \le R_0 = G
\]
of closed characteristic subgroups with the following properties:
\begin{enumerate}[(i)]
\item $\xip{\mc{S}}(G_0) < \xip{\mc{S}}(G)$;
\item $R_n = \Res{\omega}(G)$ and $G_0$ is expressible as the intersection of $R_n$ with a finite number (possibly zero) of closed normal subgroups $N$ of $R_n$ such that $R_n/N \in \mc{S} \smallsetminus \Es{}$ and $|G:\N_G(N)|<\infty$;
\item For $1 \le i \le n$, and given $H \in \Oc(R_{i-1}/R_i)$, then $H$ is a SIN group.
\end{enumerate}
\end{thm}

\begin{proof}
We start by forming the elementary rank-$2$ series $(R_i)_{i \in \Nb}$ of the compactly generated group $G \in \Es{\mc{S}}$, starting with $R_0 = G$.  Each term in this series is clearly characteristic in $G$, and by Lemma~\ref{lem:finite_rank_residual}, the series terminates after finitely many steps; let $n$ be the least $n \in \Nb$ such that $R_n = R_{n+1}$.  Then $R_n = \Res{\omega}(G)$.  By construction, $R_{i-1}/R_i \in \Es{}(2)$ for $1 \le i \le n$, so by Lemma~\ref{lem:elementary2}, every compactly generated open subgroup of $R_{i-1}/R_i$ is a SIN group.

Let $\mc{N}$ be the set of closed normal subgroups $N$ of $\Res{\omega}(G)$ such that $\Res{\omega}(G)/N \in \mc{S} \smallsetminus \Es{}$ and $|G:\N_G(N)|<\infty$, and set $G_0 = \Res{\omega}(G) \cap \bigcap_{N \in \mc{N}}N$.  By Lemma~\ref{lem:S_quotient_count}, $\mc{N}$ is finite, and it is clear that $G_0$ is a closed characteristic subgroup of $G$.  All that remains is to show $\xip{\mc{S}}(G_0) < \xip{\mc{S}}(G)$. Let $U$ be a compact open subgroup of $G$ and let $H = UG_0$; then $\xip{\mc{S}}(G_0) = \xip{\mc{S}}(H)$ by Lemma~\ref{lem:rank_inequalities}(ii).  Given $O \in \Oc(H)$ then
\[
\Res{\omega}(O) \le \Res{\omega}(UG_0) = \Res{\omega}(G_0) \le G_0,
\]
from which it is clear that $O \ll G$.  Hence $\xip{\mc{S}}(H) < \xip{\mc{S}}(G)$ by Lemma~\ref{lem:rank_inequalities}(iv).
\end{proof}

The following will be useful inductive steps in what follows.

\begin{lem}\label{lem:step}
Let $G$ be a first-countable \tdlc group.
\begin{enumerate}[(i)]
\item Let $\alpha$ be an ordinal.  Suppose that for all pairs $O_1,O_2 \in \Oc(G)$ such that $O_1 \ll O_2$, we have $O_1 \in \Es{\mc{S}}(\alpha+1)$.  Then $G \in \Es{\mc{S}}(\alpha+2)$.
\item Let $\alpha$ be a limit ordinal.  Suppose that for all pairs $O_1,O_2 \in \Oc(G)$ such that $O_1 \ll O_2$, we have $O_1 \in \Es{\mc{S}}(\alpha)$, but $\xip{\mc{S}}(G) > \alpha+1$.  Then there is $H \in \Oc(G)$ such that $\sup\{\xip{\mc{S}}(K) \mid K \in \Oc(G), K \ll H\} = \alpha$.
\end{enumerate}
\end{lem}

\begin{proof}
(i)
Given $H \in \Oc(G)$, we have
\[
\rho_G(H) = \sup\{\rho_G(K) \mid K \in \Oc(G), K \ll H\} + 1 = \sup\{\rho_K(K) \mid K \in \Oc(G), K \ll H\} + 1.
\]
Given $K \ll H$, then $\xip{\mc{S}}(K) \le \alpha+1$ by hypothesis, so $\rho_K(K) \le \alpha$.  Hence $\rho_G(H) \le \alpha+1$ for all $H \in \Oc(G)$, so $\xip{\mc{S}}(G) \le \alpha+2$.

(ii)
By part (i), $\xip{\mc{S}}(G) \le \alpha+2$, so in fact $\xip{\mc{S}}(G) = \alpha+2$.  We therefore have 
\[
\alpha+1 = \sup\{\rho_G(H) \mid H \in \Oc(G)\},
\]
so there is $H \in \Oc(G)$ such that $\rho_G(H) = \alpha+1$. In other words, we have
\[
\sup\{\rho_G(K) \mid K \in \Oc(G), K \ll H\} = \alpha,
\]
so
\[
\sup\{\xip{\mc{S}}(K) \mid K \in \Oc(G), K \ll H\} = \alpha. \qedhere
\]
\end{proof}

Note that for a first-countable \tdlc group $G$, given $[H] \in \Reg(G)$ then $[H]$ is $\ll$-minimal if and only if $H$ is compact. We can thus describe the $\mc{S}$-well-founded groups of rank $1$ (independently of the choice of $\mc{S}$) as follows.

\begin{lem}\label{lem:rank_one}
We have $\Es{\mc{S}}(1) = \Es{}(\mathrm{r})$.
\end{lem}

In general, given $G \in \Sclass$, it is not clear if $G$ is $\mc{S}$-well-founded, or if it is, what value $\xip{\mc{S}}(G)$ could take (other than the fact that it must be at least $2$); but we can at least determine $\rho_G(G)$.  (Again, the choice of $\mc{S}$ is not important here; note however that $\rho_G(G)$ only makes sense if $G$ is compactly generated, hence the assumption that $G \in \Sclass$.)

\begin{lem}\label{lem:simple_rho}
Let $G \in \Sclass$ and let $\mc{H}$ be the set of $H \in \Oc(G)$ such that $H$ is elementary and $\xi(H) < \omega$.  Then exactly one of the following holds:
\begin{enumerate}[(i)]
\item Every $H \in \mc{H}$ is compact and $\rho_G(G) = 1$;
\item Some $H \in \mc{H}$ is noncompact and $\rho_G(G)=2$.
\end{enumerate}
\end{lem}

\begin{proof}
Note that $\mc{H}$ is exactly the set of compactly generated open subgroups $H$ of $G$ such that $\Res{\omega}(H) = \triv$.  From the definition of $\ll$, we then see that $\mc{H}$ is in fact the set of $H \in \Oc(G)$ such that $H \ll G$.  Thus
\[
\rho_G(G) = \sup\{\rho_G(H) \mid H \in \mc{H}\} + 1.
\]
Given $H \in \mc{H}$, if $H$ is compact then clearly $\rho_G(H) = 0$.  On the other hand, if $H$ is not compact and $K \in \Oc(G)$, then $[K] \ll [H]$ if and only if $K$ is compact, so $\rho_G(H) = 1$.  The conclusion is now clear.
\end{proof}

As we will see in Section~\ref{sec:trees}, both cases of Lemma~\ref{lem:simple_rho} occur, giving examples of nonelementary groups of small $\Sclass$-well-founded rank; indeed for $i \in \{1,2\}$ there are examples of groups $G \in \Sclass \cap \Es{\Sclass}$ with $\rho_G(G) = i$ and $\xip{\Sclass}(G) = i+1$.

\subsection{Relationship to elementary groups}

The class of $\mc{S}$-well-founded groups contains the class of regionally elementary \tdlc groups, and the respective rank functions are each bounded by a function of the other.

\begin{prop}\label{prop:elementary}
Let $G$ be a regionally elementary \tdlc group.  Then $G$ is $\mc{S}$-well-founded and
\[
\xip{\mc{S}}(G) \le \xi(G) \le \omega.\xip{\mc{S}}(G)+1.
\]
If $\xi(G) \le \omega+1$, then $\xip{\mc{S}}(G) \le 2$.
\end{prop}

\begin{proof}
We proceed by induction on $\xi(G)$.  Let us begin with the case that $\xi(G) \le \omega+1$; in other words, every $H \in \Oc(G)$ has $\Res{\omega}(H) = \triv$.  If $G$ is regionally elliptic, then $G \in \Es{\mc{S}}(1)$ by Lemma~\ref{lem:rank_one}.  Otherwise, there exists $H \in \Oc(G)$ such that $H$ is not compact.  In the latter case we see that $\xi(G) \ge 2$, while the maximal length of a chain in $(\Reg(G),\ll)$ is achieved by $\{[H],[U]\}$ where $U$ is a compact open subgroup of $H$, so $G \in \Es{\mc{S}}(2)$.  In either case, $\xip{\mc{S}}(G) \le \xi(G)$ and $\omega.\xip{\mc{S}}(G) +1 \ge \xi(G)$.

For the rest of the proof we suppose that $\xi(G) = \alpha + n$, where $\alpha$ is an infinite limit ordinal and $n < \omega$.  We split into two cases, according to whether $n=1$ or $n > 1$.  (Note that $n \neq 0$, since $\xi$ takes only successor ordinals as values.).

Suppose $n > 1$ and let $H \in \Oc(G)$.  If $\xi(H) = 2$, then $\rho_G(H) \le 1$ by the first paragraph.  If $\xi(H) > 2$, we consider $K \in \Oc(H)$ such that $K \ll H$.  We see that
\[
\xi(\Res{\omega}(H)) \le \xi(\Res{}(H)) < \xi(H),
\]
and it follows from Lemma~\ref{lem:elementary:cocompact} that $\xi(K) < \xi(H)$.  In particular, $\xi(K) \le \alpha+(n-1)$; by the inductive hypothesis, $\xip{\mc{S}}(K) \le \alpha+(n-1)$, and hence $\rho_G(K) \le \alpha+(n-2)$.  We deduce that $\rho_G(H) \le \alpha + (n-1)$ for all $H \in \Oc(G)$, so $G \in \Es{\mc{S}}(\alpha+n)$.  For the upper bound, since $\xi(G)$ is not the successor of a limit ordinal, the elementary rank of $G$ is achieved by some $H \in \Oc(G)$.  By Lemma~\ref{lem:finite_rank_residual}, we have $H/\Res{\omega}(H) \in \Es{}(\omega)$; by Lemma~\ref{lem:elementary_extension}, it follows that $\xi(\Res{\omega}(H)) \ge \alpha+1$.  Let $U$ be a compact open subgroup of $H$ and let $K = \Res{\omega}(H)U$; then $\xi(K) \ge \alpha+1$ and $K \ll H$.  By the inductive hypothesis, $\omega.\xip{\mc{S}}(K) \ge \alpha$; by Lemma~\ref{lem:rank_inequalities} we have $\xip{\mc{S}}(G) > \xip{\mc{S}}(K)$.  Thus
\[
\omega.\xip{\mc{S}}(G) \ge \omega(\xip{\mc{S}}(K)+1) = \omega.\xip{\mc{S}}(K) + \omega \ge \alpha + \omega \ge \xi(G).
\]

Now suppose $n=1$.  Then $G$ is not compactly generated and every $H \in \Oc(G)$ has elementary rank strictly less than $\alpha$.  Hence by the inductive hypothesis we have $\xip{\mc{S}}(H) < \alpha$ for all $H \in \Oc(G)$, so $G \in \Es{\mc{S}}(\alpha+1)$.  For the upper bound on $\xi(G)$, we see that for all $\beta < \alpha$ there is $H \in \Oc(G)$ with $\xi(H) \ge \beta+1$, so by the inductive hypothesis, $\omega.\xip{\mc{S}}(H) \ge \beta$.  Since $\xip{\mc{S}}(G) \ge \xip{\mc{S}}(H)$ for all $H \in \Oc(G)$, it follows that $\omega.\xip{\mc{S}}(G) \ge \beta$ for all $\beta < \alpha$, so $\omega.\xip{\mc{S}}(G) + 1 \ge \alpha + 1$.
\end{proof}

In the other direction, we have a version of the elementary dichotomy for $G \in \Es{\mc{S}}$ that is an easy consequence of the definition of $\mc{S}$-well-founded groups.

\begin{prop}\label{prop:dichotomy}
Let $G \in \Es{\mc{S}}$.  Then one of the following holds:
\begin{enumerate}[(i)]
\item $G$ is regionally elementary;
\item There is $O \in \Oc(G)$, $S \in \mc{S} \smallsetminus \Es{}$ and a quotient map $\pi: \Res{\omega}(O) \rightarrow S$.
\end{enumerate}
\end{prop}

\begin{proof}
The two cases are clearly mutually exclusive, so let us assume (ii) is false and proceed by induction on $\xip{\mc{S}}(G)$.  Fix $H \in \Oc(G)$; we aim to show that $H$ is elementary.  If $H$ is compact, there is nothing to prove, so suppose that $H$ is not compact.  Then given $O \in \Oc(H)$, we have $O \ll H$ if and only if $O\Res{\omega}(H)/\Res{\omega}(H)$ is compact.  In particular, we see that if $U$ is a compact open subgroup of $G$ and $O \in \Oc(U\Res{\omega}(H))$, then $O \ll H$, so $\xip{\mc{S}}(O) < \xip{\mc{S}}(H) \le \xip{\mc{S}}(G)$ by Lemma~\ref{lem:rank_inequalities}, and hence $O$ is elementary by the inductive hypothesis.  From there it follows that the \tdlcsc group $U\Res{\omega}(H)$ is elementary, and hence $H$ is elementary.
\end{proof}

As we will see later (Theorem~\ref{thm:nilpotent_tree}), for $\mc{S} = \Sclass$ there are some groups in $\Sclass \smallsetminus \Es{\Sclass}$ that do not satisfy the dichotomy of Proposition~\ref{prop:dichotomy}, but the known examples still leave open Question~\ref{que:dichotomy}.

\subsection{Quotients, extensions and dense regionally normal maps}\label{sec:extensions}

Like the class of elementary groups, the class of $\mc{S}$-well-founded groups is closed under extensions and dense regionally normal maps, and we obtain bounds on the rank that are analogous to the elementary case.

\begin{thm}\label{thm:extensions}
Let $G$ be a first-countable \tdlc group.
\begin{enumerate}[(i)]
\item Let $H$ be a first-countable \tdlc group and let $\phi: G \rightarrow H$ be a dense regionally normal map.  If $G \in \Es{\mc{S}}$, then also $H \in \Es{\mc{S}}$ and $\xip{\mc{S}}(H) \le \xip{\mc{S}}(G)$.
\item Suppose that $N$ is a closed normal subgroup of $G$.  Then $G \in \Es{\mc{S}}$ if and only if $N,G/N \in \Es{\mc{S}}$.  If $G \in \Es{\mc{S}}$, then
\[
\max\{\xip{\mc{S}}(N),\xip{\mc{S}}(G/N)\} \le \xip{\mc{S}}(G) \le (\xip{\mc{S}}(N)-1) + \xip{\mc{S}}(G/N).
\]
\end{enumerate}
\end{thm}

\begin{proof}
(i)
Assume $G \in \Es{\mc{S}}$; we proceed by induction on $\xip{\mc{S}}(G)$.  In the base case, $\xip{\mc{S}}(G)=1$, so $G$ is regionally elliptic.  Given $E \in \Oc(H)$, by Lemma~\ref{lem:cocompact_gen} we can then write $E = \phi(K)U$, where $K$ is a finitely generated subgroup of $G$.  In fact, $K$ has compact closure in $G$, so $\phi(K)$ also has compact closure and hence $E$ is compact; thus $H$ is regionally elliptic, so $H \in \Es{\mc{S}}(1)$.

Now let us suppose $\xip{\mc{S}}(G) = \alpha+1$ for some $\alpha > 0$, and that the conclusions hold whenever $\xip{\mc{S}}(G) \le \alpha$.

\emph{Claim: Given $E_1,E_2 \in \Oc(H)$ such that $E_1 \ll E_2$, and given $O_2 \in \Oc(G)$ such that $E_2$ is a reduced envelope of $\phi(O_2)$, there is $O_1 \in \Oc(G)$ such that $O_1 \ll O_2$, $E_1$ is a reduced envelope of $\phi(O_1)$, and $\xip{\mc{S}}(E_1) \le \xip{\mc{S}}(O_1) \le \alpha$.}

Let $O_2 \in \Oc(G)$ be such that $E_2$ is a reduced envelope of $\phi(O_2)$.  Then by Theorem~\ref{thm:reg_quotient}(ii), there is $O_1 \in \Oc(G)$ such that $O_1 \ll O_2$ and $E_1$ is a reduced envelope of $\phi(O_1)$.  Then $\ol{\phi(O_1)}$ is a locally normal, in particular RIO, subgroup of $E_1$, and also $\ol{\phi(O_1)}$ is cocompact in $E_1$; hence $\xip{\mc{S}}(\ol{\phi(O_1)}) = \xip{\mc{S}}(E_1)$ by Lemma~\ref{lem:rank_inequalities}(ii).  At the same time we have $\xip{\mc{S}}(O_1) < \xip{\mc{S}}(G)$ by Lemma~\ref{lem:rank_inequalities}(iv), so $\xip{\mc{S}}(O_1) \le \alpha$.  We also see that $\phi$ restricts to a dense regionally normal map from $O_1$ to $\ol{\phi(O_1)}$.  By the inductive hypothesis, $\xip{\mc{S}}(\ol{\phi(O_1)}) \le \xip{\mc{S}}(O_1)$, so $\xip{\mc{S}}(E_1) \le \xip{\mc{S}}(O_1)$, proving the claim.

There are now two cases, according to whether $\alpha$ is a successor ordinal or a limit ordinal.

If $\alpha$ is a successor ordinal, by Lemma~\ref{lem:step}(i) we have $\xip{\mc{S}}(H) \le (\alpha-1)+2 = \xip{\mc{S}}(G)$.

If $\alpha$ is a limit ordinal, suppose for a contradiction that $\xip{\mc{S}}(H) > \alpha+1$.  Then by Lemma~\ref{lem:step}(ii) there is $E \in \Oc(H)$ such that 
\[
\sup\{\xip{\mc{S}}(E_1) \mid E_1 \in \Oc(H), E_1 \ll E\} = \alpha.
\]
Now take $O \in \Oc(G)$ such that $E$ is a reduced envelope of $\phi(O)$.  By the Claim, we see that
\[
\sup\{\xip{\mc{S}}(O_1) \mid O_1 \in \Oc(G), O_1 \ll O\} = \alpha.
\]
In particular, we see that $\rho_G(O) \ge \alpha+1$, so $\xip{\mc{S}}(G) \ge \alpha+2$, a contradiction.  From this contradiction we conclude that $\xip{\mc{S}}(H) \le \alpha+1 = \xip{\mc{S}}(G)$.  This proves (i).

(ii)
Let $Q = G/N$ and let $\pi:G\rightarrow Q$ be the quotient map.  If $G \in \Es{\mc{S}}$, then it follows from Lemma~\ref{lem:rank_inequalities}(i) that $N \in \Es{\mc{S}}(\xip{\mc{S}}(G))$ and from part (i) that $Q \in \Es{\mc{S}}(\xip{\mc{S}}(G))$.

It remains to prove the upper bound for $\xip{\mc{S}}(G)$ in terms of $\xip{\mc{S}}(N)$ and $\xip{\mc{S}}(Q)$, under the hypothesis that $N,Q \in \Es{\mc{S}}$.  We argue by induction on $\xip{\mc{S}}(Q)$. In the base case, $\xip{\mc{S}}(Q)= 1$, so $Q$ is regionally elliptic by Lemma~\ref{lem:rank_one}.  We then see that $G$ is a directed union of open subgroups $O$ such that $N$ is a cocompact normal subgroup of $O$.  Applying parts (ii) and (iii) of Lemma~\ref{lem:rank_inequalities}, we see that $G \in \Es{\mc{S}}$ and $\xip{\mc{S}}(G) = \xip{\mc{S}}(N)$, so $\xip{\mc{S}}(G) = (\xip{\mc{S}}(N) - 1) + \xip{\mc{S}}(Q)$ in this case.

Let $\alpha \ge 1$.  Suppose that the result holds whenever $\xip{\mc{S}}(Q) \le \alpha$, and suppose that we are now in the situation that $\xip{\mc{S}}(Q) = \alpha+1$.  Consider $O_1,O_2 \in \Oc(G)$ such that $O_1 \ll O_2$.  Then by Theorem~\ref{thm:reg_quotient}, either $\pi(O_2)$ is compact or $\pi(O_1) \ll \pi(O_2)$.  If $\pi(O_2)$ is compact, then we see that $\xip{\mc{S}}(O_2N) = \xip{\mc{S}}(N)$.  If instead $\pi(O_1) \ll \pi(O_2)$, then by Lemma~\ref{lem:rank_inequalities}(iv), $\xip{\mc{S}}(\pi(O_1)) \le \alpha$, and hence by the inductive hypothesis, $\xip{\mc{S}}(O_1) \le (\xip{\mc{S}}(N) - 1) + \alpha$.

There are now two cases, according to whether $\alpha$ is a successor ordinal or a limit ordinal.

If $\alpha$ is a successor ordinal, by Lemma~\ref{lem:step}(i) we have
\[
\xip{\mc{S}}(G) \le (((\xip{\mc{S}}(N)-1) + \alpha)-1)+2 = (\xip{\mc{S}}(N) - 1) + \alpha + 1,
\]
that is, $\xip{\mc{S}}(G) \le (\xip{\mc{S}}(N)-1) + \xip{\mc{S}}(Q)$.

If $\alpha$ is a limit ordinal, we have $\rho_G(O_1) < (\xip{\mc{S}}(N) - 1) + \alpha$ for all $O_1,O_2 \in \Oc(G)$ such that $O_1 \ll O_2$.  Suppose for a contradiction that there is $H \in \Oc(G)$ such that $\rho_G(H) = (\xip{\mc{S}}(N)-1) + \alpha+1$.  In other words, we have
\[
\sup\{\rho_G(O_1) \mid O_1 \in \Oc(G), O_1 \ll H\} = (\xip{\mc{S}}(N)-1) + \alpha.
\]
By the inductive hypothesis, we see that this can only happen if
\[
\sup\{\xip{\mc{S}}(O_1N/N) \mid O_1 \in \Oc(G), O_1 \ll H\} \ge \alpha.
\]
In particular, we see that $HN/N$ is not compact, so by by Theorem~\ref{thm:reg_quotient} we have $O_1 \ll H \Rightarrow O_1/N \ll HN/N$, and hence
\[
\sup\{\xip{\mc{S}}(O/N) \mid O/N \in \Oc(Q), O/N \ll HN/N\} \ge \alpha.
\]
It follows that $\rho_{Q}(HN/N) \ge \alpha+1$, and then $\xip{\mc{S}}(Q) \ge \alpha+2$, a contradiction.  From this contradiction we conclude that $\rho_G(H) \le (\xip{\mc{S}}(N)-1) + \alpha$ for all $H \in \Oc(G)$, and hence $\xip{\mc{S}}(G) \le (\xip{\mc{S}}(N)-1) + \xip{\mc{S}}(Q)$.  This completes the proof of the theorem.
\end{proof}

We note a special case of Theorem~\ref{thm:extensions}(ii), given Lemma~\ref{lem:rank_one}.

\begin{cor}\label{cor:extensions}
Let $G$ be a first-countable \tdlc group and let $N$ be a closed normal subgroup of $G$ such that $N \le \RadRE(G)$.  Then $G \in \Es{\mc{S}}$ if and only if $G/N \in \Es{\mc{S}}$.  If $G \in \Es{\mc{S}}$, then $\xip{\mc{S}}(G) = \xip{\mc{S}}(G/N)$.
\end{cor}

We can now prove the characterization of $\Es{\Sclass}$ given in Theorem~\ref{intro:Esp_decomp}.  In fact we can prove the analogous result for the class $\Es{\mc{S}}$.

\begin{thm}\label{thm:Esp_decomp}
Write $\mc{S}^* = (\Es{\mc{S}} \cap \mc{S}) \smallsetminus \Es{}$ and let $\mc{C}$ be the smallest class of \tdlc groups such that
\begin{enumerate}[(i)]
\item $\mc{C}$ contains $\mc{S}^*$, the discrete groups and the first-countable profinite groups; and 
\item $\mc{C}$ is closed under extensions that result in a \tdlc group and under directed unions of open subgroups.
\end{enumerate}
Then $\mc{C} = \Es{\mc{S}}$.
\end{thm}

\begin{proof}
It is clear that $\Es{\mc{S}}$ contains $\mc{S}^*$, the discrete groups and the first-countable profinite groups.  By Theorem~\ref{thm:extensions}, $\Es{\mc{S}}$ is closed under extensions, and by Lemma~\ref{lem:rank_inequalities}(iii), $\Es{\mc{S}}$ is closed under directed unions of open subgroups.  Thus $\mc{C} \subseteq \Es{\mc{S}}$.  We note also that $\mc{C}$ contains all first-countable SIN \tdlc groups, since they are compact-by-discrete.  To finish the proof, it is enough to show, given $G \in \Es{\mc{S}}$, that $G \in \mc{C}$: we proceed by induction on $\xip{\mc{S}}(G)$.  Given Theorem~\ref{thm:extensions} and the fact that $\mc{C}$ is closed under extensions, it will not disrupt the induction to pass from $G$ to some quotient of $G$, as long as we know the kernel belongs to $\mc{C}$.

Since $\mc{C}$ is closed under directed unions of open subgroups, it is enough to show $H \in \mc{C}$ for all $H \in \Oc(G)$.  So we may assume $G$ is compactly generated; clearly we may also assume that $G$ is not compact.  Form the characteristic subgroups $R_0,\dots,R_n$ as in Theorem~\ref{thm:rank_reduction}.  Let $\mc{M}$ be the set of $N \unlhd R_n$ such that $R_n/N \in \mc{S} \smallsetminus \Es{}$ and $|G:\N_G(N)|<\infty$, and let $G_0 = \bigcap_{N \in \mc{M}}N$.  Note that $\mc{M}$ is finite by Lemma~\ref{lem:S_quotient_count}, and after replacing $G$ with a finite index subgroup, we may assume that every $N \in \mc{M}$ is normal in $G$.  Then $\xip{\mc{S}}(G_0) < \xip{\mc{S}}(G)$, so by the inductive hypothesis, $G_0 \in \mc{C}$.  It is now enough to show $G/G_0 \in \mc{C}$; from now on we may assume without loss of generality that $G_0 = \triv$.  Inside $R_n$, there is a closed normal subgroup $Q$, which is the quasi-product of the lowermost representatives of the chief factors of $G$ represented by elements of $\{R_n/N \mid N \in \mc{M}\}$.  Using property $\propS$, we see that $Q$ admits a $G$-invariant series whose factors are all in $\mc{S}$; since they are normal factors of $G$ and $G \in \Es{\mc{S}}$, in fact the factors all belong to $\mc{S}^*$.  Thus $Q \in \mc{C}$.  Turning to the quotient $G/Q$, we see that $R_n = \Res{\omega}(G/Q)$, but $R_n/Q$ has no quotients in $\mc{S} \smallsetminus \Es{}$ for which the kernel is virtually normal in $G/Q$.  Applying Theorem~\ref{thm:rank_reduction} to $G/Q$, we see that either $G/Q$ is compact, or else $\xip{\mc{S}}(R_n/Q) < \xip{\mc{S}}(G/Q)$; thus $R_n/Q \in \mc{C}$.  For $1 \le i \le n$, the group $R_{i-1}/R_i$ is a directed union of compactly generated open subgroups, each of which is a SIN group; thus $R_{i-1}/R_i \in \mc{C}$.  Finally, since $\mc{C}$ is closed under extensions that result in a \tdlc group, we conclude that $G \in \mc{C}$, as required.
\end{proof}

We obtain sharper control over the rank in Theorem~\ref{thm:extensions}(i) in the case that $\phi$ is a regionally normal compression.

\begin{prop}\label{prop:normal_compression:rank}
Let $G$ and $H$ be first-countable \tdlc groups and let $\phi: G \rightarrow H$ be a regionally normal compression.  Then $G \in \Es{\mc{S}}$ if and only if $H \in \Es{\mc{S}}$.  If $G \in \Es{\mc{S}}$, then
\[
\xip{\mc{S}}(H) \le \xip{\mc{S}}(G) \le 1+\xip{\mc{S}}(H).
\]
\end{prop}

\begin{proof}
Let $O \in \Oc(G)$.  Note that $\ol{\phi(O)}$ is a locally normal, in particular RIO, subgroup of $H$, so $\phi$ restricts to a normal compression from $O$ to $\ol{\phi(O)}$, and by Lemma~\ref{lem:rank_inequalities}(i), $\xip{\mc{S}}(\ol{\phi(O)}) \le \xip{\mc{S}}(H)$.  Moreover, there is a reduced envelope $E$ of $\ol{\phi(O)}$ in $H$ such that $\ol{\phi(O)}$ is a cocompact normal subgroup of $E$; we then have $\xip{\mc{S}}(E) = \xip{\mc{S}}(\ol{\phi(O)})$.

By Theorem~\ref{thm:extensions}(i), we already know that if $G \in \Es{\mc{S}}$, then $H \in \Es{\mc{S}}$ and $\xip{\mc{S}}(H) \le \xip{\mc{S}}(G)$.  So let us suppose that $H \in \Es{\mc{S}}$; we aim to show $G \in \Es{\mc{S}}$ with $\xip{\mc{S}}(G) \le 1+\xip{\mc{S}}(H)$, and we proceed by induction on $\xip{\mc{S}}(H)$.

In the base case, $\xip{\mc{S}}(H) = 1$ and $H$ is regionally elliptic.  Given $O \in \Oc(G)$, then $\ol{\phi(O)}$ is compact, hence profinite.  It follows that $O$ is residually discrete; in particular, $\Res{\omega}(O)=\triv$.  We deduce that $G \in \Es{\mc{S}}(2)$.

Now suppose $\xip{\mc{S}}(H) = \alpha+1$ where $\alpha \ge 1$.  Let $O_1,O_2 \in \Oc(G)$ such that $O_1 \ll O_2$, and let $O'_i = \ol{\phi(O_i)}$.  Then by Theorem~\ref{thm:reg_quotient}(i), either $O'_2$ is compact or $O'_1 \ll O'_2$.  In the former case, $O_2$ is residually discrete, so $O_1$ is compact and hence $\xip{\mc{S}}(O'_1) = \xip{\mc{S}}(O_1)=1$.  Otherwise, we have $\xip{\mc{S}}(O'_1) \le \alpha$; by the inductive hypothesis, 
\[
\xip{\mc{S}}(O'_1) \le \xip{\mc{S}}(O_1) \le 1+\xip{\mc{S}}(O'_1).
\]

If $\alpha$ is a successor ordinal, then Lemma~\ref{lem:step}(i) implies that $G \in \Es{\mc{S}}(1+\alpha+1)$.  So we may suppose that $\alpha$ is a limit ordinal.  Suppose for a contradiction that $\xip{\mc{S}}(G) > \alpha+1$.  Then by Lemma~\ref{lem:step}(ii) there is $O \in \Oc(G)$ such that 
\[
\sup\{\xip{\mc{S}}(O_1) \mid O_1 \in \Oc(G), O_1 \ll O\} = \alpha.
\]
Note that in this case $O$ cannot be residually discrete, so $\ol{\phi(O)}$ is not compact.  By the inductive hypothesis and the fact that $\alpha$ is an infinite limit ordinal, we see that also
\[
\sup\{\xip{\mc{S}}(\ol{\phi(O_1)}) \mid O_1 \in \Oc(G), O_1 \ll O\} = \alpha.
\]
Now take a reduced envelope $E$ of $\ol{\phi(O)}$ in $H$.  Then for each $O_1 \in \Oc(G)$ such that $O_1 \ll O$, there is a reduced envelope $E_1$ of $\ol{\phi(O_1)}$ in $H$ such that $E_1 \ll E$, and we have $\xip{\mc{S}}(\ol{\phi(O_1)}) = \xip{\mc{S}}(E_1)$.  Thus
\[
\sup\{\xip{\mc{S}}(E_1) \mid E_1 \in \Oc(H), E_1 \ll E\} \ge \alpha.
\]
It follows that $\rho_H(E) \ge \alpha+1$, so $\xip{\mc{S}}(H) \ge \alpha+2$, a contradiction.  From this contradiction, we conclude that $\xip{\mc{S}}(G) \le \xip{\mc{S}}(H)$ as required.
\end{proof}

We now consider how the rank of $G$ can be bounded above in terms of quotients of $G$.

\begin{lem}\label{lem:subdirect}
Let $G$ be a first-countable \tdlc group.  Suppose $A_1,\dots,A_n$ are closed normal subgroups of $G$ with $\bigcap^n_{i=1}A_i = \triv$, such that $G/A_i \in \Es{\mc{S}}$ for $1 \le i \le n$, and let $\alpha+1 = \max\{\xip{\mc{S}}(G/A_i) \mid 1 \le i \le n\}$.  Then $G \in \Es{\mc{S}}$ with $\alpha+1 \le \xip{\mc{S}}(G) \le \alpha + 2$.  If the natural homomorphism from $G$ to $\prod^n_{i=1}G/A_i$ is a closed embedding, then $\xip{\mc{S}}(G) = \alpha + 1$.
\end{lem}

\begin{proof}
We proceed by induction on $\alpha$.  If $G \in \Es{\mc{S}}$, the lower bound $\alpha+1 \le \xip{\mc{S}}(G)$ is clear from Theorem~\ref{thm:extensions}, so it suffices to prove the upper bounds on $\xip{\mc{S}}(G)$.

Suppose $\alpha = 0$, that is, $G/A_i$ is regionally elliptic for $1 \le i \le n$.  Then given $H \in \Oc(G)$ and $1 \le i \le n$, we have $HA_i/A_i \cong H/(H \cap A_i)$, so $H/(H \cap A_i)$ is compact, and hence $\Res{}(H) \le A_i$.  Hence $\Res{}(H) = \triv$, so $\xip{\mc{S}}(H) \le \xi(H) \le 2$, and hence $\xip{\mc{S}}(G) \le 2$.  Now consider the special case that $G$ is embedded as a closed subgroup of $\prod^n_{i=1}G/A_i$.  In this case, we notice that $\prod^n_{i=1}G/A_i$ is regionally elliptic, so its closed subgroup $G$ is also regionally elliptic and hence $\xip{\mc{S}}(G)=1$.

Now suppose $\alpha \ge 1$ and consider $O_1,O_2 \in \Oc(G)$ such that $O_1 \ll O_2$.  Then by Theorem~\ref{thm:reg_quotient}, for each $i$, either $O_2A_i/A_i$ is compact or $O_1A_i/A_i \ll O_2A_i/A_i$.  If $O_2A_i/A_i$ is compact, then $O_1A_i/A_i$ is compact, so $\xip{\mc{S}}(O_1A_i/A_i) = 1$.  If $O_1A_i/A_i \ll O_2A_i/A_i$ then $\xip{\mc{S}}(O_1A_i/A_i) \le \alpha$ by Lemma~\ref{lem:rank_inequalities}(iv).  In either case, $\xip{\mc{S}}(O_1A_i/A_i) \le \alpha$, for $1 \le i \le n$; that is, there is an ordinal $\beta$ with $\beta+1 \le \alpha$ such that 
\[
\max\{\xip{\mc{S}}(O_1A_i/A_i) \mid 1 \le i \le n\} = \beta+1.
\]
Thus by the inductive hypothesis and the fact that $O_1A_i/A_i \cong O_1/(O_1 \cap A_i)$ as topological groups, we have $\xip{\mc{S}}(O_1) \le \beta+2 \le \alpha+1$.  We deduce by Lemma~\ref{lem:step}(i) that $\xip{\mc{S}}(G) \le \alpha+2$.  In the case that $G$ embeds as a closed subgroup of $\prod^n_{i=1}G/A_i$, then $O_1$ embeds as a closed subgroup of $\prod^n_{i=1}O_1/(O_1 \cap A_i)$, and the inductive hypothesis tells us that $\xip{\mc{S}}(O_1) \le \beta+1$; we deduce by Lemma~\ref{lem:step}(i) that $\xip{\mc{S}}(G) \le \beta+2$, so $\xip{\mc{S}}(G) \le \alpha+1$.
\end{proof}

\begin{prop}\label{prop:residual}
Let $G$ be a first-countable \tdlc group and let $\mc{N}$ be a family of closed normal subgroups of $G$.  Suppose that $\bigcap_{N \in \mc{N}}N = \triv$ and that $G/N \in \Es{\mc{S}}$ for all $N \in \mc{N}$.  Let $\alpha = \sup\{\xip{\mc{S}}(G/N) \mid N \in \mc{N}\}$.
\begin{enumerate}[(i)]
\item We have $G \in \Es{\mc{S}}$ and
\[
\alpha \le \xip{\mc{S}}(G) \le 1+ \alpha+1.
\]
\item If $\mc{N}$ is filtering, then $\xip{\mc{S}}(G) \le 1+\alpha^+$.
\item If $\mc{N}$ is filtering and $G$ is compactly generated, then there is $M \in \mc{N}$ such that
\[
\xip{\mc{S}}(G/M) \le \xip{\mc{S}}(G) \le 1 + \xip{\mc{S}}(G/M).
\]
\end{enumerate}
\end{prop}

\begin{proof}
If $G \in \Es{\mc{S}}$, then $\alpha \le \xip{\mc{S}}(G)$ by Theorem~\ref{thm:extensions}.  So it suffices to prove the upper bounds on $\xip{\mc{S}}(G)$.

The upper bounds we aim to show on $\xip{\mc{S}}(G)$ are all successor ordinals.  To show $\xip{\mc{S}}(G)$ is at most any given successor ordinal $\gamma+1$, it is enough to show $\xip{\mc{S}}(H) \le \gamma+1$ for all $H \in \Oc(G)$.  Moreover, given $H \in \Oc(G)$, we have a family $\{N \cap H \mid N \in \mc{N}\}$ of closed normal subgroups of $H$ with trivial intersection such that
\[
 \sup\{\xip{\mc{S}}(H/(N \cap H)) \mid N \in \mc{N}\} =  \sup\{\xip{\mc{S}}(HN/N) \mid N \in \mc{N}\} \le \alpha.
\]
Thus we may assume $G$ is compactly generated.

Since $\mc{N}$ has trivial intersection, it follows by \cite[Theorem~3.3]{RW-EC} that there exist $M_1,\dots,M_n \in \mc{N}$, $M := \bigcap^n_{i=1}M_i$ and a compact normal subgroup $K$ of $G$ such that $K \le M$ and $M/K$ is discrete.  By Lemma~\ref{lem:subdirect}, we see that $\xip{\mc{S}}(G/M) \le \alpha+1$.  If $\mc{N}$ is filtering, in fact there is $M \in \mc{N}$ and a compact normal subgroup $K$ of $G$ such that $K \le M$ and $M/K$ is discrete, so that $\xip{\mc{S}}(G/M) \le \alpha$.  In either case, $M \in \Es{\mc{S}}(2)$; we can apply Theorem~\ref{thm:extensions} to conclude that $G \in \Es{\mc{S}}$ and 
\[
\xip{\mc{S}}(G) \le 1 + \xip{\mc{S}}(G/M) \le 1+\alpha+1. \qedhere
\]
\end{proof}

Proposition~\ref{prop:residual} allows us to define the $\mc{S}$-well-founded residual of a first-countable \tdlc group.

\begin{cor}\label{cor:residual}
Let $G$ be a first-countable \tdlc group.  Then there is a unique smallest closed normal subgroup $\Res{\Es{\mc{S}}}(G)$ such that $G/\Res{\Es{\mc{S}}}(G) \in \Es{\mc{S}}$.  Moreover, $\Res{\Es{\mc{S}}}(G)$ has no nontrivial $\mc{S}$-well-founded quotients.
\end{cor}

\begin{proof}
It is clear from Proposition~\ref{prop:residual} that $R = \Res{\Es{\mc{S}}}(G) := \{N \unlhd G \mid G/N \in \Es{\mc{S}}\}$ is the unique smallest closed normal subgroup $R$ of $G$ such that $G/R \in \Es{\mc{S}}$.  We can then form $\Res{\Es{\mc{S}}}(R)$, which is a topologically characteristic subgroup of $R$, hence a closed normal subgroup of $G$.  By Theorem~\ref{thm:extensions}(ii), $G/\Res{\Es{\mc{S}}}(R) \in \Es{\mc{S}}$, and hence $\Res{\Es{\mc{S}}}(R) = R$ by the minimality of $R$.  In other words, $R$ has no nontrivial $\mc{S}$-well-founded quotients.
\end{proof}

We can also now compute the rank of a \tdlcsc quasi-product.

\begin{cor}\label{cor:quasiproduct}
Let $G$ be a first-countable \tdlc group.  Suppose that $G$ has a quasi-direct factorization $\mc{Q}$.  Then $G \in \Es{\mc{S}}$ if and only if $N \in \Es{\mc{S}}$ for all $N \in \mc{Q}$.  If $G \in \Es{\mc{S}}$, then 
\[
\xip{\mc{S}}(G) = {\sup_{N \in \mc{Q}}}^+\xip{\mc{S}}(N).
\]
\end{cor}

\begin{proof}
First, note that if $G \in \Es{\mc{S}}$, then $N \in \Es{\mc{S}}$ and $\xip{\mc{S}}(N) \le \xip{\mc{S}}(G)$ for all $N \in \mc{Q}$ by Lemma~\ref{lem:rank_inequalities}(i).  So we may assume that $N \in \Es{\mc{S}}$ for all $N \in \mc{Q}$, and we aim to show that $G \in \Es{\mc{S}}(\alpha)$, where $\alpha = {\sup_{N \in \mc{Q}}}^+\xip{\mc{S}}(N)$.

Let $\mc{F}$ be the set of finite subsets of $\mc{Q}$, and given $F \in \mc{F}$ write $G_F = \cgrp{N \in F}$.  Then $\{G_F \mid F \in \mc{F}\}$ is a family of closed normal subgroups of $G$, directed under inclusion, with dense union in $G$.  Given Lemma~\ref{lem:rank_inequalities}(iii), it now suffices to show
\[
G_F \in \Es{\mc{S}} \text{ with } \xip{\mc{S}}(G_F) \le \alpha.
\]
In other words, we may assume that $\mc{Q}$ is finite.  Under this assumption, the direct product $P = \prod_{N \in \mc{Q}}N$ is a first-countable \tdlc group; we then have a natural injective homomorphism $\phi: P \rightarrow G$.  Given $N \in \mc{Q}$ and $O \in \Oc(N)$, we see that $O \le B$ for some second-countable open subgroup of $G$; by considering the action of $B$ on $B \cap N$, we then see that $\N_B(O)$ is open, so $\N_G(O)$ is open.  From here, we see via Lemma~\ref{lem:sigma-normal} that $\phi$ is a regionally normal compression.  By Lemma~\ref{lem:subdirect}, we have $\xip{\mc{S}}(P) = \max\{\xip{\mc{S}}(N) \mid N \in \mc{Q}\} \le \alpha$; then by Theorem~\ref{thm:extensions}(i), it follows that $\xip{\mc{S}}(G) \le \alpha$.
\end{proof}

To conclude this subsection, here is a basic construction we can use, starting from compactly generated groups in $\Es{\mc{S}}$, to produce a group in $\Es{\mc{S}}$ of larger rank.

\begin{prop}\label{prop:wreath_rank}
Let $H$ and $K$ be compactly generated \tdlcsc groups in $\Es{\mc{S}}$, such that $K$ is a closed subgroup of $\mathrm{Sym}(X)$ for some countably infinite set $X$, and let $U$ be a compact open subgroup of $H$.  Form the semidirect product
\[
G = \bigoplus_X (H,U) \rtimes K,
\]
where $K$ acts on the local direct product by permuting the copies of $H$.  Suppose that $K$ acts transitively and that $\Res{\omega}(K) \neq \triv$.  Then $\xip{\mc{S}}(G) \ge \xip{\mc{S}}(H)+1$.
\end{prop}

For the proof we appeal to a standard double commutator argument.

\begin{lem}[{\cite[Lemma~6.8]{CRW-Part2}}]\label{lem:commutator}
Let $A, N$ be subgroups of a group $G$.  Assume $N$ is normal, and contains an element $t \in N$ such that $[A,tAt\inv]=\triv$.  Then $[A,A]\le N$.
\end{lem}

\begin{proof}[Proof of Proposition~\ref{prop:wreath_rank}]
Since $\Res{\omega}(K)$ is normal in $K$ and $K$ acts transitively on $X$, we see that $\Res{\omega}(K)$ acts without global fixed points on $X$.  The hypotheses also ensure that $G$ is compactly generated, for example $G = \grp{\prod_X U , H_x, K}$ where $x \in X$ and $H_x$ is the copy of $H$ in the $x$ coordinate of the base group.

Let $\pi$ be the natural projection map of $G$ onto $K$.  By Lemma~\ref{lem:RIO_sigma-normal}(ii), $\pi(\Res{\omega}(G))$ is dense in $\Res{\omega}(K)$, so $\Res{\omega}(G)$ also induces a permutation group on $X$ with no global fixed points.  Let $B = \Res{\omega}(G) \cap \bigoplus_X(H,U)$; Lemma~\ref{lem:commutator} shows that $B$ contains the derived group of each local direct summand of $\bigoplus_X(H,U)$.

Consider now a closed normal subgroup $N$ of $\Res{\omega}(G)$ such that $\Res{\omega}(G)/N \in \mc{S} \smallsetminus \Es{}$.  If $B \nleq N$, then $\ol{BN} = \Res{\omega}(G)$, and in particular, $\pi(N)$ is dense in $\Res{\omega}(K)$.  But then $N$ would have no global fixed points acting on $X$, so by Lemma~\ref{lem:commutator} applied to $N$ as a subgroup of $BN$, we see that $N$ would have to contain the second derived group of each local direct summand of $\bigoplus_X(H,U)$.  In that case we see that $B/(B \cap N)$ is soluble, which is incompatible with $\ol{BN}/N$ being nonabelian and topologically simple.  From this contradiction, we see that $B \le N$.

From the previous paragraph, we see that given any $C \in \Oc(B)$, then $C \ll G$ and hence by Lemma~\ref{lem:rank_inequalities}(iv), we have $\xip{\mc{S}}(G) \ge \xip{\mc{S}}(C)+1$.  It remains to find $C \in \Oc(B)$ such that $\xip{\mc{S}}(C) \ge \xip{\mc{S}}(H)$.  Let $\rho$ be the projection of $\bigoplus_X(H,U)$ onto $H_x$.  Let $r \in \Res{\omega}(G)$ be such that $r.x \neq x$ and fix a compact generating set $D$ for $H$; indeed we can take $D = F \cup U$ where $F$ is finite.  Then we see that for each $h \in D$, there is an element $g_h = ((g_h)_x)_{x \in X}$ of $B$ such that $(g_h)_x = h$ and $(g_h)_{r.x} = h\inv$; moreover, given the form of $D$, we can clearly choose the set $\{g_h \mid h \in D\}$ to have compact closure, which means it is contained in some $C \in \Oc(B)$.  In particular, $\rho(C) = H_x$.  By Theorem~\ref{thm:extensions} we deduce that $\xip{\mc{S}}(C) \ge \xip{\mc{S}}(H)$, and hence $\xip{\mc{S}}(G) \ge \xip{\mc{S}}(H)+1$.
\end{proof}

\subsection{A construction preserving the class of $\mc{S}$-well-founded groups}\label{sec:ended_tree}

So far, we have shown that the class $\Es{\mc{S}}$ is stable under various basic constructions.  We now recall a construction given in \cite[\S9]{RW-DenseLC}, which can be used to produce more complicated examples; see also \cite[\S10.2]{RW-Polish}.  Let $X$ be a countably infinite set with a distinguished element $0$, let $T$ be the regular tree of countably infinite degree and fix an end $\delta$ of $T$.  Given an arc $a$ we write $o(a)$ for the origin of $a$ and $t(a)$ for the terminus of $a$.  Let $A$ be the set of arcs of $T$ pointing towards $\delta$ and given $v \in VT$, write $A_v$ for the set of arcs in $A$ terminating at $v$.  Choose a function $c: A \rightarrow X$, called an \defbold{ended colouring}, with the following properties:
\begin{enumerate}[(i)]
\item For every $v \in VT$, then $c$ restricts to a bijection $c_v$ from $A_v$ to $X$;
\item There is a ray $(a_1,a_2,a_3,\dots)$ of arcs in $A$, representing the end $\delta$, such that $c(a_i) = 0$ for all $i$.
\end{enumerate}
Given an automorphism $g$ of $T$ fixing $\delta$ and $v \in VT$, we now define the \defbold{local action of $g$ at $v$} to be
\[
\sigma(g,v) = c_{gv} g c\inv_v \in \Sym(X).
\]

We recall a few standard facts about group actions on a tree fixing an end; all the claims in this lemma are easily verified.

\begin{lem}\label{lem:Busemann}
Let $G$ be a group acting on a tree $T$, fixing an end $\delta$.  Then there is a function $b: VT \rightarrow \Zb$, which is unique up to a constant, such that whenever $(v,w)$ is an arc pointing towards $\delta$ then $b(w)-b(v)=1$.  There is then a homomorphism $\beta: G \rightarrow \Zb$, which does not depend on the choice of $b$, such that $\beta(g) = b(g.v)-b(v)$ for all $g \in G$ and $v \in VT$.  Moreover, every finitely generated subgroup of $\ker\beta$ fixes pointwise a ray representing $\delta$, while every element of $G \smallsetminus \ker\beta$ is a translation.
\end{lem}

From now on we let $\beta$ be the function $\beta: \Aut(T)_{\delta} \rightarrow \Zb$ as in Lemma~\ref{lem:Busemann}.  The conditions imposed on $c$ ensure that there is a monochromatic doubly infinite path in $T$ with one end at $\delta$, and hence there is a \defbold{basic translation} $s \in \Aut(T)_{\delta}$ such that $\beta(s)=1$ and $\sigma(s,v)=\mathrm{id}_X$ for all $v \in VT$.  In fact, $s$ is uniquely specified by these conditions.

Now suppose that $G$ is a locally compact subgroup of $\Sym(X)$ and $U$ is a compact open subgroup of $G$.  We define $E_X(G,U)$ to be the set of all $g \in \Aut(T)_{\delta}$ such that $\sigma(g,v) \in G$ for all $v \in VT$ and $\sigma(g,v) \in U$ for all but finitely many $v \in VT$.  It is straightforward to check that $E_X(G,U)$ is a group.  Consider the subgroup $E_X(U,U)_{v_0}$ of $E_X(G,U)$, where $v_0 = o(a_1)$.  Since $U$ is compact, it has finite orbits on $X$; from this fact, we deduce that $E_X(U,U)_{v_0}$ has finite orbits on $VT$, and hence $E_X(U,U)_{v_0}$ is a compact subgroup of $\Aut(T)_{\delta}$ in the permutation topology.  Moreover, it is not hard to show that $E_X(U,U)_{v_0}$ is commensurated by $E_X(G,U)$.  Consequently, there is a unique group topology on $E_X(G,U)$ so that $E_X(U,U)_{v_0}$ is embedded in $E_X(G,U)$ as a compact open subgroup; from now on we will take $E_X(G,U)$ to be equipped with this topology.  One can then check that $E_X(G,U)$ is a \tdlcsc group.  We also define $P_X(G,U) = E_X(G,U) \cap \ker(\beta)$.  Note that $E_X(G,U)$ contains the basic translation $s$, and hence $E_X(G,U) = P_X(G,U) \rtimes \grp{s}$; also, $E_X(U,U)_{v_0} \le P_X(G,U)$, so $P_X(G,U)$ is open in $E_X(G,U)$.

We now bound the $\mc{S}$-well-founded rank of $E_X(G,U)$ in terms of that of $G$.

\begin{prop}\label{prop:ended_construction}
Let $E_X(G,U)$ be constructed as above.  Then $E_X(G,U) \in \Es{\mc{S}}$ if and only if $G \in \Es{\mc{S}}$.  If $G \in \Es{\mc{S}}$, then
\[
\xip{\mc{S}}(G) \le \xip{\mc{S}}(E_X(G,U)) \le \xip{\mc{S}}(G).\omega + 2.
\]
If $G$ is a compactly generated group in $\Es{\mc{S}}$ acting transitively on $X$ such that $\Res{\omega}(G) \neq \triv$, then
\[
\xip{\mc{S}}(G) + \omega + 1 \le \xip{\mc{S}}(E_X(G,U)).
\]
\end{prop}

\begin{proof}
Note that given Lemma~\ref{lem:Busemann}, we see that every compactly generated subgroup of $P_X(G,U)$ fixes pointwise a ray representing $\delta$.  In particular we can tell whether or not $P_X(G,U) \in \Es{\mc{S}}$, and also bound the rank, by considering vertex stabilizers.

Define the \defbold{horoballs} $X_i = \{v \in VT \mid \beta(v) \ge i\}$ and \defbold{horospheres} $Y_i = \{v \in VT \mid \beta(v) = i\}$ for $i \in \Zb$.  Note that $sX_i = X_{i+1}$ and $sY_i = Y_{i+1}$ for all $i \in \Zb$.  We see that $P_X(G,U)$ stabilizes each horosphere setwise, and so there is an action of $P_X(G,U)$ on $Y_i$ for each $i \in \Zb$.  Let $P_i$ be the set of $g \in P_X(G,U)$ such that $\sigma(g,v) \in U$ for all $v \in X_i$ and let $N_i$ be the set of $g \in P_X(G,U)$ such that $\sigma(g,v) \in U$ for all $v \in VT \smallsetminus X_i$.  Since $X_i$ is preserved by $P_X(G,U)$, we see that $P_i$ and $N_i$ are subgroups of $P_X(G,U)$; in fact they are open subgroups, since they each contain a conjugate of $E_X(U,U)_{v_0}$.  In addition, $P_i \le P_{i+1}$ and $N_i \ge N_{i+1}$ for all $i \in \Zb$.  Since any finite set of vertices is contained in the difference of some pair of horoballs, we see that $P_X(G,U) = \bigcup_{i \in \Zb}P_i = \bigcup_{i \in \Zb}N_i$.  In particular, every compactly generated open subgroup of $P_X(G,U)$ is contained in $P_i \cap N_j$ for some $i$ and $j$; more specifically, one can take $i \ge j$ here.

Fix integers $i \ge j$ and consider $R = P_i \cap N_j$.  We can decompose $R$ as follows: let $R_n$ be the set of $g \in R$ such that $\sigma(g,v) = \mathrm{id}$ for all $v \in X_n$, and let $Q_n$ be the set of $g \in R$ such that $\sigma(g,v) = \mathrm{id}$ for all $v \in VT \smallsetminus X_n$.  Note that $R_n$ can equivalently be defined as the set of elements of $R$ that fix $X_{n-1}$ pointwise.  We see that $R = R_i \rtimes Q_i$, and in turn $R_i = R_j \rtimes S$ where $S = R_i \cap Q_j$.  Given $g \in Q_i \cup R_j$, we see that $\sigma(g,v) \in U$ for all $v \in VT$, that is, $Q_i \cup R_j \subseteq P_X(U,U)$.  Since every compactly generated subgroup fixes a vertex, and since vertex stabilizers in $P_X(U,U)$ are compact, we see that $Q_i,R_j \in \Es{\mc{S}}(1)$.  Hence by Theorem~\ref{thm:extensions}, we have $R \in \Es{\mc{S}}$ if and only if $S \in \Es{\mc{S}}$, and if $R \in \Es{\mc{S}}$ then $\xip{\mc{S}}(R) = \xip{\mc{S}}(S)$.

Finally, consider the structure of $S$.  We have a faithful action of $S$ on $X_j \smallsetminus X_i$.  In the case that $i = j+1$, we see that $S$ is actually isomorphic to a local direct product $L = \bigoplus_{\aleph_0} (G,U)$ with a countably infinite number of copies of $G$.  More generally, $S$ takes the form of an iterated semidirect product of copies of $L$, that is,
\[
S \cong (\dots(L \rtimes L) \rtimes L) \dots \rtimes L),
\]
where there are $j-i$ occurrences of $L$ in the right-hand expression.

We now relate $G$ to $E_X(G,U)$.  If $G \not\in \Es{\mc{S}}$, then $L \not\in \Es{\mc{S}}$ by Lemma~\ref{lem:rank_inequalities}, so $S \not\in \Es{\mc{S}}$ by Theorem~\ref{thm:extensions} and hence $R \not\in \Es{\mc{S}}$; we then recall that $R$ is  an open subgroup of $E_X(G,U)$, so by Lemma~\ref{lem:rank_inequalities} again, $E_X(G,U) \not\in \Es{\mc{S}}$.  Now let us suppose $G \in \Es{\mc{S}}$.  Then by Corollary~\ref{cor:quasiproduct}, $L \in \Es{\mc{S}}$ and $\xip{\mc{S}}(L) = \xip{\mc{S}}(G)$.  Hence by Theorem~\ref{thm:extensions}, $S \in \Es{\mc{S}}$ and
\[
\xip{\mc{S}}(G) \le \xip{\mc{S}}(S) \le \xip{\mc{S}}(G).(j-i) < \xip{\mc{S}}(G).\omega.
\]
In turn, $\xip{\mc{S}}(R) = \xip{\mc{S}}(S)$.  By letting $i$ and $j$ vary, we obtain subgroups containing all compactly generated open subgroups of $P_X(G,U)$, and so for all $H \in \Oc(P_X(G,U))$ we have $H \in \Es{\mc{S}}$ with
\[
\xip{\mc{S}}(H) \le \xip{\mc{S}}(G).n
\]
for some $n < \omega$ depending on $H$.  We deduce that $P_X(G,U) \in \Es{\mc{S}}$ with
\[
\xip{\mc{S}}(G) \le \xip{\mc{S}}(P_X(G,U)) \le \xip{\mc{S}}(G).\omega + 1,
\]
and hence by Theorem~\ref{thm:extensions}, $E_X(G,U) \in \Es{\mc{S}}$ with
\[
\xip{\mc{S}}(G) \le \xip{\mc{S}}(E_X(G,U)) \le \xip{\mc{S}}(G).\omega + 2.
\]

Finally, suppose now that $G$ is compactly generated, $\Res{\omega}(G) \neq \triv$, and $G$ acts transitively on $X$.  Then we can obtain a better lower bound on $\xip{\mc{S}}(S)$, and hence on $\xip{\mc{S}}(E_X(G,U))$.  If $i=j$, write $S_{i,j} = \triv$.  Let us suppose now that $i < j$.  Observe that $S$ has a compactly generated closed normal subgroup of the form
\[
S_{i,j} \cong (\dots(L \rtimes L) \rtimes L) \dots \rtimes L) \rtimes G,
\]
with $j-i-1$ occurrences of $L$.  The group $S_{i,j}$ is formed as follows.  Choose a single vertex $w_0$ of $Y_{i-1}$, and for each $w \in VT$, let $Z_w$ be the set of vertices $v$ of $T$ such that the ray from $v$ to $\delta$ passes through $w$.  We then define $S_{i,j}$ to be the rigid stabilizer of $Z_{w_0}$ in $S$, that is, of all $g \in S$ such that $g$ fixes every $v \in VT \smallsetminus Z_{w_0}$.  To see that $S_{i,j}$ is indeed normal in $S$, note that $S$ fixes $w_0$, so the set $Z_{w_0}$ is $S$-invariant.  There is a natural quotient map $\pi$ from $S_{i,j}$ onto $G$; let $K_{i,j}$ be the kernel of this action.

The kernel $K_{i,j}$ naturally decomposes as a local direct product $\bigoplus_X (S_{i+1,j},V)$ of copies of $S_{i+1,j}$.  The copy of $S_{i+1,j}$ indexed by $x$ is the rigid stabilizer in $K_{i,j}$ of $Z_{w_x}$, where there is an arc from $w_x$ to $w_0$ with colour $x$; the corresponding compact open subgroup $V$ consists of those elements $g$ such that in addition, $\sigma(g,v) \in U$ at every vertex.  We now see that we can rewrite $S_{i,j}$ in the form
\[
S_{i,j} = \bigoplus_X (S_{i+1,j},V) \rtimes G.
\]
Hence by Proposition~\ref{prop:wreath_rank}, we have
\[
\xip{\mc{S}}(S_{i,j}) \ge \xip{\mc{S}}(S_{i+1,j}) + 1.
\]
Repeating this argument, we find that
\[
\xip{\mc{S}}(S_{i,j}) \ge \xip{\mc{S}}(S_{j-1,j}) + (j-1-i) = \xip{\mc{S}}(G) + (j-i-1).
\]
The quantity $j-i-1$ is unbounded below $\omega$, so we conclude that
\[
\xip{\mc{S}}(E_X(G,U)) \ge \xip{\mc{S}}(G) + \omega + 1. \qedhere
\]
\end{proof}

\begin{rmk}
We can iterate the construction in Proposition~\ref{prop:ended_construction} to produce a group in $\Es{\Sclass}$ of rank at least $\omega^2+1$.  Start with some $G_0 \in \Es{\Sclass} \cap \Sclass$ (see below for examples), let $U_0$ be a compact open subgroup of $G_0$ and let $X_0 = G_0/U_0$.  We then produce a sequence $(G_i)_{i < \omega}$ in $\Es{\Sclass}$, each with a transitive permutation action on a countable set $X_i$, so that $G_{i+1} = E_{X_i}(G_i,U_i)$; $U_{i+1}$ is the standard compact open subgroup $E_{X_i}(U_i,U_i)_{v_0}$; and $X_{i+1} = G_{i+1}/U_{i+1}$.  We see then that $G_i$ is $\Sclass$-well-founded for all $i < \omega$ and that $\xip{\Sclass}(G_{i+1}) \ge \xip{\Sclass}(G_i) + \omega + 1$, so $\xip{\Sclass}(G_i) \ge \omega.i + 1$.  Now form the local direct product $H = \bigoplus_{i < \omega}(G_i,U_i)$; then $H \in \Es{\Sclass}$ and $\xip{\Sclass}(H) \ge \omega^2+1$.
\end{rmk}

\subsection{Regionally near-simple classes}

Recall the class $\ms{R}$ of robustly monolithic groups from Section~\ref{sec:simple_def}.  The class $\ms{R}$ does not consist solely of topologically simple groups, but we note that $G \in \ms{R}$ if and only if $\Mon(G) \in \ms{R}$, and $\Mon(G)$ is topologically simple.  By \cite[Proposition~5.1.2]{CRWes}, every group in $\ms{R}$ is regionally expansive, and hence first-countable.  Moreover, being robustly monolithic is a regional property.  This was used in \cite{CRWes} to show that groups in $\ms{R}$ are not elementary, by an infinite descent argument.  A similar argument can potentially be applied to other classes of \tdlc group; for the applications it is useful to expand the notion of ``monolith'' somewhat.

\begin{defn}\label{def:near-simple}
Let $G$ be a \tdlc group.  A \tdlc group is \defbold{regionally elliptic} if every compact subset is contained in a compact subgroup; note that all first-countable regionally elliptic \tdlc groups belong to $\Es{}(2)$.  We define the \defbold{regionally elliptic radical} $\RadRE(G)$ to be the smallest closed subgroup containing all regionally elliptic closed normal subgroups of $G$.  We refer to $\Mon(G/\RadRE(G))$ as the \defbold{near-monolith} of $G$.  Write $\NMon(G)$ for the preimage of $\Mon(G/\RadRE(G))$ in $G$.

Let $\ms{A}$ be a class of first-countable \tdlc groups.  We say $\ms{A}$ is a \defbold{regionally near-simple class} if the following conditions hold.
\begin{enumerate}[(a)]
\item Given $G \in \ms{A}$, then $\Mon(G/\RadRE(G))$ is nondiscrete and topologically simple.
\item Given $G \in \ms{A}$ and given an open subgroup $H$ of $G$ such that $\NMon(G) \le H\RadRE(G)$, then $H \in \ms{A}$.
\item Given $G \in \ms{A}$, then there exists $H \in \Oc(G)$ such that for all $K \in \Oc(G)$ with $K \ge H$, we have $K \in \ms{A}$.
\end{enumerate}
\end{defn}

The class $\Sclass$ is clearly a regionally near-simple class, and it follows from \cite[Proposition~5.1.2 and Theorem~5.2.2]{CRWes} that $\ms{R}$ is a regionally near-simple class.  Here are some features of regionally near-simple classes with respect to well-foundedness.

\begin{lem}\label{lem:near-simple_nonelementary}
Let $\ms{A}$ be a regionally near-simple class and let $G \in \ms{A}$.  Then the near-monolith of $G$ is not regionally elementary.
\end{lem}

\begin{proof}
Let $G \in \ms{A}$ and suppose for contradiction that the near-monolith of $G$ is regionally elementary.  Letting $U$ be a compact open subgroup of $G$, we see that $H = \NMon(G)U$ is also regionally elementary.  From parts (b) and (c) of the definition, there is $G_0 \in \Oc(H) \cap \ms{A}$; now $G_0$ is compactly generated and elementary.  We can then construct an infinite descending sequence of compactly generated open subgroups in $\ms{A}$ as follows: for $i \ge 0$, set $G_{i+1} \in \Oc(\NMon(G_i)U_i)$ such that $U_i$ is a compact open subgroup of $G_i$ and $G_{i+1} \in \ms{A}$.  At each stage we see that $\NMon(G_i) \le \Res{}(G_i)$, since the near-monolith is not discrete, and hence $\xi(G_{i+1}) \le \xi(\Res{}(G_i)) < \xi(G_i)$ for all $i$, which is absurd, since $\xi$ takes ordinal values on elementary groups.  This contradiction finishes the proof.
\end{proof}

\begin{prop}\label{prop:R_descent}
Let $\ms{A}$ be a regionally near-simple class, let $G$ be a compactly generated group in $\ms{A}$, let $U$ be a compact open subgroup of $G$ and let $H = \NMon(G)U$.  Suppose there exists $K \in \Oc(H)$ such that $K \not\ll G$.  Then $\Mon(G/\RadRE(G)) \in \mc{S} \smallsetminus \Es{}$ and $G/\NMon(G) \in \Es{}(\omega)$.
\end{prop}

\begin{proof}
Write $R = \RadRE(G)$ and let $\mc{M}$ be the set of closed normal subgroups $M$ of $\Res{\omega}(G)$ such that $|G:\N_G(M)|<\infty$ and $\Res{\omega}(G)/M \in \mc{S} \smallsetminus \Es{}$.  Note that $\mc{M}$ is finite by Lemma~\ref{lem:S_quotient_count}(ii), and similar to the proof of Lemma~\ref{lem:S_quotient_count}, the elements $M$ of $\mc{M}$ correspond to nonelementary topologically simple minimal closed normal subgroups $K_M/N$ of $\Res{\omega}(G)/N$ where $N = \bigcap_{M \in \mc{M}}M$.  Using Lemma~\ref{lem:commutator} and the fact that $K_M/N \not\in \Es{}$ for each $M \in \mc{M}$, we see that $R$ normalizes each $K_M$; moreover, using the fact that $[R,K_M]$ has regionally elementary closure and $K_M/N$ is topologically simple, we see that $[R,K_M] \le N$ for every $M \in \mc{M}$ and hence $R$ centralizes $\Res{\omega}(G)/N$.  In particular, $R \le G_0$, where $G_0 = \bigcap_{M \in \mc{M}}\N_G(M)$.  Given Lemma~\ref{lem:near-simple_nonelementary}, we see that $\NMon(G)/R \not\in \Es{}$, ensuring that
\[
\NMon(G) \le \ol{\Res{\omega}(G)R} \le G_0,
\]
and from here it is easy to see that $R = \RadRE(G_0)$, $\NMon(G) = \NMon(G_0)$, $\Res{\omega}(G) = \Res{\omega}(G_0)$ and $G_0 \in \ms{A}$.  We can also replace $K$ with $K \cap G_0$ with no effect on the partial order $\ll$.  So from now on we may assume $G = G_0$, that is, every $M \in \mc{M}$ is normal.

We claim that $RU$ is regionally elliptic; it is enough to show every $E \in \Oc(RU)$ containing $U$ is regionally elliptic.  Given such an $E$, then by Lemma~\ref{lem:cocompact_gen} we have $E = \grp{F}U$ where $F$ is a finite subset of $R$.  By the regionally elliptic property, $\grp{F}$ has compact closure; hence $E$ is compact, as required.

Suppose $K \in \Oc(H)$ is such that $K \not\ll G$.  Given Lemma~\ref{lem:near-simple_nonelementary} we see that $\NMon(G) \le \ol{\Res{\omega}(G)R}$, so $K \le \Res{\omega}(G)RU$.  By the previous paragraph, it follows that $K\Res{\omega}(G)/\Res{\omega}(G)$ is compact.  So the only way to have $K \not\ll G$ is if there is a closed normal subgroup $M \in \mc{M}$ such that $\Res{\omega}(K) \nleq M$.  In particular, the image of $\Res{\omega}(K)$ in $\Res{\omega}(G)/M$ contains a noncentral element; since $R$ centralizes $\Res{\omega}(G)/M$ we deduce that $\Res{\omega}(K) \nleq \ol{MR}$.  We also see that $\Res{\omega}(K) \le \NMon(G)$, so in fact $\NMon(G) \nleq \ol{MR}$.  Now $\ol{MR}/R$ is a closed normal subgroup of $G/R$ that does not contain the monolith, so it is trivial, that is, $M \le \RadRE(G)$.  Since $\NMon(G)/M$ and $\Res{\omega}(G)/M$ have nontrivial intersection (namely they both contain $\Res{\omega}(K)M/M$) and $\Res{\omega}(G)/M$ is topologically simple, in fact $\Res{\omega}(G) \le \NMon(G)$; clearly also $\Res{\omega}(G) \nleq R$, so considering the monolith of $G/R$ again, we have
\[
\NMon(G) = \ol{\Res{\omega}(G)R}.
\]
As a result, we have a normal compression from $\Res{\omega}(G)/M$ to $\Mon(G/R)$.  Using property $\propS$, it follows that $\Mon(G/R) \in \mc{S} \smallsetminus \Es{}$.  At the same time, the fact that $\Res{\omega}(G) \le \NMon(G)$ ensures that $G/\NMon(G) \in \Es{}(\omega)$.
\end{proof}

Consider the contrapositive of Proposition~\ref{prop:R_descent}: if $\Mon(G/\RadRE(G)) \not\in \mc{S} \smallsetminus \Es{}$, or if $G/\NMon(G)$ does not have finite elementary rank, then $K \ll G$ for all $K \in \Oc(H)$; from the definition of regionally near-simple class, we see that we can choose $K \in \ms{A}$.  An infinite descent argument then immediately leads to the following conclusion.

\begin{cor}\label{cor:R_descent}
Let $\ms{A}$ be a regionally near-simple class and let $G \in \ms{A} \cap \Es{\mc{S}}$.  Then there is $H \in \Oc(G)$ such that $\Mon(H/\RadRE(H)) \in \mc{S} \smallsetminus \Es{}$ and $H/\NMon(H) \in \Es{}(\omega)$.
\end{cor}

To finish this subsection, we note a condition related to near-simplicity that ensures $\Bigsclass$-well-foundedness.  As we will see later, this condition occurs in some natural families of examples.

\begin{lem}\label{lem:simple_well-founded}
Let $G$ be a first-countable \tdlc group.  Suppose that for every $H \in \Oc(G)$, there are closed normal subgroups $K \le L$ such that $K, H/L \in \Es{}(\omega)$ and $L/K$ is a (possibly trivial) finite direct product of nonelementary topologically simple \tdlc groups.  Then $G \in \Es{\Bigsclass}(3)$.
\end{lem}

\begin{proof}
Write $\rho$ for the rank function on $\Oc(G)$ associated to $\ll_{\Bigsclass}$.  Let $H \in \Oc(G)$ and let $K$ and $L$ be as in the statement.  Then we see that $L/K = \Res{\omega}(H/K)$, so given $H'/K \in \Oc(G/K)$ such that $H'/K \ll_{\Bigsclass} H/K$, then $H'/K \in \Es{}(\omega)$.  By Theorem~\ref{thm:reg_quotient}, and since $K \in \Es{}(\omega)$, we see that whenever $H' \in \Oc(G)$ is such that $H' \ll_{\Bigsclass} H$, then $H' \in \Es{}(\omega)$, meaning that $\rho(H') \le 1$.  From the freedom of choice of $H$ and $H'$, it follows that $\rho(H) \le 2$ and hence $\xip{\Bigsclass}(G) \le 3$.
\end{proof}

\section{Special families of groups with regard to $\Sclass$-well-foundedness}\label{sec:special}

\subsection{Noetherian groups}

in this section we will give some sufficient conditions for a group to be $\Sclass$-well-founded, that take us beyond the elementary case, but also give an example of a family of groups that is $\Bigsclass$-well-founded but not $\Sclass$-well-founded.  Before discussing these conditions, we briefly recall a property related to \emph{ascending} chains in $(\Oc(G),\le)$.

\begin{defn}
A \tdlc group $G$ is \defbold{Noetherian} if every open subgroup of $G$ is compactly generated.  (Note: some authors impose the stronger condition that every closed subgroup should be compactly generated.)
\end{defn}

\begin{lem}
Let $G$ be a \tdlc group.  Then $G$ is Noetherian if and only if there is no infinite ascending chain in $(\Oc(G),\le)$.
\end{lem}

\begin{proof}
Suppose that $G$ is Noetherian; let $(H_i)_{i \in \Nb}$ be an ascending sequence in $\Oc(G)$, and let $H = \bigcup_{i \in \Nb}H_i$.  Then $H$ is an open subgroup of $G$, so it has a compact generating set $X$.  Since $X$ is compact, we must have $X \subseteq H_i$ for some $i$.  But then for all $j \ge i$, we have
\[
H = \grp{X} \le H_i \le H_j \le H,
\]
so $H_j = H_i$; in other words, the sequence $(H_i)_{i \in \Nb}$ stabilizes.

Conversely, suppose that $G$ is not Noetherian, and let $H$ be an open subgroup of $G$ that is not compactly generated.  Then we can produce a strictly ascending sequence $(H_i)_{i \in \Nb}$ in $\Oc(G)$ by taking $H_0$ to be a compact open subgroup of $H$, and thereafter $H_{i+1} = \grp{H_i,x_{i+1}}$ where $x_{i+1} \in H \smallsetminus H_i$; the existence of the elements $x_{i+1}$ is ensured by the fact that $H$ is not compactly generated.
\end{proof}

We also note that in a Noetherian \tdlc group, every locally normal subgroup is cocompact in an open subgroup, hence also compactly generated.

A general structure theorem for Noetherian locally compact groups was proved in \cite{CM}.

\begin{thm}[{See \cite[Theorem~C]{CM}}]\label{thm:CM_Noetherian}
Let $G$ be a Noetherian \tdlc group.  Then there is an open normal subgroup $G_k$ and closed subnormal subgroups
\[
\triv = G_0 \lhd G_1 \lhd \dots \lhd G_k,
\]
such that, for each $1 \le i \le k$, the factor $G_i/G_{i-1}$ is compact, isomorphic to $\Zb$, or in $\Sclass$.
\end{thm}

\begin{cor}\label{cor:CM_Noetherian}
Let $G$ be a Noetherian \tdlcsc group and let $G_0,G_1,\dots,G_k$ be as in Theorem~\ref{thm:CM_Noetherian}.  Let $\mc{F} = \{G_i/G_{i-1} \mid 1 \le i \le k, G_i/G_{i-1} \in \Sclass\}$.
\begin{enumerate}[(i)]
\item $G$ is elementary if and only if $\mc{F}$ is empty.  If $G$ is elementary, then $\xi(G)<\omega$.
\item $G$ is $\Sclass$-well-founded if and only if $\mc{F} \subseteq \Es{\Sclass}$.
\end{enumerate}
\end{cor}

There are certainly compactly generated $\Sclass$-well-founded groups that are not Noetherian.  For instance, one sees that the group $E_X(G,U)$ constructed in Section~\ref{sec:ended_tree} is non-Noetherian provided that $G$ does not fix the distinguished point $0 \in X$: consider the ascending chain of stabilizers $E_X(G,U)_v$ as $v$ ranges over a ray of vertices representing the fixed end.  However, it is not clear whether there are any Noetherian \tdlc groups that are not $\Sclass$-well-founded. Many of the examples below are both Noetherian and $\Sclass$-well-founded, whereas the non-$\Sclass$-well-founded examples we find are also non-Noetherian.

\subsection{Groups acting on trees with Tits' independence property}\label{sec:trees}

\begin{defn}
Let $\Gamma$ be an undirected graph.  We write $V\Gamma$ for the set of vertices of $\Gamma$.  Each edge in $\Gamma$ is considered as a pair of arcs $\{a,\ol{a}\}$, where $a$ has origin $o(a) \in VT$ and terminus $t(a) \in VT$, and the reverse arc $\ol{a}$ has $o(\ol{a}) = t(a)$ and $t(\ol{a}) = o(a)$.  Write $A\Gamma$ for the set of arcs of $\Gamma$.  A \defbold{tree} is then a connected undirected graph with no loops, multiple edges or cycles.

Let $T$ be a tree (not necessarily locally finite) and let $G \le \Aut(T)$.  We define the \defbold{rigid stabilizer} $\rist_G(T')$ of a subgraph $T'$ to be the pointwise fixator of $VT \smallsetminus VT'$.  Given a path $L$ in $T$, let $\pi_L$ be the closest point projection from $T$ to $L$.  Then $G$ has \defbold{property $\propP{}$} (with respect to a class of paths $\mc{L}$) if, for every finite or infinite path $L$ in $T$ (or every $L \in \mc{L}$), the pointwise fixator $G_{(L)}$ of $L$ splits as a direct product 
\[
G_{(L)} = \prod_{v \in VL}\rist_G(\pi\inv_L(v)).
\]
We say the action is \defbold{geometrically dense} if it does not preserve any proper subtree or fix any end of $T$.
\end{defn}

Property $\propP{}$ was introduced by J. Tits in \cite{Tits70}, where it was used to prove the following simplicity theorem.

\begin{thm}[{\cite[Th\'{e}or\`{e}me~4.5]{Tits70}}]\label{thm:Tits}
Let $G$ be a group acting geometrically densely with property $\propP{}$ on a tree $T$.  Let $G^+$ be the subgroup of $G$ generated by the arc stabilizers.  Then every subgroup $H$ of $G$ normalized by $G^+$ contains $G^+$.  In particular, $G^+$ is trivial or abstractly simple.
\end{thm}

Given $v \in VT$, write $G_v$ for the vertex stabilizer; $G_{v,1}$ for the subgroup fixing pointwise the set $o\inv(v)$ of arcs originating at $v$; and $G(v)$ for the \defbold{local action at $v$}, which is the quotient $G_v/G_{v,1}$ regarded as a permutation group acting on $o\inv(v)$.  We may also equivalently regard $G(v)$ as acting on the neighbours of $v$; the correspondence between these actions is given by the bijection $a \mapsto t(a)$ from $o\inv(v)$ to the set of neighbours of $v$.  We will generally be restricting attention to the case that $T$ is countable and $G(v)$ has compact point stabilizers, in order for the group acting on the tree to be a \tdlcsc group with compact arc stabilizers.  Note, however, that property $\propP{}$ by itself does not impose any further restriction on $G(v)$; so up to a compact normal subgroup, $G(v)$ could be any \tdlcsc group.  (For more on why this is the case, see \cite{RS}, where a general construction is given for all groups acting on trees with property $\propP{}$.)  Thus all results in this subsection will be conditional on some hypothesis about $G(v)$, or reducing questions about $G$ to questions about the permutation groups $G(v)$.

If some $G(v)$ is not $\Sclass$-well-founded, then the open subgroup $G_v$ of $G$ is not $\Sclass$-well-founded and hence $G$ is not $\Sclass$-well-founded.  If $G(v)$ is $\Sclass$-well-founded for every $v \in VT$, we will see that sometimes $G$ is $\Sclass$-well-founded and sometimes not, and this is true even when the tree is locally finite.  In any case though, we have quite a good general picture of how the partially ordered set $(\Oc(G),\ll)$ relates to the tree structure, and in turn to the structure of vertex stabilizers.

In this subsection we will make use of a few results from an upcoming revision of the preprint \cite{RS}.  For clarity, the relevant results including proofs are provided here in an appendix.

We say a group $G$ acting on a tree $T$ is of \defbold{general type} if it contains a translation and does not fix any end.  If $G$ is of general type, then there is a unique smallest $G$-invariant subtree $T^*$, which is the union of the axes of translation of $G$ (see \cite[Corollaire~3.5]{Tits70} and \cite[Lemma~2.1(iii)]{MollerVonk}); in this case, we will write $G^{*+}$ for the subgroup generated by the stabilizers of arcs of $T^*$.  The action of $G$ on $T^*$ then has property~$\propP{}$ (see Lemma~\ref{lem:propP_subtree}).  Theorem~\ref{thm:Tits} thus generalizes to actions of general type as follows.

\begin{cor}\label{cor:Tits}Let $G$ be a group acting on a tree $T$ with property $\propP{}$, such that $G$ is of general type, and let $T^*$ be the smallest invariant subtree.  Let $K$ be the kernel of the action of $G$ on $T^*$.  Then $G^{*+}/K$ is trivial or abstractly simple, and every nontrivial subgroup of $G/K$ normalized by $G^{*+}/K$ contains $G^{*+}/K$.\end{cor}

In a \tdlc group acting on a tree with property $\propP{}$, the structure of RIO subgroups is also quite restricted.  In particular, RIO subgroups of general type are open and inherit $\propP{}$.  Subgroups not of general type do occur, but we can regard them as having ``small rank'', at least relative to a vertex stabilizer, and so they are not so interesting as far as $\Sclass$-well-foundedness is concerned.

\begin{prop}\label{prop:tree_RIO}
Let $T$ be a tree and let $G$ be a closed subgroup of $\Aut(T)$ with property $\propP{}$, such that arc stabilizers in $G$ are compact.  Let $H$ be a RIO subgroup of $G$.
\begin{enumerate}[(i)]
\item If $H$ contains a translation, then for all arcs $a$ of $T$ belonging to an axis of translation of $H$, then $G_a \le \Res{}(H)$.  In particular, $H$ is open in $G$ and has property $\propP{}$.
\item If $H$ is not of general type, then $H$ has an open normal subgroup $K$, such that $H/K$ is isomorphic to $\triv$, $\Zb/2\Zb$ or $\Zb$, and such that every compactly generated subgroup of $K$ fixes a vertex.
\end{enumerate}
\end{prop}

\begin{proof}
Suppose $H$ contains a translation $g \in H$.  Then there is $K \in \Oc(H)$, where $K$ is an intersection of open subgroups of $G$, such that $g \in K$.  Let $a$ be an arc on the axis of $g$.  We claim that every open subgroup $O$ of $G$ normalized by $g$ also contains $G_a$, which would then imply $G_a \le \Res{}(K) \le \Res{}(H)$.  To prove this, let $\pi_a$ be the closest point projection of $VT$ onto $a$, let $T_a = \pi\inv_a(t(a))$, and decompose $G_a$ as $G_a = \rist_G(T_a) \times \rist_G(T_{\ol{a}})$.  By replacing $g$ with $g\inv$ if needed we may ensure $gT_a$ is contained in $T_a$.  Then $\bigcap_{n \ge 0}g^nT_a = \emptyset$; since $O$ is open, there are $v_1,\dots,v_n \in VT$ such that $\bigcap^n_{i=1}G_{v_i} \le O$.  Now take $n$ large enough that $g^nT_a$ is disjoint from $\{v_1,\dots,v_n\}$.  Then $g^n\rist_G(T_a)g^{-n}$ fixes each of $v_1,\dots,v_n$ and hence is contained in $O$, so also $\rist_G(T_a) \le O$.  The argument to show $\rist_G(T_{\ol{a}}) \le O$ is similar.  In particular, $G_a \le H$, so $H$ is open in $G$.

Now consider a different arc $e$ of $T$.  We have a similar decomposition $G_e = \rist_G(T_e) \times \rist_G(T_{\ol{e}})$, and the arc $a$ belongs to one of the half-trees, say $a \in T_{\ol{e}}$; hence
\[
\rist_G(T_e) \le G_a \le \Res{}(K) \le H.
\]
Thus
\[
H_e =  (\rist_G(T_e) \times \rist_G(T_{\ol{e}})) \cap H = \rist_G(T_e) \times (\rist_G(T_{\ol{e}}) \cap H),
\]
in other words, $H$ has property~$\propP{}$ with respect to edges of $T$, and hence has property~$\propP{}$ by Theorem~\ref{propertyP_oneclosure}.  This completes the proof of (i).

We may assume from now on that $H$ is not of general type.  Given the desired conclusion, we may assume that there is no subgroup $H'$ of $H$ fixing a vertex such that $|H:H'|\le 2$.  Thus $H$ does not fix a vertex or preserve an undirected edge.  Since $H$ is also not of general type, it follows by \cite[Proposition~3.4]{Tits70} that $H$ fixes at least one end $\delta$.  Let $\beta$ be a Busemann function for $\delta$ as in Lemma~\ref{lem:Busemann} and let $K = \ker\beta$.  By Lemma~\ref{lem:Busemann} and the fact that vertex stabilizers are open, we see that $H/K$ is trivial or cyclic and that every compactly generated subgroup of $K$ fixes a vertex.  This completes the proof of (ii).
\end{proof}

From Proposition~\ref{prop:tree_RIO}, we obtain a regionally near-simple class.  Let $\ms{T}_{ne}$ denote the class of \tdlc groups $G$ appearing as closed groups of automorphisms of a countable tree with property $\propP{}$, such that arc stabilizers are compact, local actions are elementary, and $G$ is not elementary.

\begin{cor}\label{cor:tree:near-simple}
The class $\ms{T}_{ne}$ is a regionally near-simple class.
\end{cor}

\begin{proof}
Let $G \in \ms{T}_{ne}$.  The fact that $G$ is a closed subgroup of the automorphism group of a countable tree ensures that $G$ is second-countable.  Since $G$ is not elementary, every sufficiently large $H \in \Oc(G)$ is nonelementary; moreover, given Proposition~\ref{prop:tree_RIO} and the fact that vertex stabilizers are elementary, every nonelementary $H \in \mc{IO}_c(G)$ has a general type action on the tree and inherits property $\propP{}$, from which we deduce that $H \in \ms{T}_{ne}$.  In particular, $G$ itself must have general type action on the tree.  We have proven property (c) of Definition~\ref{def:near-simple}; it remains to check properties (a) and (b).

For property (a): write $T^*$ for the smallest invariant tree for $G$, and $K$ for the kernel of the action of $G$ on $T^*$; then $K$ is compact, since $G$ has compact arc stabilizers.  By Corollary~\ref{cor:Tits}, the group $G^{*+}/K$ is topologically simple and is the monolith of $G/K$; since $G$ is not elementary, while $G^{*+}$ is an open normal subgroup of $G$, we see that $G^{*+}/K$ is not elementary, so certainly it is not discrete, nor is it regionally elliptic.  Thus $K = \RadRE(G)$ and $G^{*+}/K$ is the near-monolith of $G$.

For property (b): given an open subgroup $H$ of $G$ such that $G^{*+} \le HK$, then clearly $H$ is not elementary.  We deduce that $H \in \ms{T}_{ne}$ as required.
\end{proof}

More generally, given a class $\mc{S}$ with property $\propS$, we can divide the compactly generated open subgroups into a few cases as far as the ordering $\ll_{\mc{S}}$ is concerned.  First, we need a lemma to rule out a certain possibility for the (near-)monolith.

\begin{lem}\label{lem:tree_aniso}
Let $T$ be a tree and let $G$ be a geometrically dense closed locally compact subgroup of $\Aut(T)$ with property $\propP{}$.  Suppose also that $G^+ \in \Es{}(2)$.  Then $G^+ = \triv$, in other words, $G$ acts freely on the arcs of $T$.
\end{lem}

\begin{proof}
Note that by Theorem~\ref{thm:comp_gen+geom_dense}, the arc stabilizers of $G$ are compact.

Suppose that $G^+$ contains a translation and let $a$ be an arc of $T$.  Then the smallest invariant subtree of $G^+$, being $G$-invariant, must be equal to $T$.  Thus $a$ lies on the axis of some translation $g \in G^+$, pointing towards the attracting end.  Given $K \in \Oc(G^+)$ containing $g$, then by Proposition~\ref{prop:tree_RIO} we have $G_a \le \Res{}(K)$, but $\Res{}(K) = \triv$ since $G^+ \in \Es{}(2)$; thus $G_a = \triv$.  Since the arc $a$ was arbitrary, we conclude that $G^+ = \triv$, a contradiction.  Thus $G^+$ does not contain any translations.

Since $G^+$ is normal in $G$, we see that $G^+$ cannot preserve an undirected edge and invert this edge (since the edge would be unique and hence $G$-invariant).  Thus by \cite[Proposition~3.4]{Tits70}, $G^+$ fixes a vertex or end.  Moreover, if $G^+$ fixes the end $\xi$, then $G^+$ also fixes $g\xi \neq \xi$ where $g$ is some element of $G$ that does not fix $\xi$, and then $G^+$ fixes pointwise the line from $\xi$ to $g\xi$.  So in fact $G^+$ fixes a vertex.

Given that $G^+$ is normal in $G$, it now follows that $G^+$ fixes pointwise a set of vertices that is not contained in any proper subtree, which ensures that $G^+ = \triv$.
\end{proof}

\begin{prop}\label{prop:tree_ll}
Let $T$ be a countable tree and let $G$ be a closed subgroup of $\Aut(T)$ with property $\propP{}$, such that arc stabilizers in $G$ are compact.  Let $H \in \mc{IO}_c(G)$, and if $H$ is of general type, write $T^*$ for the smallest invariant subtree for $H$ and $K$ for the kernel of the action of $H$ on $T^*$.  Given a class $\mc{S}$ of \tdlcsc groups with property $\propS$, write $\rho_{\mc{S}}$ for the rank function associated to the partial order $\ll_{\mc{S}}$ on $\Oc(G)$.  Then exactly one of the following holds:
\begin{enumerate}[(i)]
\item (Bounded type) $H$ is not of general type, and given $H',H'' \in \Oc(H)$ such that $H' \ll_{\mc{S}} H''$, then $|H':H'_v| \le 2$ for some $v \in VT$.  In particular, if $G(v) \in \Es{\mc{S}}$ for all $v \in VT$, then
\[
H \in \Es{\mc{S}}, \; \rho_{\mc{S}}(H) \le {\sup_{v \in VT}}^+\xip{\mc{S}}(G(v)) \text{ and } \xip{\mc{S}}(H) \le {\sup_{v \in VT}}^+\xip{\mc{S}}(G(v)) + 1.
\]
\item (Free general type) $H$ is of general type and $H^{*+} = K$.  In this case, $H \in \Es{}(3)$, so $\Res{\omega}(H)=\triv$ and we have $\rho_\Sclass(H) \le 1$ and $\xip{\Sclass}(H) \le 2$.
\item (Almost $\Sclass$ type) $H$ and $H^{*+}$ are both of general type and $H^{*+}/K \in \Sclass$.  In this case, we have $\rho_\Sclass(H) \le 2$.
\
\item (Non-Noetherian type) $H$ and $H^{*+}$ are both of general type and $H^{*+} = \Res{\omega}(H)K$, but $H^{*+}/K$ is not compactly generated.  In this case, we have $L \ll_{\Sclass} H$ for all $L \in \Oc(H^{*+})$; moreover, $H$ does not have any normal factor in $\Sclass$.
\end{enumerate}
\end{prop}

\begin{proof}
Suppose $H$ is not of general type; we claim in this case that $H$ is of bounded type.  Let $H'' \in \Oc(H)$.  Then $H''$ is not of general type either, so Proposition~\ref{prop:tree_RIO}(ii) applies and there is an open normal subgroup $H^\sharp$ of $H''$ such that every compactly generated subgroup of $H^\sharp$ fixes a vertex, with $H''/H^\sharp \in \{\triv, \Zb/2\Zb,\Zb\}$.  In particular, note that $\Res{\omega}(H'') \le H^\sharp$ and that the largest possible compact subgroup of $H''/H^\sharp$ has order $2$.  Consider now $H' \in \Oc(H)$ such that $H' \ll_{\mc{S}} H''$.  The fact that $H'\Res{\omega}(H'')/\Res{\omega}(H'')$ is compact ensures that $|H'H^\sharp/H^\sharp| \le 2$, and hence $H^\sharp \cap H'$ has index at most $2$ in $H'$.  In particular, $H^\sharp \cap H'$ is compactly generated, so it fixes a vertex $v \in VT$, and hence $|H':H'_v| \le 2$.  Suppose now that $G(v) \in \Es{\mc{S}}$ for all $v \in VT$.  We note that for all $v \in VT$, the group $G_{v,1}$ is compact, since arc stabilizers are compact; thus $\xip{\mc{S}}(G_v) = \xip{\mc{S}}(G(v))$ by Corollary~\ref{cor:extensions}.  The conclusions about $\rho_G(H)$ and $\xip{\mc{S}}(H)$ now follow by the same argument as in Lemma~\ref{lem:step}(i).

From now on we may assume that $H$ is of general type; by Proposition~\ref{prop:tree_RIO}(i), it follows that $H$ is open and has property $\propP{}$.  By Lemma~\ref{lem:propP_subtree}, the action of $H/K$ on $T^*$ also has property $\propP{}$.

If $H^{*+}/K \in \Es{}(2)$, then $H^{*+}= K$ by Lemma~\ref{lem:tree_aniso}, so $H$ is compact-by-discrete; from here it is easy to see that $H$ is of free general type.  So from now on we may assume that $H^{*+}/K \not\in \Es{}(2)$; indeed by Lemma~\ref{lem:chief_rank}, we then have $H^{*+}/K \not\in \Es{}(\omega)$.  In particular, $H^{*+}/K$ is nondiscrete, and we see that $\Res{\omega}(H)K = H^{*+}$, that is, $\Res{\omega}(H)$ is cocompact in $H^{*+}$.

We see now that $H^{*+}$ acts nontrivially on $T^*$ and $T^*$ has more than two ends; it follows that the action of $H^{*+}$ is of general type, preserving no proper subtree of $T^*$ (see for instance \cite[Lemme~4.4]{Tits70}).  Thus $H$ and $H^{*+}$ are both of general type, and by Theorem~\ref{thm:Tits}, the group $H^{*+}/K$ is simple.  We see that $H^{*+}/K$ represents the only robust chief block of $H_0$, for any finite index open subgroup $H_0$ of $H$; given Lemma~\ref{lem:semisimple:basic}, it follows that $H_0$ has a normal factor in $\Sclass$ if and only if $H^{*+}/K \in \Sclass$.

If $H^{*+}/K$ is compactly generated, then $H^{*+}/K \in \Sclass$.  It follows in this case that given $H' \in \Oc(H)$ such that $H' \ll_{\Sclass} H$, then $\Res{\omega}(H') \le K$, so $\Res{\omega}(H')$ is compact and hence trivial.  From here it is clear that $\rho_\Sclass(H) \le 2$, and we conclude that $H$ is of almost $\Sclass$ type.

In the final case, $H^{*+}/K$ is not compactly generated, so $H^{*+}/K \not\in \Sclass$.  Given a finite index open subgroup $H_0$ of $H$, then $H_0$ has no normal factors in $\Sclass$.  In particular, we see that $L \ll_{\Sclass} H$ for all $L \in \Oc(H^{*+})$, and conclude that $H$ is of non-Noetherian type.
\end{proof}

We refer to the last case as ``non-Noetherian type'' since it is the only case that directly contradicts the Noetherian property.  In particular, notice that if the non-Noetherian type case is ruled out and the local actions are well-founded, then $G$ itself is $\Sclass$-well-founded.

If we have some (possibly infinite) bound on the elementary rank of local actions of $G$, we can bound the elementary rank of elementary members of $\mc{IO}_c(G)$.

\begin{cor}\label{cor:tree_elementary}
Let $T$ be a countable tree and let $G$ be a closed subgroup of $\Aut(T)$ with property $\propP{}$, such that arc stabilizers in $G$ are compact.  Suppose that the local actions of $G$ are all elementary and that $\alpha$ is an ordinal such that $\xi(G(v)) \le \alpha+1$ for all $v \in VT$.  Let $H \in \mc{IO}_c(G)$.  Then exactly one of the following holds:
\begin{enumerate}[(i)]
\item $H \in \Es{}(1+\alpha+2)$;
\item $H$ is nonelementary, of almost $\ms{S}$ type or non-Noetherian type, and belongs to the class $\ms{T}_{ne}$.
\end{enumerate}
\end{cor}

\begin{proof}
Retain the notation of Proposition~\ref{prop:tree_ll} and let $H \in \mc{IO}_c(G)$, considering the types from Proposition~\ref{prop:tree_ll}.  If $H$ is nonelementary, then $H \in \ms{T}_{ne}$ and from the proof of Corollary~\ref{cor:tree:near-simple}, we see that $H$ has general type action on $T$; in Proposition~\ref{prop:tree_ll} we see the possible types for $H$ are almost $\ms{S}$ type and non-Noetherian type.  Thus we may assume from now on that $H$ is elementary.

We see that every vertex stabilizer in $G$ belongs to $\Es{}(1+\alpha+1)$.  If $H$ is not of general type, then $H$ has an open normal subgroup $N$ such that every compactly generated subgroup of $N$ fixes a vertex; in particular, $\xi(N) \le 1+\alpha+1$, so $\xi(H) \le 1+\alpha+2$.  If $H$ is of free general type, we see that $H$ is compact-by-discrete, so $\xi(H) \le 3$.  Since $H$ is elementary, it cannot be of almost $\Sclass$ type.  The only remaining case to consider is that $H$ is of non-Noetherian type; in this case $H$ has general type action on $T$ and is elementary of infinite rank.

Suppose for a contradiction that $H$ is elementary and $\xi(H) > 1+\alpha + 2$.  Then there is $L \in \Oc(H)$ of rank exactly $1+\alpha+3$, by Lemma~\ref{lem:ele_rank_witnesses}.  In particular, we see that $L$ must be of non-Noetherian type, as all the other types have been ruled out.  Let $M$ be the kernel of the action of $L$ on its smallest invariant tree.  By Lemma~\ref{lem:chief_rank}, we have $\xi(L^{*+}/M) = \beta+1$, where $\beta$ is an infinite limit ordinal; by standard properties of the elementary rank, we then have $\xi(L^{*+}) = \xi(\Res{}(L)) = \beta+1$ and hence $\xi(L) = \beta+2$.  In particular, $\xi(L) \neq 1+\alpha+3$, a contradiction.
\end{proof}

Here are some sufficient conditions for $G$ to be $\Bigsclass$-well-founded, respectively $\Sclass$-well-founded.

\begin{prop}\label{prop:tree_Noetherian}
Let $T$ be a countable tree and let $G$ be a closed subgroup of $\Aut(T)$ with property $\propP{}$, such that arc stabilizers in $G$ are compact and $G(v) \in \Es{\Bigsclass}$ for all $v \in VT$.  Let $\alpha$ be an ordinal such that $\xi(E) \le \alpha+1$ for all $v \in VT$ and elementary subgroups $E$ of $G(v)$.
\begin{enumerate}[(i)]
\item We have $G \in \Es{\Bigsclass}$ and
\[
\xip{\Bigsclass}(G) \le \max\{{\sup_{v \in VT}}^+\xip{\Bigsclass}(G(v)) + 1,1+\alpha+2\}.
\]
\item Suppose that $G$ is Noetherian and $G(v) \in \Es{\Sclass}$ for all $v \in VT$.  Then $G \in \Es{\Sclass}$ and
\[
\xip{\Sclass}(G) \le \max\{{\sup_{v \in VT}}^+\xip{\Sclass}(G(v)) + 1,3\}.
\]
\end{enumerate}
\end{prop}

\begin{proof}
Retain the notation of Proposition~\ref{prop:tree_ll} and let $H \in \Oc(G)$.

Consider first the case that $H$ is elementary.  Then by Corollary~\ref{cor:tree_elementary} and Proposition~\ref{prop:elementary} we have $H \in \Es{}(1+\alpha+2) \subseteq \Es{\Bigsclass}(1+\alpha+2)$.

If $H$ is nonelementary and of non-Noetherian type, we are in a situation where $H$ has closed normal subgroups $K < H^*$, where $K$ and $H/H^*$ belong to $\Es{}(\omega)$ (indeed $K$ and $H/H^*$ are compact and discrete respectively) and $H^*/K$ is simple but not elementary (since otherwise $H$ would be elementary).  The same argument as Lemma~\ref{lem:simple_well-founded} shows that $\rho_{\Bigsclass}(H) \le 2$.

If $H$ is not of non-Noetherian type, then we see that
\[
\rho_{\mc{S}}(H) \le \max\{{\sup_{v \in VT}}^+\xip{\mc{S}}(G(v)),2\},
\]
where $\mc{S}$ is either $\Sclass$ or $\Bigsclass$ as appropriate.  The assertion (i) is now clear.

Now suppose that $G$ is Noetherian.  Then $H^{*+}$ must be compactly generated for all $H \in \Oc(G)$ of general type, so the non-Noetherian type case of Proposition~\ref{prop:tree_ll} does not occur.  The desired bound for $\xip{\Sclass}(G)$ as in (ii) then follows immediately from the previous paragraph.
\end{proof}

We now follow an approach inspired by an article of Caprace--Wesolek (\cite{CWes}), in which the first examples were found of nonelementary groups $G$ such that $H/\Res{}(H)$ is infinite for every compactly generated subgroup $H$ of $G$.

If we are in a situation where the local actions cannot be generated by point stabilizers, we can eliminate almost $\Sclass$ type from Proposition~\ref{prop:tree_ll}.  Rather than trying to exhaust all the possibilities, let us focus on the case where the local actions are all nilpotent; this is a convenient assumption that rules out certain types of permutation action and also clearly passes to subgroups.  In this context we find a strong dichotomy in the RIO subgroups.

\begin{lem}\label{lem:nilpotent_intransitive}
Let $G$ be a nontrivial nilpotent group acting on a set $X$.  Then the subgroup of $G$ generated by point stabilizers is intransitive.
\end{lem}

\begin{proof}
We may assume $G$ is transitive.  Choose $x \in X$; then the subgroup generated by point stabilizers is
\[
H = \grp{G_x \mid x \in X} = \grp{gG_xg\inv \mid g \in G}.
\]
Using the upper central series, we see that every proper subgroup of $G$ is contained in a proper normal subgroup of $G$.  In particular, since $G_x < G$, we have $G_x \le N$ where $N$ is a proper normal subgroup of $G$.  Since $N$ contains $G_x$ but is a proper subgroup, we see that the orbit $Nx$ is properly contained in $Gx$.  Now $H \le N$, so $Hx$ is properly contained in $Gx$, as required.
\end{proof}

\begin{thm}\label{thm:nilpotent_tree}
Let $T$ be a countable tree and let $G$ be a closed subgroup of $\Aut(T)$ with property $\propP{}$, such that arc stabilizers in $G$ are compact.  Suppose that $G(v)$ is nilpotent for all $v \in VT$, and let $H$ be a RIO subgroup of $G$.  Then $H$ has no normal factor in $\Sclass$, and exactly one of the following holds:
\begin{enumerate}[(i)]
\item $H$ is not of general type and belongs to $\Es{}(4)$.
\item $H$ is of free general type; in particular, $H \in \Es{}(3)$.
\item $H$ is of non-Noetherian type and $H$ is not $\Sclass$-well-founded.
\end{enumerate}
\end{thm}

\begin{proof}
Let us note why proving that all RIO subgroups of $G$ satisfy (i)--(iii) implies $H$ cannot have any normal factor in $\Sclass$.  If $H$ has a normal factor in $\Sclass$, then so does some $K \in \Oc(H)$; we then have $K \in \mc{IO}_c(G)$.  By Proposition~\ref{prop:tree_ll}, $K$ cannot be of non-Noetherian type.  But then $K$ also cannot be elementary, so cases (i) and (ii) are ruled out, a contradiction.

It now suffices to prove the cases when $H$ is a compactly generated IO subgroup of $G$, since each of (i)--(iii) is stable under directed unions of open subgroups.  So let us assume that $H$ is compactly generated.

By \cite{WillisNilpotent}, all nilpotent \tdlcsc groups belong to $\Es{}(2)$; thus by Corollary~\ref{cor:tree_elementary}, either $H \in \Es{}(4)$ or $H$ is not elementary, and in the latter case we see that $H$ must be of general type.  So if $H$ is not of general type, we see that case (i) holds; we may assume from now on that $H$ is of general type.  Let $T^*$ be the smallest invariant subtree for $H$ and let $K$ be the kernel of the action of $H$ on $T^*$.  Note that by Proposition~\ref{prop:tree_RIO}, this action has property $\propP{}$.  If $H$ is of free general type, we are in case (ii), so let us assume this is not the case; hence the action of $H^{*+}$ on $T^*$ is of general type.

The action of $H/K$ on $T^*$ has property $\propP{}$, leaves no proper subtree invariant and has nilpotent local actions $(H/K)(v)$.  By Lemma~\ref{lem:nilpotent_intransitive}, in all of the local actions, the subgroup $(H/K)(v)$ generated by point stabilizers is intransitive.  By Corollary~\ref{cor:Gplus:comp_gen} it follows that $H^{*+}/K$ cannot be compactly generated.  Thus $H$ is not of almost $\Sclass$ type.

By Proposition~\ref{prop:tree_ll}, the only remaining possibility is that $H$ is of non-Noetherian type.  In particular $H \not\in \Es{}(\omega)$; thus by Corollary~\ref{cor:tree_elementary}, we have $H \in \ms{T}_{ne}$.  Thus in all cases, either $H \in \Es{}(4)$ or $H \not\in \Es{}$, and moreover $H$ is nonelementary if and only if it is of non-Noetherian type.

It remains to show that if $H \not\in \Es{}$, then $H$ is not $\Sclass$-well-founded.  Set $H_0 = H$.  Suppose we have chosen $H_i \in \Oc(H)$ of non-Noetherian type.  We then take $H_{i+1} \in \Oc(H^{*+}_i)$; by Proposition~\ref{prop:tree_ll}, we have $H_{i+1} \ll_{\Sclass} H_i$.  We see that $H^{*+}_i$ is nonelementary; hence by taking $H_{i+1}$ large enough, we may ensure that $H_{i+1}$ is nonelementary, so it is itself of non-Noetherian type.  Thus we obtain an infinite descending chain in $(\Oc(H),\ll_{\Sclass})$, showing that $H$ is not $\Sclass$-well-founded.
\end{proof}

\begin{cor}\label{cor:nilpotent_tree}
Let $T$ be a countable tree and let $G$ be a compactly generated closed geometrically dense subgroup of $\Aut(T)$ with property $\propP{}$, such that arc stabilizers in $G$ are compact.  Suppose that $G(v)$ is nilpotent for all $v \in VT$ and that for some $v \in VT$, $G(v)$ does not act freely on $o\inv(v)$.  Then $G \in \Es{\Bigsclass}(3)$, but $G \not\in \Es{\Sclass}$.  Moreover, there is no RIO subgroup of $G$ with a quotient in $\Sclass$.
\end{cor}

For examples of nonelementary groups satisfying the hypotheses of Corollary~\ref{cor:nilpotent_tree}, we can turn to the well-studied family of groups first introduced by M. Burger and Sh. Mozes in \cite{BurgerMozes}, as stated in Theorem~\ref{intro:notEsp}.  Given $d \ge 3$ finite, and a subgroup $F$ of $\Sym(d)$, the \defbold{Burger--Mozes group} $\Univ(F)$ is a subgroup of $\Aut(T_d)$ characterized up to conjugacy in $\Aut(T_d)$ by the following properties: $\Univ(F)$ has property $\propP{}$; the local action at every vertex is isomorphic to $F$; and for every edge of $T_d$, there is an involution in $\Univ(F)$ inverting that edge.  In particular, we see that if $F$ is nilpotent and does not act freely, then Corollary~\ref{cor:nilpotent_tree} applies to $G = \Univ(F)$.

In contrast to Theorem~\ref{thm:nilpotent_tree}, Theorem~\ref{intro:tree} will provide us with examples of groups $G$ acting on trees with property $\propP{}$, such that $G$ belongs to both $\Sclass$ and $\Es{\Sclass}$.  Theorem~\ref{intro:tree} is based on a sufficient condition, given a \tdlc group $G$ acting on a tree, for every proper open subgroup of $G$ to fix a vertex, which we prove in the appendix (Theorem~\ref{thm:open_primitive}).

\begin{lem}\label{lem:primitive_elementary}
Let $X$ be a countable set and let $G$ be a \tdlc group that acts faithfully and primitively on $X$ with compact open stabilizers.  Suppose that $G$ is not discrete.  Then $G/\Res{\omega}(G)$ is compact and $G$ is nonelementary.
\end{lem}

\begin{proof}
Let $x \in X$ and write $G_x$ for the point stabilizer.  Then $G_x$ is both compact and a maximal subgroup of $G$.  In particular, we see that $G = \grp{g,G_x}$ for any $g \in G \smallsetminus G_x$, so $G$ is compactly generated.

Suppose that $G/\Res{\omega}(G)$ is not compact.  Then we cannot have $\Res{\omega}(G)G_x = G$, so we must have $\Res{\omega}(G) \le G_x$.  Moreover, the intersection of conjugates of $G_x$ in $G$ is trivial, so in fact $\Res{\omega}(G) = \triv$.  Thus $G \in \Es{}$ and $\xi(G) < \omega$.  Since $G$ is compactly generated, it follows that $G/\Res{}(G)$ is not compact, so the same argument shows $\Res{}(G) = \triv$.  Thus $G$ is a SIN group by Lemma~\ref{lem:elementary2}; in particular, $G$ has a compact open normal subgroup $N$.  We then see that $G \neq NG_x$, so $N \le G_x$, and hence $N = \triv$, so $G$ is discrete.

Now consider the case that $G/\Res{\omega}(G)$ is compact.  If $G$ is elementary, we see that $G$ must be compact; in this case, $G_x$ has only finitely many conjugates in $G$, from which we conclude that $G$ is finite, so in particular $G$ is discrete.  We have now proved via the contrapositive that if $G$ is nondiscrete, then $G/\Res{\omega}(G)$ is compact and $G$ is nonelementary.
\end{proof}

\begin{proof}[Proof of Theorem~\ref{intro:tree}]
Since $G$ has transitive local action at every vertex, we see that $G$ has at most two orbits on arcs of $T$ and the action on $T$ is geometrically dense.  The arc stabilizers are thus compact by Theorem~\ref{thm:comp_gen+geom_dense}; in particular, the local action at each vertex has compact point stabilizers.  In particular, for each vertex $v$, then $G_{v,1}$ is a compact normal subgroup of $G_v$.  Thus $G_v \in \Es{\Sclass}$ with $\xip{\Sclass}(G_v) = \xip{\Sclass}(G(v))$.  Moreover, by Lemma~\ref{lem:primitive_elementary}, we see that $G_v \in \Es{}$ if and only if $G(v)$ is discrete.  Note also that $G(v)$, and hence also $G_v$, is compactly generated for all $v \in VT$.

Note that $(G^{+})^{+} = G^{+}$.  Given $v \in VT$, since the action of $G_v$ on $o\inv(v)$ is primitive but not regular, it is generated by point stabilizers; thus $G_v \le G^+$.  All the conclusions of Theorem~\ref{thm:open_primitive} now apply to $G^+$; in particular, $G^+$ acts transitively on each part of the natural bipartition of $VT$, so $|G:G^+| \le 2$.  Thus from now on we may assume $G = G^+$.  Then $G$ is a group in $\Sclass$ of the form $G = G_x \ast_{G_{(x,y)}} G_y$ for adjacent vertices $x$ and $y$, such that every proper open subgroup of $G$ is contained in a conjugate of $G_x$ or $G_y$.  Write 
\[
\max\{\xip{\Sclass}(G(x)),\xip{\Sclass}(G(y))\} = \alpha+1;
\]
then we have $\xip{\Sclass}(H) \le \alpha+1$ for every proper open subgroup $H$ of $G$, so $\rho_G(H) \le \alpha$ for every proper open subgroup of $G$.  As for $G$ itself, we have $G \in \Sclass$, so $\rho_G(G) \le 2$ by Lemma~\ref{lem:simple_rho}.
Thus for all $H \in \Oc(G)$ we have $\rho_G(H) \le \max\{\alpha,2\}$, and hence $G \in \Es{\Sclass}$ with $\xip{\Sclass}(G) \le \max\{\alpha+1,3\}$.

On the other hand, it is clear that $\xip{\Sclass}(G) \ge \xip{\Sclass}(G_v)$ for all $v \in VT$, so $\xip{\Sclass}(G) \ge \alpha+1$.  From the fact that $G$ only has two orbits on vertices, we see that $G$ contains a translation $s$; in particular, $\grp{s}$ is an infinite discrete cyclic subgroup of $G$.  Thus $G$ is not regionally elliptic, so $\xip{\Sclass}(G) \ge 2$.

The rank is now determined exactly except when $\alpha \le 1$, so let us assume that $\alpha \le 1$; the only question for $\xip{\Sclass}(G)$ is whether we have $\rho_G(G)=1$ or $\rho_G(G)=2$.  As in Lemma~\ref{lem:simple_rho}, let $\mc{H}$ be the set of $H \in \Oc(G)$ such that $H$ is elementary and $\xi(H) < \omega$; we are interested in possible noncompact elements of $\mc{H}$.  Let $H \in \mc{H}$ and suppose that $H$ is not compact.  Then see that $H$ is contained in one of $G_x$ and $G_y$; say $H \le G_x$.  If $H \ll_{\Sclass} G_x$, then $HG_{x,1}/G_{x,1} \ll_{\Sclass} G(x)$; since $\xip{\Sclass}(G(x)) \le 2$ it follows that $HG_{x,1}/G_{x,1}$ is compact, so $H$ is compact, a contradiction.  Thus $H \not\ll_{\Sclass} G_x$.  Since $H \in \mc{H}$, we know that $\Res{\omega}(H) = \triv$; thus the only way to have $H \not\ll_{\Sclass} G_x$ is for $H\Res{\omega}(G_x)/\Res{\omega}(G_x)$ to be noncompact, which then implies $G(x)/\Res{\omega}(G(x))$ is noncompact.  By Lemma~\ref{lem:primitive_elementary} we deduce that $G(x)$ is an infinite discrete group.  Conversely, if $G(x)$ is infinite and discrete, then $G_x$ is a noncompact element of $\mc{H}$.  Lemma~\ref{lem:simple_rho} now determines $\rho_G(G)$, and hence $\xip{\Sclass}(G)$, in all cases.
\end{proof}

Theorem~\ref{intro:tree} includes many known examples of groups in $\Sclass$.  As a basic example, for $d \ge 3$, the automorphism group $\Aut(T_d)$ of the regular tree $T$ of degree $d$ has primitive local action $\Sym(d)$, so we conclude that $\Aut(T_d) \in \Es{\Sclass}(2)$; the subgroup $\Aut(T_d)^+$ belongs to $\Sclass$.  Moving away from locally finite trees, we have a continuum of examples of groups $S \in \Sclass$ constructed by Smith in \cite{SmithDuke}; in these examples, there are two local actions, and we can take one the actions to be $\Sym(3)$ and the other to be a countably infinite discrete group acting primitively with nontrivial finite stabilizers.  In this case, Theorem~\ref{intro:tree} shows that $\xip{\Sclass}(S)=3$.  Thus there are $2^{\aleph_0}$ groups in $\Sclass \cap \Es{\Sclass}$ of rank $3$.

\begin{rmk}\label{rmk:bad_rank}
The most basic example of a nonelementary group in $\Es{\Sclass}$, namely $G = \Aut(T_3)$, demonstrates that the rank function $\xip{\Sclass}$ does not behave well on closed subgroups, even closed cocompact subgroups.  Specifically, we have $G \in \Es{\Sclass}$ with $\xip{\Sclass}(G)=2$; indeed, by Theorem~\ref{thm:open_primitive}, every noncompact open subgroup contains $G^+$, so $|\Reg(G)|=2$.  On the other hand, taking $F$ to be the subgroup $\grp{(12)}$ of $\Sym(3)$, then $\Univ(F)$ is a cocompact subgroup of $G$, but by Theorem~\ref{intro:notEsp}, the poset $(\Reg(\Univ(F)),\ll_{\Sclass})$ is not even well-founded.  More broadly, it is unlikely that there is any way of assigning a well-founded ordinal rank to the groups $\Aut(T_d)$, in a manner that assigns small rank to discrete groups and is stable under (cocompact) closed subgroups as well as directed unions, extensions and quotients, for the following reason.  Given an arbitrary \tdlcsc group $A$ (which can presumably have ``large rank''), take $H \in \Oc(A)$ of ``large rank'', and then take a Cayley--Abels graph $\Gamma$ for $H$, of degree $d$ say; the kernel $K$ of the action of $H$ on $\Gamma$ is then compact.  The universal cover of $\Gamma$ is then a regular tree $T$, and the action of $H$ lifts to a group $\widetilde{H}$ acting on $T$, equipped with a quotient map $\theta: \widetilde{H} \rightarrow H/K$.  The kernel of $\theta$ is discrete.  Thus we expect the cocompact closed subgroups of $\Aut(T_d)$ (allowing $d$ to range over the natural numbers) to witness all possible ranks of \tdlcsc groups.
\end{rmk}

\subsection{Locally linear groups}\label{sec:locally_linear}

We recall a structure theorem of Caprace and T. Stulemeijer, on locally linear \tdlc groups.  (Some details that are not relevant for the present purposes are omitted.)  Let $k$ be a \tdlc local field.  We say a \tdlc group $G$ is \defbold{locally linear (over $k$)} if there is a compact open subgroup $U$ of $G$ admitting a continuous faithful finite-dimensional representation over a local field (over $k$).  For Caprace--Stulemeijer, a \defbold{topologically simple algebraic group over $k$} is a topologically simple \tdlc group that is topologically isomorphic to $H(k)^+/\Z(H(k)^+)$, where $H$ is a $k$-simple algebraic $k$-group and $H(k)^+$ is the normal subgroup generated by $k$-rational points of split unipotent $k$-subgroups of $H$.

\begin{thm}[{\cite[Theorem~1.1]{CS}}]\label{thm:CS_linear}
Let $G$ be a locally linear \tdlc group.  Then $G$ has a series of closed normal subgroups 
\[
\triv \le R \le G_1 \le G_0 \le G
\]
enjoying the following properties.

The group $R$ is a closed characteristic subgroup and is locally soluble.  The group $G_0$ is an open characteristic subgroup of finite index in $G$.  The quotient group $H_0 = G_0/R$, if nontrivial, has finitely many nontrivial closed normal subgroups, say $M_1, \dots, M_m$, satisfying the following properties.
\begin{enumerate}[(i)]
\item For some $l \le m$ and all $i \le l$, the group $M_i$ is a topologically simple algebraic group over a local field $k_i$.  In particular $M_i$ is compactly generated and abstractly simple.
\item For all $j > l$, the group $M_j$ is compact.
\item The group $H_1 = G_1/R$ coincides with the product $M_1 \dots M_m \cong M_1 \times \dots \times M_m$, which is closed in $H_0$.  Moreover $H_0/H_1 = G_0/G_1$ is locally abelian.
\end{enumerate}
\end{thm}

From this theorem and a couple of other facts from the literature, we can easily deduce Theorem~\ref{intro:linear}.

\begin{proof}[Proof of Theorem~\ref{intro:linear}]
Given Lemma~\ref{lem:rank_inequalities}(iii), it is enough to prove the result in the case that $G$ is compactly generated; in particular, we may assume that $G$ is second-countable.  Decompose $G$ as in Theorem~\ref{thm:CS_linear}.  By Lemma~\ref{lem:rank_inequalities}(ii), replacing $G$ with $G_0$ does not change the rank, so we may assume $G = G_0$.   The groups $R$ and $G/G_1$ are locally soluble; by \cite[Theorem~8.1]{WesEle}, every locally soluble \tdlcsc group $A$ belongs to $\Es{}(\omega)$, and hence $A \in \Es{\Sclass}(2)$.  The groups $M_{l+1}, M_{l+2},\dots, M_m$ are obviously $\Sclass$-well-founded of rank $1$, since they are compact.  This leaves the topologically simple algebraic groups $M_1,\dots,M_l$.  Given such a group $M_i$, then by \cite[Theorem~(T)]{Prasad82}, every proper open subgroup of $M_i$ is compact, so $\xip{\Sclass}(M_i) \le 2$.  Given the direct product decomposition of $H_1$, we then have $\xip{\Sclass}(H_1) \le 2$.

Putting the factors together and applying Theorem~\ref{thm:extensions}, we see that $G \in \Es{\Sclass}$, with
\[
\xip{\Sclass}(G) \le (\xip{\Sclass}(R) - 1) + (\xip{\Sclass}(G_1/R) - 1) + \xip{\Sclass}(G/G_1)  \le 1 + 1 + 2 = 4. \qedhere
\]
\end{proof}

\begin{rmk}
Every closed subgroup of a locally linear \tdlcsc group is locally linear.  In particular, the locally linear \tdlcsc groups form an interesting class of $\Sclass$-well-founded groups that are not all elementary, but still have the property that every closed subgroup is $\Sclass$-well-founded.  Given a linear \tdlcsc group $G$ such that no nontrivial element of $G$ has open centralizer, then there is a first-countable \tdlc group $L = \ms{L}(G)$, which is the universal \tdlc group locally isomorphic to $G$ with no nontrivial discrete normal subgroup (see \cite{BEW} and \cite{CapDeM}).  The group $L$ is then locally linear, so $L \in \Es{\Sclass}(4)$.
\end{rmk}

\subsection{Kac--Moody groups}\label{sec:KM}

Kac--Moody groups over finite fields are another important source of examples of groups in $\Sclass$.  Significantly for the present context, a theorem of Caprace and T. Marquis (\cite[Theorem~A]{CapMar}) gives a classification up to finite index of the open subgroups of a complete geometric Kac--Moody group $G$ over a finite field.  In particular, the structure of $\Reg(G)$ is known and we can deduce that $G$ is $\Sclass$-well-founded of finite rank.

We will not recall the construction of complete geometric Kac--Moody groups here, as we do not need to know anything about Kac--Moody groups \textit{per se} except for results from \cite{CapMar}.  (See \cite{Marquis} for a detailed account of these groups.)  However, we will need to use the more basic concept of $BN$-pairs (see \cite{AB} for reference).

A \defbold{$BN$-pair} in a group $G$ is a pair $(B,N)$ of subgroups such that $G = \grp{B,N}$ and $T:= B \cap N$ is normal in $N$, and $W = N/T$ admits a generating set $S$ such that the following holds:
\begin{enumerate}[(i)]
\item Given $s \in S$ and $w \in W$, then $BwB . BsB \subseteq BwB \cup BwsB$;
\item For all $s \in S$ we have $sBs\inv \not\subseteq B$.
\end{enumerate}
It turns out that if such a set $S$ exists, it is uniquely determined by $B$ and $N$, and $(W,S)$ is a Coxeter group, which is the \defbold{Weyl group} of the $BN$-pair.  Given $J \subseteq S$, write $W_J = \grp{J}$.  We call the subgroups $W_J$ for $J \subseteq S$ the \defbold{standard parabolic subgroups} of $(W,S)$, and a \defbold{parabolic subgroup} of $W$ is a subgroup of the form $wW_Jw\inv$ for some $w \in W$.  We note that, although the Weyl group is not a subgroup of $G$, it is still true that sets of the form $wP$ and $Pw$ are cosets of $P$ for $w \in W$ and any $P \ge B$, namely $wP = w'P$ and $Pw = Pw'$ for any $w' \in N$ such that $w'T = w$.

We also recall that $G$ admits a double coset decomposition, the \defbold{Bruhat decomposition} $G = \bigsqcup_{w \in W} BwB$.  Analogously to subgroups of the Weyl group, subgroups of $G$ of the form $P_J = BW_JB$ are called \defbold{standard parabolic subgroups} of $G$ (or of the $BN$-pair); these are precisely the subgroups of $G$ that contain $B$.  A \defbold{parabolic subgroup} is then a conjugate of a standard parabolic subgroup.

One of the ways to form a locally compact Kac--Moody group, which we call a \defbold{complete geometric Kac--Moody group} here, is to let a minimal Kac--Moody group $\mc{G}$ over a finite field act on its positive building $\Delta_+$, and then take the closure $G$ of the image of $\mc{G}$ in $\Aut(\Delta_+)$, where $\Aut(\Delta_+)$ carries the permutation topology on chambers.  Writing $B$ and $N$ for the stabilizers in $G$ of the fundamental chamber and fundamental apartment respectively, one finds that $(B,N)$ is a $BN$-pair for $G$ such that $B$ is a compact open subgroup and the generating set $S$ of the Weyl group is finite.  In particular, with the given $BN$-pair structure, every parabolic subgroup of $G$ is compactly generated and open.  The result of Caprace--Marquis is essentially the converse to this fact, up to finite index.

\begin{thm}[{\cite[Theorem~A and Corollary~B]{CapMar}}]\label{thm:capmar}
Let $G$ be a complete geometric Kac--Moody group over a finite field and let $O$ be an open subgroup of $G$.  Then the set of parabolic subgroups of $G$ containing $O$ as a finite index subgroup is finite but nonempty.  In particular, $O$ is compactly generated, so $G$ is Noetherian.
\end{thm}

For the present purposes we need to also take account of which parabolic subgroups are virtually contained in one another.  The following is an expanded version of an argument outlined to me by P.-E. Caprace in correspondence.

Given a Coxeter group $(W,S)$, let $\Lambda_{(W,S)}$ be the set of standard parabolic subgroups of $W$ ordered by inclusion, and let $\Lambda^f_{(W,S)}$ be the same poset taking subgroups up to finite index, that is, the quotient poset $\Lambda_{(W,S)}/\sim_f$.

\begin{thm}\label{thm:building_commensurate}
Let $G$ be a group with a $BN$-pair with Weyl group $(W,S)$, such that $B$ is a commensurated subgroup of $G$ and $|S| < \infty$; given $J \subseteq S$, let $P_J$ be the standard parabolic subgroup $BW_JB$ of $G$.  Suppose $J,J' \subseteq S$ and $g \in G$ are such that $gP_{J}g\inv$ is virtually contained in $P_{J'}$.  Then $[W_{J}] \le [W_{J'}]$ as elements of $\Lambda^f_{(W,S)}$; we have $[W_J] = [W_{J'}]$ if and only if $gP_Jg\inv$ is commensurate with $P_{J'}$.

Conversely, if $[W_J] \le [W_{J'}]$ then $P_J$ is virtually contained in $P_{J'}$.
\end{thm}

We begin the proof with some lemmas.  For convenience we take ``Coxeter group'' to include the assumption that the Coxeter generating set is finite, and similarly for the Weyl group of a $BN$-pair.

Given a Coxeter group $(W,S)$ and $J \subseteq S$, write $J^{\sph}$ for the union of the spherical components of $J$ (in terms of its Coxeter diagram) and $J^{\infty} = J \smallsetminus J^{\sph}$.  In particular, $(W_{J^{\sph}},J^{\sph})$ is finite.  Say $J$ is \defbold{essential} if $J = J^{\infty}$.

Standard parabolics of a Coxeter group can be conjugate, but the conjugating element is necessarily of a special form given by V. Deodhar (\cite{Deo}); in particular, the essential part is preserved.

\begin{lem}\label{lem:Deo}
Let $(W,S)$ be a Coxeter group and let $w \in W$ and $J \subseteq S$ such that $wJw\inv \subseteq S$.  Then $w \in W_{S^{\sph}}$.  In particular, $J^{\infty} = wJ^{\infty}w\inv = (wJw\inv)^{\infty}$.
\end{lem}

\begin{proof}
By \cite[Proposition~5.5]{Deo}, one has a decomposition of $w$ as $w = v_m v_{m-1} \dots v_1$, where for each $v_i$, there is some $s_i \in S$ such that $v_i(S \smallsetminus \{s_i\})v\inv_i \subseteq S$.  Moreover, in the discussion before \cite[Proposition~5.5]{Deo}, it is shown that any such element $v_i$ belongs to a spherical component of $(W,S)$.  Thus $w \in W_{S^{\sph}}$.  In particular, since clearly $J^{\infty} \subseteq S^{\infty}$, we see that $w$ commutes with $W_{J^{\infty}}$, so $J^{\infty}$ is fixed under conjugation by $w$.
\end{proof}

We deduce a stability property of the partial order on $\Lambda^f_{(W,S)}$ with respect to conjugation in $W$.

\begin{lem}\label{lem:parabolic_inclusion}
Let $(W,S)$ be a Coxeter group.  Then the following are equivalent, for $J,J' \subseteq S$:
\begin{enumerate}[(i)]
\item $W_J$ is virtually contained in $W_{J'}$;
\item $J^{\infty} \subseteq (J')^{\infty}$;
\item There is $w \in W$ such that $wJ'w\inv \cap J$ generates a subgroup of finite index in $W_J$.
\end{enumerate}
\end{lem}

\begin{proof}
Suppose (i) holds.  We see that $W_{J} = W_{J^{\sph}} \times W_{J^{\infty}}$ and $W_{J^{\sph}}$ is finite; since $J^{\infty}$ is essential, $W_{J^{\infty}}$ has no proper parabolic subgroups of finite index.  So we must have $W_{J^{\infty}} \le W_{J'}$ and hence $J^{\infty} \subseteq J'$.  It is then clear that no $s \in J^{\infty}$ belongs to a spherical component of $J'$; hence $J^{\infty} \subseteq (J')^{\infty}$.  Thus (i) implies (ii).

If (ii) holds, then clearly (iii) holds with $w=1$.

Now suppose (iii) holds; let $J'' = wJ'w\inv \cap J$.  Then $W_{J''}$ has finite index in $W_J$, so $(J'')^{\infty} = J^{\infty}$, and $w\inv (J'')^{\infty}w \subseteq J' \subseteq S$, so by Lemma~\ref{lem:Deo}, $w\inv (J'')^{\infty}w = (J'')^{\infty}$.  Thus $J^{\infty} \subseteq J'$, so $W_{J'}$ contains the finite index subgroup $W_{J^{\infty}}$ of $W_J$, showing that (i) holds.
\end{proof}

We now obtain some conditions under which parabolic subgroups of a group with a $BN$-pair are commensurate with each other.

\begin{lem}\label{lem:parallel}
Let $G$ be a group with a $BN$-pair with Weyl group $(W,S)$, such that $B$ is a commensurated subgroup of $G$; let $\Delta$ be the associated building.
\begin{enumerate}[(i)]
\item Given $J,J' \subseteq S$, then $P_J$ is virtually contained in $P_{J'}$ if and only if $[W_J] \le [W_{J'}]$.
\item Let $R$ and $R'$ be a pair of parallel residues in $\Delta$.  Then $\mathrm{Stab}_G(R)$ is commensurate with $\mathrm{Stab}_G(R')$.
\end{enumerate}
\end{lem}

\begin{proof}
(i)
We can write $P_J$ as a product $P_J = P_{J^{\sph}}P_{J^{\infty}}$ and similarly for $P_{J'}$.  Since $B$ is commensurated, $P_{J^{\sph}} = BW_{J^{\sph}}B$ is a union of finitely many left cosets of $B \le P_{J^{\infty}}$, so $J^{\sph}$ does not contribute to the commensurability class of $P_J$.  Thus we may assume $J$ and $J'$ are essential.  We then see by the Bruhat decomposition that if $P_J$ is virtually contained in $P_{J'}$, then $[W_J] \le [W_{J'}]$; conversely if $[W_J] \le [W_{J'}]$, the lack of spherical components means that $W_J \le W_{J'}$ and hence $P_J \le P_{J'}$.

(ii)
Let $P = \mathrm{Stab}_G(R)$ and $P' = \mathrm{Stab}_G(R')$, let $J$ be the type of $R$ and let $J'$ be the type of $R'$.  By \cite[Proposition~21.8]{MPW}, the fact that $R$ and $R'$ are parallel means that the Weyl distance $w$ from $R$ to $R'$ satisfies $J = wJ'w\inv$.  In particular, by Lemma~\ref{lem:Deo} we have $w \in W_{S^{\sph}}$ and $J^{\infty} = (J')^{\infty}$.  Take residues $R^{\infty}$ and $(R')^{\infty}$ of type $J^{\infty}$ contained in $R$ and $R'$ respectively.  Then by (i), $\mathrm{Stab}_G(R^{\infty})$ has finite index in $P$ and $\mathrm{Stab}_G((R')^{\infty})$ has finite index in $P'$; moreover, the Weyl distance from $R^{\infty}$ to $(R')^{\infty}$ is still $w$ (see \cite[Proposition~21.10]{MPW}), so $R^{\infty}$ and $(R')^{\infty}$ are parallel.  Thus we may assume $J = J^{\infty} = J'$, that is, that $R$ and $R'$ are of the same type.

Since $w \in W_{S^{\sph}}$, we see that $R$ and $R'$ are both contained in a residue $R''$ of type $J \times S^{\sph}$.  As a building, we can write $R'' = X \times Y$ where $X$ is a building of type $J$ and $Y$ is a building of type $S^{\sph}$, so that $R = X \times \{y\}$ and $R' = X \times \{y'\}$ for some chambers $y,y' \in Y$.  Since $B$ is commensurated, we see that $\Delta$ is locally finite, so for any $x \in X$ the $P$-orbit of $(x,y')$ is finite, and hence $P$ is virtually contained in $P'$; similarly $P'$ is virtually contained in $P$.
\end{proof}

The next lemma is essentially \cite[Lemma~4.2]{CapMin}.

\begin{lem}\label{lem:essential}
Let $(W,S)$ be a Coxeter group with no spherical components.  Then there is $w \in W$ that is \defbold{$S$-essential}, that is, $w$ is not contained in any proper parabolic subgroup of $W$.  Moreover, if $w$ is $S$-essential then so is $w^n$ for all $n > 0$.
\end{lem}

\begin{proof}
We can write $W = \prod^n_{i=1}W_{S_i}$ where $(W_{S_i},S_i)$ is irreducible and nonspherical, hence $W_{S_i}$ is infinite.  By \cite[Lemma~4.2]{CapMin}, there is $w_i \in W_{S_i}$ such that all positive powers of $w_i$ are $S_i$-essential.  Finally, we observe that any product of $S_i$-essential elements of $W_{S_i}$, taking one element for each $i$, will form an $S$-essential element of $(W,S)$, due to the fact that every parabolic subgroup of $(W,S)$ decomposes according to its intersections with the irreducible components of $(W,S)$.  In particular, any positive power of $w = w_1w_2 \dots w_n$ is essential.
\end{proof}

With these lemmas in hand, we can finish the proof of the theorem.

\begin{proof}[Proof of Theorem~\ref{thm:building_commensurate}]
Suppose $J,J' \subseteq S$ and $g \in G$ are such that $gP_{J}g\inv$ is virtually contained in $P_{J'}$; we aim to show $[W_{J}] \le [W_{J'}]$.  Given Lemma~\ref{lem:parabolic_inclusion} and Lemma~\ref{lem:parallel}(i), we are free to assume $J$ and $J'$ are essential.  By Lemma~\ref{lem:essential} we can then take an element $v \in W_J$, all of whose positive powers are $J$-essential.

Let $\Delta$ be the building associated to the $BN$-pair.  Then $P_J$ and $P_{J'}$ are the stabilizers of the standard $J$-residues $R_J$ and $R_{J'}$ respectively in $\Delta$.  Let $A$ be an apartment containing a chamber of $gR_J$ and a chamber of $R_{J'}$.  Then we can regard $A$ as a Coxeter complex for the copy of $W$ induced by $\mathrm{Stab}_G(A)$; in particular, there is some element $v'$ of $\mathrm{Stab}_G(A) \cap gP_Jg\inv$ that acts as $v$ on $A$.  By replacing $v'$ by a positive power we may assume that $v' \in P_{J'}$.  By \cite[Theorem~4.97]{AB}, $R_1 = \mathrm{proj}_{gR_{J}}(R_{J'})$ is a subresidue of $gR_J$ of type $J \cap wJ'w\inv$, where $w$ is the Weyl distance from the simplex defining $gR_J$ to the simplex defining $R_{J'}$.  Now we see that $v'$ stabilizes $R_1$.  Since $v'$ is $J$-essential on $A$, we must have $J \cap wJ'w\inv = J$ and hence $R_1 = gR_J$.  In particular, we have $[W_{J}] \le [W_{J'}]$ by Lemma~\ref{lem:parabolic_inclusion}.

Suppose now that $[W_J] = [W_{J'}]$.  Then by Lemma~\ref{lem:parabolic_inclusion}, in fact $J = J'$, and we must also have $J = wJw\inv$.  It follows that the residues $gR_J$ and $R_{J'}$ are parallel, so by Lemma~\ref{lem:parallel}(ii), $gP_Jg\inv$ is commensurate with $P_{J'}$.  Conversely, if we suppose that $gP_Jg\inv$ is commensurate with $P_{J'}$, then the argument of the previous paragraph shows that $J \cap wJ'w\inv = J$ and $J' \cap w\inv Jw = J'$, so by Lemma~\ref{lem:parabolic_inclusion} we have $[W_J] = [W_{J'}]$.

On the other hand, given $J,J' \subseteq S$ such that $[W_J] \le [W_{J'}]$, then $P_J$ is virtually contained in $P_{J'}$ by Lemma~\ref{lem:parallel}(i).
\end{proof}

In particular, Theorem~\ref{thm:building_commensurate} applies to the parabolic subgroups of complete geometric Kac--Moody groups, and hence via Theorem~\ref{thm:capmar} we obtain a complete classification of which open subgroups can be virtually conjugated inside one another.

\begin{cor}\label{cor:capmar:exact}
Let $G$ be a complete geometric Kac--Moody group over a finite field, equipped with its defining $BN$-pair structure.  Then there is a unique surjective map $\theta: \mc{O}(G) \rightarrow \Lambda^f_{(W,S)}$, with the following properties:
\begin{enumerate}[(i)]
\item For all $J \subseteq S$, then $\theta(P_J) = [W_J]$;
\item Given $H,K \in \mc{O}(G)$, then $\theta(H) = \theta(K)$ if and only if there is $g \in G$ such that $gHg\inv$ is commensurate with $K$;
\item Given $H,K \in \mc{O}(G)$, then $\theta(H) < \theta(K)$ if and only if there is $g \in G$ such that $gHg\inv$ is commensurate with a subgroup of $K$ of infinite index.
\end{enumerate}
\end{cor}

With this level of control over commensurability classes of open subgroups, we immediately see that $G \in \Es{\Sclass}$.

\begin{cor}\label{cor:KM:Esp}
Let $G$ be a complete geometric Kac--Moody group over a finite field.  Then $(\Reg(G),\le)$ is well-founded and has rank equal to that of the finite poset $\Lambda^f_{(W,S)}$, where $(W,S)$ is the Weyl group of $G$.  In particular, $G$ is $\Sclass$-well-founded and $\xip{\Sclass}(G) \le \rho(\Lambda^f_{(W,S)}) < \omega$.
\end{cor}

\subsection{Groups with just infinite locally normal subgroups}\label{sec:hji}

A profinite group $U$ is \defbold{just infinite} if every nontrivial closed normal subgroup of $U$ is open, and \defbold{hereditarily just infinite} (\hji) if every open subgroup is just infinite.  A \tdlc group is \defbold{locally \hji} if it has a compact open \hji subgroup.

A \tdlc group $G$ has \defbold{property (LD)} if, given a closed locally normal subgroup $K$ of $G$, there is an open subgroup $K_0$ of $K$ that is a direct factor of an open subgroup of $G$.  The class of groups with property (LD) is introduced in \cite{Reid-LD}; in particular, it includes all \tdlc groups locally isomorphic to a just infinite profinite group.  We summarize some structural results about these groups.  (See also \cite[Proposition~5.1 and Theorem~5.4]{BEW}, which gave similar results for the structure of locally \hji \tdlc groups.)

\begin{thm}[{\cite[Theorem~1.7]{Reid-LD}}]\label{thm:ld_char}Let $G$ be a \tdlc group.  The following are equivalent:
\begin{enumerate}[(i)]
\item $G$ has (LD);
\item $G$ locally isomorphic to a profinite group of the form
$$ \prod_{i \in I}L_i,$$
such that finitely many factors $L_i$ (possibly none) are just infinite profinite groups and the remaining factors are finite simple groups.
\end{enumerate}
\end{thm}

\begin{thm}\label{thm:ld}Let $G$ be a \tdlc group with property~(LD).
\begin{enumerate}[(i)]
\item (\cite[Lemma~1.10]{Reid-LD}) There is a characteristic closed subgroup $Q$ of $G$ such that $G/Q$ has property (LD), $\QZ(G/Q) = \triv$, and $\ol{\QZ(Q)}$ is open in $Q$.
\item (See \cite[Theorems~1.9, 1.11, 4.8 and 4.11]{Reid-LD}) Suppose $\QZ(G) = \triv$.  Then $G$ is first-countable and $\Res{}(G)$ is a direct factor of an open subgroup of $G$.  In addition, $\Res{}(G)$ is a direct product of finitely many groups $P_1,\dots,P_n$, each of which is locally isomorphic to a just infinite profinite (\hji or branch) group and is either residually discrete or topologically simple.
\end{enumerate}
\end{thm}

From Theorem~\ref{thm:ld} we can easily deduce the following.

\begin{cor}\label{cor:ld_wellfounded}
Let $G$ be a first-countable \tdlc group with property~(LD).  Then $G \in \Es{\Bigsclass}(4)$; if $\QZ(G) = \triv$, then $G \in \Es{\Bigsclass}(3)$.  The same conclusions apply with $\Es{\Sclass}$ in place of $\Es{\Bigsclass}$ if we additionally assume $G$ is Noetherian.
\end{cor}

\begin{proof}
If $\QZ(G) = \triv$, then by Theorem~\ref{thm:ld}(ii) we see that the situation of Lemma~\ref{lem:simple_well-founded} applies, so $G \in \Es{\Bigsclass}(3)$.  In the general case, we see that the characteristic closed subgroup $Q$ as in Theorem~\ref{thm:ld}(i) belongs to $\Es{}(3)$, hence also to $\Es{\Bigsclass}(2)$.  Theorem~\ref{thm:extensions} now yields $G \in \Es{\Bigsclass}(4)$.

If $G$ is Noetherian, then all the nonelementary topologically simple groups that occur in the argument are compactly generated, since they appear as quotients of locally normal subgroups.  The same arguments therefore yield $G \in \Es{\Sclass}(3)$ if $\QZ(G) = \triv$ and $G \in \Es{\Sclass}(4)$ in general.
\end{proof}

In the absence of the Noetherian property, the question of whether \tdlc groups with property (LD) are $\Sclass$-well-founded is more complicated.  Given Theorem~\ref{thm:ld}(i), we can reduce to the case that $\QZ(G) = \triv$.  In that case, given $H \in \Oc(G)$, we know from Theorems~\ref{thm:ld_char} and~\ref{thm:ld}(ii) that $\Res{\omega}(H)$ is a direct product of finitely many topologically simple groups, each of which also has (LD).  We can split up the possible topologically simple groups as follows.

\begin{lem}\label{lem:ld:simple}
Let $G$ be a first-countable topologically simple \tdlc group with property (LD).  Then $G \in \Es{}(3) \sqcup \ms{R}$.
\end{lem}

\begin{proof}
Note that the classes $\Es{}(3)$ and $\ms{R}$ are indeed disjoint.  If $\QZ(G) > \triv$, then $\QZ(G)$ is dense in $G$ by topological simplicity, so $G \in \Es{}(2)$ by Lemma~\ref{lem:quasi-discrete}.  We may therefore assume from now on that $\QZ(G)=\triv$, ensuring that $G$ is locally isomorphic to a profinite group that is either \hji or a just infinite branch group.

Suppose $G$ is not regionally elementary.  Then for a sufficiently large $H \in \Oc(G)$, then $H$ is not elementary, so $\Res{\omega}(H)$ is nontrivial.  If $G$ is locally \hji, it is then clear that $\Res{\omega}(H)$ is open; set $H' = H$.  If $G$ is locally isomorphic to a just infinite profinite branch group, we can enlarge $H$ to some $H' \in \Oc(G)$ with $\Res{\omega}(H')$ open, as follows.  The group $\Res{\omega}(H)$ represents an element $\alpha$ of the local decomposition lattice $\mc{L}$ of $G$.  Since $G$ is locally isomorphic to a just infinite profinite group, $G$ acts with clopen orbits on the Stone space of $\mc{L}$; since in addition, $G$ has no nontrivial fixed points in $\mc{L}$ (see \cite[Theorem~4.11]{Reid-LD}), in fact $G$ acts transitively on the Stone space of $\mc{L}$.  There is therefore $g_1,\dots,g_n \in G$ such that $\bigvee^n_{i=1}g_i\alpha = \infty$ in $\mc{L}$.  It then follows that the group $\grp{g_i\Res{\omega}(H)g\inv_i \mid 1 \le i \le n}$ is open in $G$, and hence $\Res{\omega}(H')$ is open in $G$ where $H' = \grp{H,g_1,\dots,g_n}$.  In particular, we see that if $U$ is a compact open subgroup of $\Res{\omega}(H')$, then $U$ has trivial core in $H'$, showing that $H'$ is expansive.  Thus $G \in \ms{R}$.

We may now suppose that $G \in \Es{}$.  Suppose $\xi(G) > 3$.  Then by Lemma~\ref{lem:ele_rank_witnesses}, $G$ has a compactly generated open subgroup $H$ with $\xi(H) = 4$; note that $H$ is second-countable.  By Theorem~\ref{thm:ld}(ii) $\Res{}(H)$ is a direct product of factors $P_1,\dots,P_n$, each of which is either residually discrete group or topologically simple group in $\Es{}(4)$; in fact, by Lemma~\ref{lem:chief_rank}, the groups $P_1,\dots,P_n$ are all in $\Es{}(2)$.  This implies $H \in \Es{}(3)$, a contradiction.  Thus $\xi(G) \le 3$.
\end{proof}

Most likely, regionally elementary topologically simple groups $G$ with property (LD) have $\xi(G) = 2$; certainly this is true by Lemma~\ref{lem:chief_rank} if $G$ is second-countable.  Combining Lemma~\ref{lem:ld:simple} with Theorem~\ref{thm:ld}(ii) gives the following.

\begin{cor}\label{cor:ld:simple}
Let $G$ be a \tdlc group with property (LD), such that $\QZ(G) = \triv$.  Then $\Res{\omega}(G)$ is a direct product of finitely many (possibly none) robustly monolithic topologically simple \tdlc groups.
\end{cor}

There are many examples of topologically simple groups that are locally isomorphic to profinite branch groups, including compactly generated examples.  See for example the groups $\Univ(F)$, or the groups $\ms{N}_F$ of \cite{Lederle}; in the case that  $F$ is a finite $2$-transitive permutation group with perfect point stabilizers, both $\Univ(F)$ and $\ms{N}_F$ are locally isomorphic to a just infinite profinite branch group.  Given a profinite branch group $U$ that is commensurable with one of its proper direct powers, one can also follow the approach of \cite{Roe} to produce a group in $\Sclass$ locally isomorphic to $U$.  The situation with locally \hji \tdlc groups is less well-understood.  There are certainly topologically simple examples, however, the only examples I am aware of that are either directed unions of compact open subgroups, or simple algebraic groups over local fields; in the latter case, by \cite[Theorem~T]{Prasad82}, all proper open subgroups are compact, so the group is both $\Sclass$-well-founded and Noetherian.  I do not know if there exists a \tdlc group with property (LD) that is not $\Sclass$-well-founded.  We can at least establish a relationship between $\Sclass$-well-foundedness and two other properties of $G$.

\begin{prop}\label{prop:hji:open_in_S}
Let $G$ be a first-countable group with property (LD), such that $\QZ(G) = \triv$.  Suppose that $G$ is $\Sclass$-well-founded but not regionally elementary.  Then there is an open subgroup of $G$ that is a finite direct product of groups in $\Sclass$.
\end{prop}

\begin{proof}
Since $G$ is not regionally elementary, by Theorem~\ref{thm:ld}(ii) there is a topologically simple locally normal subgroup $K$ that is a direct factor of an open subgroup of $G$; hence by Lemma~\ref{lem:ld:simple} we have $K \in \ms{R}$.  By Corollary~\ref{cor:R_descent} applied to the regionally near-simple class $\ms{R}$, there is $H \in \Oc(K)$ such that $\Mon(H/\RadRE(H)) \in \Sclass$.  In fact, from Theorem~\ref{thm:ld}(ii) we see that $H$ has a characteristic closed subgroup $M$ that is a direct product of finitely many groups in $\Sclass$.  Now $H$ is locally normal in $G$, so $M$ is also locally normal in $G$; as in the proof of Lemma~\ref{lem:ld:simple}, there are $g_1,\dots,g_n \in G$ such that the group $N = \grp{g_iMg\inv_i \mid 1 \le i \le n}$ is open in $G$.  From the way $N$ is generated it is clear that $\Res{\omega}(N) = N$, so by Corollary~\ref{cor:ld:simple}, $N$ is a finite direct product of robustly monolithic topologically simple groups $N_1,\dots,N_k$.  At the same time, $N$ is compactly generated; hence each of its topologically simple direct factors $N_i$ belongs to $\Sclass$.
\end{proof}

\subsection{Other classes of \tdlc group to investigate}\label{sec:other}

Here are four families of \tdlcsc groups in which some structural information is available from the existing literature and it could be interesting to determine which groups in the family are $\Sclass$-well-founded or at least $\Bigsclass$-well-founded.  For the first three, nonelementary examples are known; for the last, no such examples are known, but proving $\Bigsclass$-well-foundedness could be a first step to showing these groups are elementary.

\begin{enumerate}

\item Let $T$ be a locally finite tree of minimum degree at least $6$, and let $G$ be a closed subgroup of $\Aut(T)$ be such that $G$ acts $2$-transitively on the boundary and has local action at each vertex containing the alternating group.  All such groups have been classified by N. Radu (\cite{Radu}), and all have a finite index open subgroup in $\Sclass$.  These groups are not exactly within the scope of Theorem~\ref{thm:open_primitive}, but since the local action has a high degree of transitivity, it seems likely that strong restrictions on their open subgroups apply.  More generally, it seems likely that open subgroups of boundary-$2$-transitive subgroups of $\Aut(T)$ have a special structure, and even more so for the boundary-$3$-transitive subgroups of $\Aut(T)$ when $T$ is a regular locally finite tree.  (Note that in this case, $\Aut(T)$ itself acts $3$-transitively but not $4$-transitively on the boundary.)

\item The Burger--Mozes construction has been generalized in several directions; one direction is to actions on right-angled buildings that are not necessarily trees.  The most general construction of this kind to date is by J. Bossaert and T. De Medts (\cite{BosDeM}), building on work of Caprace and of De Medts, A. Silva and K. Struyve (\cite{CapBuild}, \cite{DSS}).  With sufficiently strong assumptions about the local actions, is conceivable that in some situations one could use the geometry of the building to control open subgroups with respect to stabilizers of residues, in a similar manner to Theorem~\ref{thm:capmar}, and from there, obtain more examples of nonelementary $\Sclass$-well-founded groups.  On the other hand, one could hope to find more non-$\Sclass$-well-founded groups by generalizing the situation of Theorem~\ref{thm:nilpotent_tree}.

\item Let $X$ be the Cantor set and let $G \le \Homeo(X)$.  The \defbold{piecewise full group} $\mathrm{F}(G)$ of $G$ consists of all $h \in \Homeo(X)$ for which there exists a clopen partition $X = X_1 \sqcup \dots \sqcup X_n$ of $X$ and $g_1,\dots,g_n \in G$, such that $hx = g_ix$ for all $x \in X_i$.  We say $G$ is \defbold{piecewise full} if $G = \mathrm{F}(G)$.  Piecewise full groups have been a productive source of examples of simple groups; indeed, by a theorem of Nekrasheyvych (\cite[Theorem~1.1]{Nekra}), every piecewise full group with minimal action has a simple normal subgroup.  This is also true for nondiscrete locally compact groups.  In particular, if we start with a group of automorphism of a locally finite tree and then consider the piecewise full group of its action on the boundary of the tree, we have examples of piecewise full groups $G$ with derived group $\Der(G) \in \Sclass$, such as Neretin's groups $N_d = \mathrm{F}(\Aut(T_d))$ of almost automorphisms of $T_d$ (\cite{Neretin}) and the coloured Neretin groups constructed in \cite{Lederle}.  More general properties and examples of locally compact piecewise full groups will be given in work under preparation by A. Garrido, D. Robertson and the author.

As noted in the introduction, some of the coloured Neretin groups certainly fall outside of $\Es{\Sclass}$, because they are piecewise full groups of non-$\Sclass$-well-founded groups acting on trees.  What is not clear is whether there is a non-$\Sclass$-well-founded group arising as the piecewise full group of an $\Sclass$-well-founded group.  Specifically, it would be interesting to know whether the groups $N_d$ are $\Sclass$-well-founded.

\item Given a topological group $G$ and $g \in G$, the \defbold{contraction group} $\con(g)$ of $g$ consists of all elements $x \in G$ such that $g^nxg^{-n} \rightarrow 1$ as $n \rightarrow +\infty$, and $G$ is \defbold{anisotropic} if $\con(g)=\triv$ for all $g \in G$.  Write $[\mathrm{Aniso}]$ for the class of anisotropic \tdlcsc groups.

Anisotropic groups are the degenerate case of the theory of dynamics of automorphisms of \tdlc groups, known as ``scale theory'', that has been developed principally by G. Willis and his collaborators, starting with \cite{Willis94}; see for instance the preliminaries section of the article \cite{ReidNY} and its references for examples of what is known about contraction groups and related topics in \tdlc group theory.  The inability to use this theory in a meaningful way has left the structure of groups in $[\mathrm{Aniso}]$ quite mysterious up to the present.

Within the class of \tdlcsc groups, $[\mathrm{Aniso}]$ is closed under taking closed subgroups, quotients and directed unions of open subgroups, but not extensions; one also has $\Es{}_{\aleph_0}(2) \subseteq [\mathrm{Aniso}]$.  (Closure of $[\mathrm{Aniso}]$ under quotients is a special case of \cite[Theorem~3.8]{BaumgartnerWillis}; the other closure properties are trivial to observe.)  Examples of elementary rank $3$ can be constructed following an approach outlined in \cite{KepertWillis}; however, as far as I am aware, all known examples of anisotropic first-countable \tdlc groups belong to $\Es{}(\omega)$.  At the same time it is hard to show that global properties of a group exclude it from the class $[\mathrm{Aniso}]$; for example it is still not known whether the intersection $\Sclass \cap [\mathrm{Aniso}]$ is empty.

On the other hand, groups in $[\mathrm{Aniso}]$ do have the following property, which is of interest from the perspective of $\Bigsclass$-well-foundedness: given $G \in [\mathrm{Aniso}]$ and $H \in \Oc(G)$, then $H$ has a finite normal series in which every factor belongs to $\Es{}(\omega+1)$ or is a quasi-product of copies of a topologically simple group (see \cite[Proposition~16]{ReidMATRIX}).  Further results along these lines might lead to a proof that $[\mathrm{Aniso}] \subseteq \Es{\Bigsclass}$.  On the other hand if $[\mathrm{Aniso}] \not\subseteq \Es{\Bigsclass}$, then there exists a \tdlcsc group $G$ that is aniostropic, not $\Bigsclass$-well-founded, but also topologically simple, which is a striking combination of properties.
\end{enumerate}

\appendix

\section{Some properties of groups acting on trees with Tits' independence property}

\begin{defn}
Given $G \le \Aut(T)$, the \defbold{$\propP{1}$-closure} of $G$, denoted by $G^{\propP{1}}$, is the set of automorphisms $g \in \Aut(T)$ such that for all $v \in VT$, and every finite set of neighbours $X$ of $v$, there exists $g_X \in G$ such that $gv = g_Xv$ and $gw = g_Xw$ for every vertex $w \in X$.  We say $G$ is \defbold{$\propP{1}$-closed} if $G = G^{\propP{1}}$.
\end{defn}

For closed subgroups of $\Aut(T)$, the $\propP{1}$-closure is the same as the $\propP{}$-closure, and in fact it is enough to consider property~$\propP{}$ with respect to edges.

\begin{thm}[See \cite{BanksElderWillis} Theorem 5.4 and Corollary 6.4]\label{propertyP_oneclosure}Let $T$ be a tree and let $G$ be a closed subgroup of $\Aut(T)$.  Then $G = G^{\propP{1}}$ if and only if $G$ has property $\propP{}$.  Furthermore, if $G$ has property $\propP{}$ with respect to the edges of $T$, then $G = G^{\propP{1}}$, so $G$ has property $\propP{}$ with respect to all simple paths.\end{thm}

\begin{proof}
Let $L$ be a nonempty simple path in $T$, and let $g \in \Aut(T)_{(L)}$ such that for each $v \in VL$, there exists $s_v \in G_{(L)}$ such that $g$ acts as $s_v$ on $\pi\inv_L(v)$.  We now claim that $g \in G^{\propP{1}}$ (indeed, $g \in (G_{(L)})^{\propP{1}}$).  Let $X$ be a finite set of vertices, all adjacent to some vertex $w$ of $T$ and let $X' = X \smallsetminus VL$.  We observe that since $X'$ is disjoint from $VL$ and the set $X' \cup \{w\}$ spans a subtree, any path from $X'$ to $L$ must pass through $x$, in other words $X' \subseteq \pi^{-1}_L(x)$.  In particular, $g$ agrees with $s_x$ on $X'$.  Since both $g$ and $s_x$ fix $VL$ pointwise, in fact $gv = s_xv$ for all $v \in X$.  Given the freedom of choice of $X$, we conclude that $g \in G^{\propP{1}}$ as claimed.  Thus if $G = G^{\propP{1}}$, then $G$ has property $\propP{}$.

Conversely, suppose that $G$ is closed and satisfies property $\propP{}$ with respect to the edges of $T$.  Suppose that $G \neq G^{\propP{1}}$ and let $g \in G^{\propP{1}} \smallsetminus G$.  Since $G$ is closed, the set $G^{\propP{1}} \smallsetminus G$ is a neighbourhood of $g$ in $G^{\propP{1}}$, so there is a finite set $X$ of vertices such that $g(G^{\propP{1}})_{(X)} \cap G = \emptyset$.  Let $S$ be the smallest subtree of $T$ containing $X$; note that $\Aut(T)_{(X)} = \Aut(T)_{(S)}$, since every vertex of $S$ lies on the shortest path between a pair of vertices in $X$.  Let us suppose that $X$ has been chosen so that $|S|$ is minimized.

By the definition of $G^{\propP{1}}$, we see that $S$ is not a star, so for every $x \in S$, there is a vertex in $S$ at distance $2$ from $x$. Hence there exist adjacent vertices $x$ and $y$ of $S$ such that neither $x$ nor $y$ is a leaf of $S$.  Let $L$ be the path formed by  the single arc $(x,y)$.  By the minimality of $|S|$, there is some $h \in G$ such that $gx = hx$ and $gy = hy$, so that $h^{-1}g$ fixes $L$ pointwise.  Let
\[
S_1 = (S \cap \pi^{-1}_L(x)) \cup \{y\} \text{ and } S_2 = (S \cap \pi^{-1}_L(y)) \cup \{x\}.
\]
Note that for $i=1,2$, then $S_i$ is the set of vertices of a subtree of $S$ that contains $L$.  The condition that neither $x$ nor $y$ is a leaf of $S$ ensures that there is some neighbour of $x$ in $S$ that is not contained in $S_2$, and similarly there is some neighbour of $y$ in $S$ that is not contained in $S_1$.  Hence $S_1$ and $S_2$ are both proper subtrees of $S$, so by the minimality of $|S|$, there exist $h_1,h_2 \in G$ such that $h_iw_i = h^{-1}gw_i$ for all $w_i \in S_i$ ($i=1,2$).  Indeed, $h_1$ and $h_2$ are elements of $G_{(L)}$, since $h_1$ and $h_2$ both agree with $h^{-1}g$ on $L$.  By (the restricted) property $\propP{}$, we can take $h_1$ to fix $S_2$ pointwise and $h_2$ to fix $S_1$ pointwise.  But then $h_1h_2$ agrees with $h\inv g$ on $X$, that is,
\[
hh_1h_2 \in g(G^{\propP{1}})_{(X)},
\]
which contradicts the hypothesis that $g(G^{\propP{1}})_{(X)} \cap G = \emptyset$.
\end{proof}

Under mild assumptions, property~$\propP{}$ passes to actions on subtrees.

\begin{lem}\label{lem:propP_subtree}
Let $G$ be a group acting on a tree $T$ with property~$\propP{}$, such that $G$ has compact arc stabilizers, and let $T'$ be a $G$-invariant subtree of $T$.  Then the action of $G$ on $T'$ has property~$\propP{}$.
\end{lem}

\begin{proof}
Let $K$ be the kernel of the action of $G$ on $T'$.  We first argue that $G$ acts on $T'$ as a closed subgroup of $\Aut(T')$. The action of $G$ on $T'$ is continuous because $G$ is Hausdorff and vertex stabilizers $G_v$ for $v \in VT'$ are open in $G$. Given $a \in AT'$, the continuous image of the compact group $G_a$ is compact in $G/K$, and therefore closed in the Hausdorff group $\Aut(T')$. In other words, the stabilizer in $G/K$ of $a$ is a closed subgroup of $\Aut(T')$, and it follows then that $G/K$ is a closed subgroup of $\Aut(T')$.

Now let $L$ be a  simple path in $T'$.  Let $\pi_L$ and $\pi'_L$ be the closest point projections for $L$ as a simple path in $T$ and $T'$ respectively; note that $\pi_L$ and $\pi'_L$ agree on $VT'$.  Let $\phi_x$ be the action homomorphism from $G_{(L)}$ to $\Aut(\pi\inv_L(x))$.  Since $G$ has property~$\propP{}$ on $T$, the natural homomorphism
\[
\phi_L: G_{(L)} \rightarrow \prod_{x \in VL} \phi_x(G_{(L)}); \; g \mapsto (\phi_x(g))_{x \in VL}
\]
is surjective.  Now consider what happens if we replace $\phi_L$ and $\phi_x$ with $\phi'_L$ and $\phi'_x$ respectively, which are now defined with respect to $T'$.  If we choose $g_x \in \phi'_x(G_{(L)})$ for each $x \in VL$, then there is some $h_x \in G_{(L)}$ such that $\phi'_x(h_x) = g_x$, and then by the surjectivity of $\phi_L$, there is $g \in G_{(L)}$ such that $\phi_L(g) = (\phi_x(h_x))_{x \in VL}$.  But then since
\[
(\pi')\inv_L(x) = \pi\inv_L(x) \cap VT' \subseteq \pi\inv_L(x),
\]
we immediately see that
\[
\phi'_L(g) = (\phi'_x(h_x))_{x \in VL} = (g_x)_{x \in VL}.
\]
Thus $\phi'_L$ is surjective, so $G$ has property~$\propP{}$ on $T'$. 
\end{proof}

The properties of being locally compact and compactly generated impose certain restrictions on the orbits and stabilizers of arcs of the tree.

\begin{thm}\label{thm:comp_gen+geom_dense}
Let $T$ be a tree and let $G$ be a closed subgroup of $\Aut(T)$ with property $\propP{}$.  Suppose that $G$ is locally compact and that the action of $G$ on $T$ is geometrically dense.
\begin{enumerate}[(i)]
\item Every arc stabilizer of $G$ is compact.
\item If $G$ is compactly generated, then $G$ has finitely many orbits on $VT \sqcup AT$ and every vertex stabilizer of $G$ is compactly generated.
\end{enumerate}
\end{thm}

\begin{proof}
Since $G$ is locally compact, there is some finite set $B$ of vertices of $T$, such that the pointwise stabilizer $H$ of $B$ in $G$ is compact.  In particular, $H(v)$ is compact Hausdorff, hence closed, for every $v \in VT$; since $H$ is open in $G$, it follows that $G$ has closed local actions.  Let $a \in AT$.  Then $a$ belongs to the axis of some translation $h \in G$; we can choose $h$ so that $a$ is oriented towards the repelling end of $h$.  Write $T_b = \pi\inv_b(t(b))$ for $b \in AT$.  Then the half-trees $T_{h^na}$ form an increasing family whose union is $T$; thus there is some $n$ such that $T_{h^na}$ contains $B$ and $t(h^na) \not\in B$.

Set $g = h^n$ and consider the stabilizer of $ga$; we have
\[
G_{ga} = \rist_G(T_{ga}) \times \rist_G(T_{\ol{ga}}).
\]
Our choice of $g$ ensures that $B \subseteq T_{ga}$; thus $\rist_G(T_{\ol{ga}}) \le H$, so $\rist_G(T_{\ol{ga}})$ is compact.  In particular, $G_{ga}$ has finite orbits on $T_{\ol{ga}}$.  After conjugating by $g\inv$ we see that $G_a$ has finite orbits on the half-tree $T_{\ol{a}}$.  A similar argument using $\ol{a}$ in place of $a$ shows that $G_a = G_{\ol{a}}$ also has finite orbits on the complementary half-tree $T_{a}$.  Thus at every vertex $v \in VT$, the point stabilizers of the closed local action $G(v)$ have finite orbits, that is, they are compact.  We then see that for each arc $b \in AT$, the group $G_b$ is an inverse limit of compact groups, so it is compact, proving (i).

Now suppose $G$ is compactly generated.  We can take a symmetric generating set for $G$ of the form $S = F \cup G_v$, where $F$ is finite.  Then the set $\{sv \mid s \in S\}$ is finite.  Let $T'$ be the subtree spanned by the paths $[v,sv]$ as $s$ ranges over $S$.  Then $T'$ is finite and for each $s \in S$, the graph $T' \cup sT'$ is connected: specifically, both $T'$ and $sT'$ are connected and contain $sv$.  From here, we see that the graph $\bigcup_{g \in G} gT'$ is also connected, and hence equal to $T'$.  This shows that $G$ has finitely many orbits on $VT \sqcup AT$.  The rest of part (ii) now follows by \cite[Proposition~4.1]{Castellano}.
\end{proof}

Our next aim is to obtain a restriction on the local action of the subgroup $G^+$ generated by arc stabilizers, in terms of the local action of $G$.  Our approach is inspired by \cite[Proposition~2.8 and Theorem~3.1]{CWes}.

\begin{defn}
Let $G$ be a group acting on a tree $T$ and let $v \in VT$.  A \defbold{$(G,v)$-focused colouring} of $T$ is a function $c: AT \rightarrow o\inv(v)$ with the following properties:
\begin{enumerate}[(a)]
\item $c$ restricts to the identity on $o\inv(v)$.
\item We have $c(a) = c(\ol{a})$ for all $a \in AT$.
\item Given $w \in Gv$, there exists $g_w \in G$ such that $g_ww = v$ and $c(a) = g_wa$ for all $a \in o\inv(w)$.
\item Given $w \in VT \smallsetminus Gv$, then $c$ is constant on $o\inv(w)$.
\end{enumerate}
\end{defn}

\begin{lem}\label{lem:focused_colouring}
Let $G$ be a group acting on a tree $T$ and let $v \in VT$.  Suppose that $G_v$ acts transitively on $o\inv(v)$.  Then there is a $(G,v)$-focused colouring of $T$.
\end{lem}

\begin{proof}
We construct the colouring inductively for the arcs in the ball $B_n(v)$ of radius $n$ around $v$.  To start with, we set $c(a) = a$ for all $a \in o\inv(v)$, and then $c(a) = \ol{a}$ for all $a \in t\inv(v)$; this defines the colouring on $B_1(v)$.  Suppose we have coloured $B_n(v)$ for $n \ge 1$ and let $w \in VT$ be such that $d(v,w) = n$.  There is then a unique arc $a \in o\inv(w)$ such that $d(v,t(a)) = n-1$; by the inductive hypothesis, we have already coloured $a$ and $c(a) = c(\ol{a})$.  If $w \not\in Gv$, we set 
\[
c(b) = c(\ol{b}) = c(a) \text{ for all } b \in o\inv(w).
\]
If instead $w \in Gv$, choose $g_w \in G$ such that $g_wa = c(a)$; this is possible because the local action at $v$ is transitive.  We then set
\[
c(b) = c(\ol{b}) = g_wb \text{ for all } b \in o\inv(w).
\]
Continuing in this way, we extend the colouring to all of $AT$; it is then clear that $c$ has the required properties.
\end{proof}

\begin{prop}\label{prop:G+_vertex}
Let $G$ be a group acting on a tree $T$ and let $v \in VT$ be such that $G_v$ acts transitively on $o\inv(v)$.  Then
\[
(G^+)_v = \grp{G_a \mid a \in o\inv(v)}.
\]
\end{prop}

\begin{proof}
Let $G(v)$ be the local action of $G$ at $v$ and let $G(v)^+$ be the subgroup of $G(v)$ generated by its point stabilizers.  Then given $a \in o\inv(v)$, then $G(v)^+$ is exactly the set of elements of $G(v)$ that send $a$ to an element of $G(v)^+a$; in other words, we have a free action of $G(v)/G(v)^+$ on the set of $G(v)^+$-orbits.

Construct a $(G,v)$-focused colouring $c$ as in Lemma~\ref{lem:focused_colouring} and let $g \in G_a$ for some arc $a \in AT$.  We claim that for all $b \in AT$, $c(gb)$ belongs to the same $G(v)^+$-orbit as $c(b)$.  We proceed by induction on 
\[
n = \min\{d(o(a),o(b)),d(o(a),t(b))\}.
\]
In the case $o(a) = o(b)$, there are two possibilities.  If $w \in Gv$, then the colouring of $o\inv(o(a))$ is given by $c(b) = hb$ for some $h \in G$ sending $w$ to $v$.  Since $g$ fixes $a$, we see that $hgh\inv$ fixes $ha \in o\inv(v)$, and hence $hgh\inv$ acts as an element of $G(v)^+$.  We then have $c(gb) = hgb = (hgh\inv)hb$, so $c(b)$ and $c(gb)$ lie in the same $G(v)^+$-orbit.  If instead $w \not\in Gv$, then $c(b) = c(gb)$ for every $b \in o\inv(o(a))$.  In either case, given $b \in t\inv(o(a))$, we have $c(b) = c(\ol{b})$ in the same $G(v)^+$-orbit as $c(gb) = c(g\ol{b})$.

Now consider an arc $b$ such that $d(o(a),o(b)) = n$ and $d(o(a),t(b)) = n+1$, for $n \ge 1$.  There is a unique arc $b'$ with $o(b') = o(b)$, such that $d(o(a),t(b')) = n-1$.  By the inductive hypothesis we know that $c(gb') \in G(v)^+ c(b')$.  If $o(b) \not\in Gv$, then $c(b) = c(b')$ and also $c(gb) = c(gb')$, so $c(gb) \in G(v)^+ c(b)$.  If instead $o(b) \in Gv$, we take $h,h' \in G$ such that $c(e) = he$ and $c(e') = h'e'$ for all $e \in o\inv(o(b))$ and $e' \in o\inv(o(gb))$ respectively.  Thus, for $e \in o\inv(o(b))$, we have
\[
 h'gh\inv c(e) = (h'gh\inv) he = h'ge = c(ge).
\]
Since $c(b')$ and $c(gb')$ are in the same $G(v)^+$-orbit, and since $G(v)/G(v)^+$ acts freely on the $G(v)^+$-orbits, it follows that $h'gh\inv$ acts on $o\inv(v)$ as an element of $G(v)^+$.  But then $c(e)$ and $c(ge)$ are in the same $G(v)^+$-orbit for all $e \in o\inv(o(b))$, so $c(gb) \in G(v)^+ c(b)$; hence also $c(g\ol{b}) \in G(v)^+ c(\ol{b})$. This completes the inductive step, so $c(gb) \in G(v)^+ c(b)$ for all $b \in AT$.

Since $G^+$ is generated by arc stabilizers, we conclude that for all $g \in G^+$ and $b \in AT$, then $c(gb)$ is in the same $G(v)^+$-orbit as $c(b)$.  In particular, $(G^+)_v$ stabilizes each $G(v)^+$-orbit, or to put it another way, every $(G^+)_v$-orbit on $o\inv(v)$ is contained in a $G^*$-orbit, where
\[
G^* =  \grp{G_a \mid a \in o\inv(v)}.
\]
On the other hand, it is clear that $G^* \le (G^+)_v$.  Since $G^*$ contains the point stabilizers of the action of $G_v$ on $o\inv(v)$, we conclude that $G^* = (G^+)_v$ as required.
%
%
%
%
\end{proof} 

\begin{cor}\label{cor:Gplus:comp_gen}
Let $T$ be a countable tree and let $G$ be a closed subgroup of $\Aut(T)$, such that arc stabilizers in $G$ are compact and such that $G$ leaves no proper subtree invariant.  Suppose that $G^+$ is compactly generated and of general type, and let $G(v)$ be the local action of $G$ at $v \in VT$.  Then there is at least one $v \in VT$ such that $G(v)$ is transitive and generated by point stabilizers.
\end{cor}

\begin{proof}
We see that the smallest $G^+$-invariant subtree is $G$-invariant and hence equal to $T$; note also that $G^+$ has property $\propP{}$ with respect to edges, since it contains all the arc stabilizers of $G$, and hence has property~$\propP{}$ in general by Theorem~\ref{propertyP_oneclosure}.  Since $G^+$ is generated by elements that fix a vertex, the quotient $T''$ of the action of $G^+$ on $T$ is a tree (see \cite[I.5.4, Exercise~2]{Serre:trees}).  By Theorem~\ref{thm:comp_gen+geom_dense}, $G^+$ has finitely many orbits on $T$, so $T''$ is a finite tree; in particular, $T''$ has at least one leaf.  Given $v \in VT$ such that $G^+v$ is a leaf of $T''$, then $(G^+)_v$ clearly acts transitively on $o\inv(v)$, so $G_v$ also acts transitively on $o\inv(v)$.  By Proposition~\ref{prop:G+_vertex}, the local action of $G^+$ at $v$ is the subgroup $G(v)^+$ generated by point stabilizers, so $G(v)^+$ is transitive.  Since $G(v)^+$ contains the point stabilizers, it follows that in fact $G(v)^+ = G(v)$.
\end{proof}

We finish this appendix with a sufficient condition to ensure all proper open subgroups of a group acting on a tree fix a vertex.

\begin{thm}\label{thm:open_primitive}
Let $T$ be a tree such that every vertex of $T$ has at least three neighbours, and let $G$ be a closed locally compact subgroup of $\Aut(T)$ with property $\propP{}$.  Let $V_+T$ and $V_-T$ be the two parts of the natural bipartition of $VT$, and suppose that $gV_+T = V_+T$ for all $g \in G$.  Suppose also that for every $v \in VT$, the local action $G(v)$ is primitive but not regular.  Then $G$ is a nondiscrete simple group taking the form $G = G_v \ast_{G_{(v,w)}} G_w$, where $v$ and $w$ are any pair of adjacent vertices, and $G$ has primitive action on $V_+T$ and $V_-T$.  In addition, every proper open subgroup of $G$ is contained in a conjugate of $G_v$ or $G_w$.
\end{thm}

\begin{proof}
Since $G$ has transitive local action at every vertex, we see that $G$ acts transitively on $V_+T$ and $V_-T$, so we have a Bass--Serre decomposition of $G$ as $G = G_v \ast_{G_{(v,w)}} G_w$ where $v \in V_+T$ and $w \in V_-T$ are any pair of adjacent vertices.  In particular, the action of $G$ is geometrically dense, so arc stabilizers are compact by Theorem~\ref{thm:comp_gen+geom_dense}.  Since $G$ also has property $\propP{}$, it follows that $G$ is the universal group $\Univ(G(v),G(w))$ defined in \cite{SmithDuke}.  By \cite[Theorem~26]{SmithDuke}, $G$ acts primitively on both parts of the natural bipartition; by \cite[Theorem~1]{SmithDuke}, $G$ is simple and nondiscrete.

We next note that for all ends $\xi$ of the tree, the stabilizer $G_{\xi}$ is not open.  We see this by considering a ray $(v_0,v_1,\dots)$ of vertices representing $\xi$: let $H_n = \bigcap^n_{i=0}G_{v_i}$.  By property~$\propP{}$, the local action of $H_n$ at $v_n$ is a point stabilizer of $G(v_n)$, and since $G(v_n)$ is primitive but not regular, no two point stabilizers are equal; thus the only neighbour of $v_n$ fixed by $H_n$ is $v_{n-1}$.  In particular, given $n \ge 1$, then $H_n$ does not fix $v_{n+1}$, that is, $H_{n+1}$ is properly contained in $H_n$.  Since $H_1$ is compact, it follows that the intersection $\bigcap_{n \ge 0}H_n$ is not open.  Since $G_{\xi} \cap G_{v_0} = G_{\xi} \cap \bigcap_{n \ge 0}H_n$ and $G_{v_0}$ is open, we deduce that $G_{\xi}$ is not open.

Now consider an open subgroup $H$ of $G$ that does not fix a vertex.  By the previous paragraph, $H$ does not fix any end, and hence by \cite[Corollaire~3.5]{Tits70}, there is a unique smallest $H$-invariant subtree $T'$ of $H$, which is infinite.  By \cite[Lemma~2.1(iii)]{MollerVonk}, $T'$ is the union of the axes of translation of $H$.  Given an axis $L$, we can argue as in the proof of Proposition~\ref{prop:tree_RIO} that $H$ contains $G_a$ for all $a \in AL$.  In particular, given a vertex $v' \in VL$, then there are two arcs $a,b \in AL$ incident with $v'$, and we have $\grp{G_a,G_b} \le H$.  Since the local action is primitive but not regular, we see that in fact $\grp{G_a,G_b} = G_{v'}$.  In particular, for each $v' \in VT'$ then $H$ has transitive local action at $v'$, ensuring that $T'$ contains all neighbours of $v'$.  Thus $T' = T$ and in fact $H$ contains every vertex stabilizer of $G$, ensuring that $H = G$.  In other words, every proper open subgroup of $G$ fixes a vertex; since $G$ acts transitively on $V_+T$ and $V_-T$, every vertex stabilizer is conjugate to one of $v$ and $w$.
\end{proof}

\end{document}